\documentclass[aop]{imsart}

\RequirePackage{amsthm,amsmath,amsfonts,amssymb}
\RequirePackage[numbers]{natbib}
\RequirePackage[colorlinks,citecolor=blue,urlcolor=blue]{hyperref}%% uncomment this for coloring bibliography citations and linked URLs

\startlocaldefs
\numberwithin{equation}{section}
\theoremstyle{plain}
\newtheorem{theorem}{Theorem}[section]
\newtheorem{lemma}[theorem]{Lemma}
\newtheorem{proposition}[theorem]{Proposition}
\newtheorem{corollary}[theorem]{Corollary}
\theoremstyle{definition}
\newtheorem{definition}[theorem]{Definition}
\newtheorem{remark}[theorem]{Remark}
\newtheorem{example}[theorem]{Example}
\usepackage{mathrsfs,mathtools,dirtytalk}

\DeclareMathOperator{\spn}{span}
\DeclareMathOperator{\ospn}{\overline{span}}
\DeclareMathOperator*{\loc}{loc}
\DeclareMathOperator*{\esssup}{ess\,sup}
\DeclareMathOperator*{\stk}{stk}
\DeclareMathOperator*{\rank}{rank}
\DeclareMathOperator*{\tr}{tr}
\DeclareMathOperator*{\Var}{Var}
\DeclareMathOperator*{\ran}{ran}
\DeclareMathOperator*{\Unif}{Unif}
\DeclareMathOperator*{\TV}{TV}
\DeclareMathOperator*{\sw}{\mathsf{sw}}
\newcommand{\sgn}{\normalfont{\text{sgn}}}

\newcommand{\pred}{\normalfont{\text{pred}}}

\newcommand{\Lip}{\normalfont{\text{Lip}}}
\newcommand{\Law}{\normalfont{\text{Law}}}
\newcommand{\lex}{\normalfont{\text{lex}}}
\newcommand{\tgamma}{\widetilde{\gamma}}
\newcommand{\cL}{\mathcal{L}}
\newcommand{\cE}{\mathcal{E}}
\newcommand{\hB}{\widehat{B}}
\newcommand{\hQ}{\widehat{Q}}
\newcommand{\heta}{\widehat{\eta}}
\newcommand{\hN}{\widehat{N}}
\newcommand{\hP}{\widehat{P}}
\newcommand{\homega}{\widehat{\omega}}
\newcommand{\hz}{\widehat{z}}
\newcommand{\hOmega}{\widehat{\Omega}}
\newcommand{\hbG}{\widehat{\mathbb{G}}}
\newcommand{\hbH}{\widehat{\mathbb{H}}}
\newcommand{\hsH}{\widehat{\mathscr{H}}}
\newcommand{\hE}{\widehat{E}}
\newcommand{\dint}{\,\normalfont{\text{d}}}
\newcommand{\di}{\normalfont{\text{d}}}
\newcommand{\mS}{\mathsf{S}}
\newcommand{\mB}{\mathsf{B}}
\newcommand{\hmB}{\widehat{\mathsf{B}}}
\newcommand{\tmB}{\widetilde{\mathsf{B}}}
\newcommand{\mM}{\mathsf{M}}
\newcommand{\mH}{\mathsf{H}}
\newcommand{\mA}{\mathsf{A}}
\newcommand{\mU}{\mathsf{U}}
\newcommand{\mP}{\mathsf{P}}
\newcommand{\ma}{\mathsf{a}}
\newcommand{\me}{\mathsf{e}}
\newcommand{\cS}{\mathcal{S}}

\newcommand{\cF}{\mathcal{F}}
\newcommand{\cQ}{\mathcal{Q}}
\newcommand{\cU}{\mathcal{U}}
\newcommand{\cP}{\mathcal{P}}
\newcommand{\cC}{\mathcal{C}}
\newcommand{\cW}{\mathcal{W}}
\newcommand{\cB}{\mathcal{B}}
\newcommand{\cY}{\mathcal{Y}}
\newcommand{\cR}{\mathcal{R}}
\newcommand{\cI}{\mathcal{I}}

\newcommand{\bF}{\mathbb{F}}
\newcommand{\bG}{\mathbb{G}}
\newcommand{\bH}{\mathbb{H}}
\newcommand{\bR}{\mathbb{R}}
\newcommand{\bN}{\mathbb{N}}
\newcommand{\bS}{\mathbb{S}}
\newcommand{\bRd}{\mathbb{R}^d}
\newcommand{\bRm}{\mathbb{R}^m}
\newcommand{\tb}{\widetilde{b}}
\newcommand{\tsigma}{\widetilde{\sigma}}
\newcommand{\sK}{\mathscr{K}}
\newcommand{\sD}{\mathscr{D}}
\newcommand{\cA}{\mathcal{A}}
\newcommand{\cD}{\mathcal{D}}
\newcommand{\cK}{\mathcal{K}}
\newcommand{\bK}{\mathbb{K}}
\newcommand{\bQ}{\mathbb{Q}}
\newcommand{\fP}{\mathfrak{P}}
\newcommand{\fp}{\mathfrak{p}}
\newcommand{\fC}{\mathfrak{C}}
\newcommand{\fr}{\mathfrak{r}}
\newcommand{\ofr}{\overline{\mathfrak{r}}}
\newcommand{\fM}{\mathfrak{M}}
\newcommand{\fm}{\mathfrak{m}}
\newcommand{\sF}{\mathscr{F}}
\newcommand{\hsF}{\widehat{\mathscr{F}}}
\newcommand{\sB}{\mathscr{B}}
\newcommand{\sG}{\mathscr{G}}
\newcommand{\hsG}{\widehat{\mathscr{G}}}
\newcommand{\sP}{\mathscr{P}}
\newcommand{\sH}{\mathscr{H}}
\def\giv{\,|\,} 
\newcommand{\cH}{\mathcal{H}}
\newcommand{\cM}{\mathcal{M}}
\newcommand{\BL}{\normalfont{\text{BL}}}

\newcommand{\step}{\normalfont{\text{step}}}
\newcommand{\hvarphi}{\widehat{\varphi}}
\newcommand{\HS}{\normalfont{\text{HS}}}
\newcommand{\id}{\normalfont{\text{id}}}
\newcommand{\ob}{\overline{b}}
\newcommand{\osigma}{\overline{\sigma}}
\newcommand{\hX}{\widehat{X}}
\newcommand{\hpi}{\widehat{\pi}}
\newcommand{\hnu}{\widehat{\nu}}
\newcommand{\hbeta}{\widehat{\beta}}
\newcommand{\tB}{\widetilde{B}}
\newcommand{\oOmega}{\overline{\Omega}}
\newcommand{\oP}{\overline{P}}
\newcommand{\pr}{\mathsf{\text{pr}}}
\newcommand{\st}{\normalfont{\text{st}}}
\newcommand{\col}{\mathsf{col}}
\newcommand{\hlambda}{\widehat{\lambda}}
\newcommand{\hSigma}{\widehat{\Sigma}}
\newcommand{\hL}{\widehat{L}}
\usepackage{xspace}
\newcommand{\cadlag}{c\`adl\`ag\xspace}

\makeatletter
\newcommand*\bigcdot{\mathpalette\bigcdot@{.5}}
\newcommand*\bigcdot@[2]{\mathbin{\vcenter{\hbox{\scalebox{#2}{$\m@th#1\bullet$}}}}}
\makeatother

\endlocaldefs

\begin{document}

\begin{frontmatter}
\title{Strong relaxed dynamics in stochastic control}
\runtitle{Strong relaxed dynamics in stochastic control}
\runauthor{B.~Kotlov}

\begin{aug}
%%%%%%%%%%%%%%%%%%%%%%%%%%%%%%%%%%%%%%%%%%%%%%%
%% Only one address is permitted per author. %%
%% Only division, organization and e-mail is %%
%% included in the address.                  %%
%% Additional information such as            %%
%% identifying the corresponding author must %%
%% be included in in the Acknowledgments     %%
%% section if necessary.                     %%
\author[A]{\fnms{Benjamin}~\snm{Kotlov}\ead[label=e1]{benjamin.kotlov@math.ethz.ch}}
%\author[B]{\fnms{???}~\snm{???}\ead[label=e2]{???@???}}
%\author[B]{\fnms{???}~\snm{???}\ead[label=e3]{???@???}}
%%%%%%%%%%%%%%%%%%%%%%%%%%%%%%%%%%%%%%%%%%%%%%
%% Addresses                                %%
%%%%%%%%%%%%%%%%%%%%%%%%%%%%%%%%%%%%%%%%%%%%%%
\address[A]{Department of Mathematics, ETH Zurich\printead[presep={,\ }]{e1}}

%\address[B]{???\printead[presep={,\ }]{e2,e3}}
\end{aug}

\begin{abstract}
    We characterize relaxed dynamics for which any finite family of strong state solutions, together with their relaxed inputs, can be approximated by the prescribed controlled dynamics, ground-stably in law. 
    
    Relaxed inputs include parametrizations of relaxed controls over the unit interval. Unlike relaxed controls, these parametrizations determine action couplings. The proposed relaxed dynamics is driven by a single cylindrical Brownian motion in a Hilbert space. For each relaxed input, its state is the unique strong solution of the corresponding equation. The collection of these states, over all relaxed inputs, satisfies the approximation requirement.

    Conversely, suppose one assigns to each relaxed input a continuous adapted state process on a fixed probability space. Pairwise approximation forces the assignment to be the strong solution map of the same dynamics, for a unique cylindrical Brownian driver. In fact, even for a state map indexed by (non-parametrized) relaxed controls, the relevant relaxed control class admits coherent parametrizations under separability and finite-family approximation assumptions, and the same type of dynamics is recovered again.
\end{abstract}



\end{frontmatter}
%%%%%%%%%%%%%%%%%%%%%%%%%%%%%%%%%%%%%%%%%%%%%%
%% Please use \tableofcontents for articles %%
%% with 50 pages and more                   %%
%%%%%%%%%%%%%%%%%%%%%%%%%%%%%%%%%%%%%%%%%%%%%%
\tableofcontents

%%%%%%%%%%%%%%%%%%%%%%%%%%%%%%%%%%%%%%%%%%%%%%
%%%% Main text entry area:

\section{Introduction}
 
Fix a stochastic basis $(\Omega, \sG_T,\bG = (\sG_t)_{0 \leqslant t \leqslant T},P)$ over a finite time horizon $T > 0$ carrying a $d$-dimensional Brownian motion $W$. Consider the controlled It\^o equation, for $m$-dimensional~$X^u$,
\begin{equation}
    \label{intro:strict_SDE}
    X^u_t = x_0 + \int_0^t b(s,X^u_s, u_s) \dint s + \int_0^t \sigma(s,X^u_s,u_s) \dint W_s,
\end{equation}
where $u$, a \textit{strict control}, is a $\bG$-predictable process taking values in a Polish action space~$\mA$ with at least two points. The \textit{system} $\mS \coloneq (m,x_0,b,\sigma)$ is implicit;
write $X^{u;\mS} \coloneq X^u$ when ambiguity arises. Maintaining standard assumptions (Section~\ref{sec:strict_controls_and_systems}) on the \textit{input} $(u,\mS)$, the corresponding \textit{state} $X^u$ is the unique strong solution of~(\ref{intro:strict_SDE}). Varying $(u,\mS)$ in (\ref{intro:strict_SDE}) yields a family of $X^{u;\mS}$'s, defined simultaneously on the fixed basis and driven by the common $W$. We seek to extend the map $(u,\mS) \mapsto X^{u;\mS}$ to relaxed inputs on a fixed extension of the ground basis, in a \textit{consistent} way. Informally, this means that there is an extension of the ground basis on which (1) every relaxed input has a corresponding state, required to be a continuous adapted process of dimension specified by the input, and (2) every finite family of relaxed states, together with their relaxed inputs, can be jointly approximated by~(\ref{intro:strict_SDE}), $\sG_T$-stably in law. Our main findings are that such state assignments exist, and are precisely the strong solution maps of the dynamics~(\ref{intro:parametrized_SDE}) below.

Classical relaxation replaces $u_t(\omega)$ by a predictable probability kernel $\rho_t(\omega,\di a) \in \cP(\mA)$, with strict controls embedded through $\rho_t = \delta_{u_t}$ \cite{ElKarouiHuuNguyenJeanblancPicque1987,HaussmannLepeltier1990}. For two \textit{relaxed controls} $\rho,\rho'$, a jointly approximating sequence of strict-control pairs implicitly selects a coupling of $\rho_t,\rho_t'$ on~$\mA^2$, which may depend on the choice of the sequence and, in particular, is not determined by the marginals $\rho_t, \rho_t'$ alone. To carry couplings intrinsically, we parametrize by $r \in [0,1]$. Indeed, a \textit{parametrized relaxed control} $\pi$ predictably specifies a profile $\pi_t : [0,1] \rightarrow \mA$, which induces $\rho^\pi_t \coloneq (\pi_t)_\# \dint r$; for any finite family $\Pi = (\pi^i)_{i \in I}$, a coupling of~$(\rho^{\pi^i}_t)_{i \in I}$ on~$\mA^I$ is canonically given by $\lambda^\Pi_t \coloneq (\Pi_t)_\# \dint r$. The $\lambda^\Pi$ are consistent with subfamily projections.

Fix an extension $(\Omega,\sF_T,\bF,P)$ of the ground basis in which $\bG$ is immersed, carrying a cylindrical $\bF$-Brownian motion $\mB$ in $\mH \coloneq L^2([0,1];\bRd)$ such that, identifying $\bRd$ with the constant functions, $\mB|_{\bRd}$ coincides with $W$ and $\mB|_{(\bRd)^\perp}$ is independent of $\sG_T$. For an $\bF$-predictable $\pi$ that is sufficiently integrable (Section~\ref{sec:parametrized_relaxed_controls}), let $X^{\pi;\mS} \coloneq X^{\pi}$ be the unique $\bF$-strong solution~of
\begin{equation}
    \label{intro:parametrized_SDE}
    X^\pi_t = x_0 + \int_0^t \ob(s,X^\pi_s, \rho^\pi_s) \dint s + \int_0^t \Sigma(s,X^\pi_s, \pi_s) \dint \mB_s,
\end{equation}
with drift averaged and diffusion identified with its Hilbert--Schmidt operator in $\cL_2(\mH,\bR^m)$:
$$
    \ob(s,X^\pi_s,\rho^\pi_s) \coloneq \int_\mA  b(s,X^\pi_s,a) \, \rho^\pi_s(\di a), \qquad \Sigma(s,X^\pi_s, \pi_s) h \coloneq \int_0^1 \sigma\big(s,X^\pi_s, \pi_s(r)\big)h(r) \dint r.
$$
As $\mB|_{\bRd} = W$,
(\ref{intro:parametrized_SDE}) reduces to (\ref{intro:strict_SDE}) for $\pi^u_t \equiv u_t$. Varying $(\pi,\mS)$ in~(\ref{intro:parametrized_SDE}) yields an assignment $(\pi,\mS) \mapsto X^{\pi;\mS}$ on the extension, with each $X^{\pi;\mS}$ driven by the common~$\mB$. Of note, for a finite family $\Pi = (\pi^i)_{i \in I}$ and systems $\mS^i, \mS^{i'}$ with diffusion coefficients $\sigma_i,\sigma_{i'}$, (\ref{intro:parametrized_SDE}) gives
$$
    \di [X^{\pi^i;\mS^i},X^{\pi^{i'};\mS^{i'}}]_t = \bigg( \int_{\mA^I} \sigma_i(t,X^{\pi^i;\mS^i}_t,a_i) \sigma_{i'}(t,X^{\pi^{i'};\mS^{i'}}_t, a_{i'})^\top \, \lambda^\Pi_t(\di a) \bigg) \dint t.
$$
This demonstrates how mixed state characteristics depend on the joint-action couplings. 

\subsection{Main results}

For any finite families $U = (u^i)_{i \in I}$, $\Pi = (\pi^i)_{i \in I}$, let
$$
    \eta^U(\di t, \di a) \coloneq \di t\, \delta_{U_t}(\di a), \qquad \eta^\Pi(\di t, \di a) \coloneq \di t\, \lambda^\Pi_t(\di a)
$$
denote their joint Young measures on $[0,T] \times \mA^I$, respectively.

\begin{longlist}
    \item\label{intro:(i)} \textit{For any finite family $\Pi = (\pi^i)_{i \in I}$ and $(\mS^i)_{i \in I}$, there are $U^n = (u^{n,i})_{i \in I}$ with
    \begin{equation}
        \label{intro:joint_approximation}
        \big(\eta^{U^n}, (X^{u^{n,i};\mS^i})_{i \in I}\big) \Longrightarrow_{\sG_T} \big(\eta^{\Pi}, (X^{\pi^i;\mS^i})_{i \in I}\big).
    \end{equation}
    Here Young measures carry a Wasserstein topology, states carry the uniform topology, and the limit holds $\sG_T$-stably in law} (Theorem~\ref{theorem:cylindrical_approximation}).
    \item\label{intro:(ii)} \textit{Let $(\pi,\mS) \mapsto Y^{\pi;\mS}$ be an assignment on an arbitrary filtered extension of the ground basis, meaning that $Y^{\pi;\mS}$ is a continuous adapted $\bR^{m_\mS}$-valued process on the extension for all $(\pi,\mS)$. Assume the analogue of~(\ref{intro:joint_approximation}) holds for all pairs $(\pi,\mS), (\pi',\mS')$. Then $\bG$ is immersed. On the extension, there is a unique cylindrical Brownian motion $\mB$ in $\mH$, with $\mB|_{\bRd} = W$, $\mB|_{(\bRd)^\perp}$ independent of $\sG_T$, and
    \begin{equation}
        \label{intro:parametrized_SDE_reconstructed}
        Y^\pi_t = x_0 + \int_0^t \ob(s,Y^\pi_s, \rho^\pi_s) \dint s + \int_0^t \Sigma(s,Y^\pi_s, \pi_s) \dint \mB_s
    \end{equation}
    for any $\pi$ and $\mS$}
    (Theorem~\ref{theorem:cylindrical_converse}). \textit{Thus $(\pi,\mS) \mapsto Y^{\pi;\mS}$ is the strong solution map of (\ref{intro:parametrized_SDE_reconstructed}).}
    \item\label{intro:(iii)} \textit{Let $(\rho,\mS) \mapsto Y^{\rho;\mS}$ be an assignment on an arbitrary extension (Definition~\ref{def:relaxed_state_assignment}) with the pairwise marginal approximation property (Definition~\ref{def:pair_approximation}). Consider a class~$\cC$ of relaxed controls satisfying a certain separability assumption (Definition~\ref{def:countable_residual_multiplicity}) and the finite-family marginal approximation property (Definition~\ref{def:finite_family_approximation}). On the extension, there exist  $(\pi^\rho)_{\rho \in \cC}$ parametrizing~$\cC$, and a cylindrical martingale $\mM$ in $\mH$, with $\mM|_{\bRd} = W$, $[\mM(h),Z] \equiv 0$ for all $h \in {(\bRd)^\perp}$ and continuous local $\bG$-martingales $Z$, and
    \begin{equation}
        \label{intro:parametrized_SDE_reconstructed_ordinary}
        Y^\rho_t = x_0 + \int_0^t \ob(s,Y^\rho_s, \rho_s) \dint s + \int_0^t \Sigma(s,Y^\rho_s, \pi^\rho_s) \dint \mM_s
    \end{equation}
    for all $\rho \in \cC$ and $\mS$}
    (Theorem~\ref{theorem:ordinary_relaxed_assignment_common_label_dynamics}). \textit{Moreover, $\mM$ can be completed to a cylindrical Brownian motion in $\mH$ on a further Wiener extension while preserving~(\ref{intro:parametrized_SDE_reconstructed_ordinary})} (Corollary~\ref{cor:complete_M_to_B}).
\end{longlist}

Assertions (\ref{intro:(i)}) and (\ref{intro:(ii)}) constitute the claimed characterization of strong relaxed dynamics, and (\ref{intro:(iii)}) provides justification for parametrizing ordinary relaxed controls in the first place.

The reconstructions (\ref{intro:(ii)})--(\ref{intro:(iii)}) also hold if the approximating strict controls are predictable for enlarged filtrations in which $\bG$ is immersed. In (\ref{intro:(i)}), they may be chosen to be randomized sample-and-hold executions (Remark~\ref{remark:sample_and_hold_cylindrical}).

\subsection{Key observation}

On an arbitrary extension of the ground basis, let $\rho$ be a relaxed control and $Y$ a continuous $\bRm$-valued process.  Set $\eta^\rho(\di t,\di a) \coloneq \di t \, \rho_t(\di a)$ on $[0,T] \times \mA$, and suppose that for some $u^n$,
\begin{equation}
    \label{intro:single_approximation}
    (\eta^{u^n}, X^{u^n}) \Longrightarrow_{\sG_T} (\eta^\rho, Y).
\end{equation}
Then $\bG$ is immersed in the limiting filtration $\bF^{\eta^\rho,Y}$ of Proposition~\ref{prop:martingale_characteristics}.
Writing $\osigma$ for the averaged diffusion coefficient like $\ob$, set
$$
    L_t \coloneq Y_t - x_0 - \int_0^t \ob(s,Y_s,\rho_s) \dint s, \qquad N_t \coloneq L_t - \int_0^t \osigma(s,Y_s,\rho_s) \dint W_s.
$$
Computing from~(\ref{intro:strict_SDE}) the drift of $X^{u^n}$ and the brackets $[X^{u^n}]$ and $[X^{u^n},Z]$ for any local $\bG$-martingale $Z$ and passing to the limit shows $L$ and $N$ are local $\bF^{\eta^\rho,Y}$-martingales and
$$
    \di [L]_t = \bigg( \int_\mA \sigma(t,Y_t,a) \sigma(t,Y_t,a)^\top \, \rho_t(\di a) \bigg) \dint t, \qquad \di [L,Z]_t = \osigma(t,Y_t,\rho_t) \dint [W,Z]_t
$$
(Proposition~\ref{prop:martingale_characteristics}). Thus (\ref{intro:single_approximation}) determines the diffusion covariance of $Y$ and its covariation with all local $\bG$-martingales. Moreover, $\di [N]_t = R(t,Y_t,\rho_t) \dint t$ and $[N,Z] \equiv 0$, where
$$
    R(t,Y_t,\rho_t) \coloneq \int_\mA \sigma(t,Y_t,a) \sigma(t,Y_t,a)^\top \, \rho_t(\di a) - \osigma(t,Y_t,\rho_t)\osigma(t,Y_t,\rho_t)^\top.
$$
A Brownian representation of $N$ shows that, possibly on a further Wiener extension carrying an $m$-dimensional Brownian motion $B$, independent of $\sG_T$,
\begin{equation}
    \label{intro:residual_SDE}
    Y_t = x_0 + \int_0^t \ob(s,Y_s,\rho_s) \dint s + \int_0^t \osigma(s,Y_s,\rho_s) \dint W_s + \int_0^t R(s,Y_s,\rho_s)^{1/2} \dint B_s,
\end{equation}
where $R^{1/2}$ is the principal square root (Theorem~\ref{theorem:residual_representation}). The extension and driver~$B$ in~(\ref{intro:residual_SDE}) may be particular to $(\eta^\rho,Y)$. Conversely, given an equation of the form~(\ref{intro:residual_SDE}) on an extension in which $\bG$ is immersed, its unique solution $Y$ satisfies~(\ref{intro:single_approximation}) for some $u^n$ (Theorem~\ref{theorem:residual_approximation}). 

There is a correspondence between state solutions of (\ref{intro:parametrized_SDE}) and (\ref{intro:residual_SDE}) up to further Wiener extensions  (Proposition~\ref{prop:residual_cylindrical_correspondence}). A finite family $\Pi = (\pi^i)_{i \in I}$ is a single $\mA^I$-valued parametrized relaxed control, inducing the $\mA^I$-relaxed control $\lambda^\Pi$. Systems $(\mS^i)_{i \in I}$ can be stacked as one system $\mS^{\stk}$ on $\mA^I$. Denote by $X^{\Pi;\mS^{\stk}}$ the solution to the analogue of~(\ref{intro:parametrized_SDE}) over $\mA^I$ for the input $(\Pi,\mS^{\stk})$, which is equivalently the vector of $(X^{\pi^i;\mS^i})_{i \in I}$. The preceding correspondence expresses $X^{\Pi;\mS^{\stk}}$ in the form~(\ref{intro:residual_SDE}) for $(\lambda^\Pi,\mS^{\stk})$, and Theorem~\ref{theorem:residual_approximation} yields claim~(\ref{intro:(i)}) above.

\subsection{Relation to existing work}

Classical relaxation of controlled diffusions is formulated with controlled martingale problems~\cite{ElKarouiHuuNguyenJeanblancPicque1987, HaussmannLepeltier1990}, and also through stochastic equations driven by action-space martingale measures~\cite{ElKarouiMeleard1990}. In the corresponding weak control formulations, the stochastic basis and driving noise form part of an admissible control-state realization, rather than being prescribed uniformly across all inputs; see also~\cite{ElKarouiTan2024PartII}. These formulations describe control-state laws, whereas we are concerned with a simultaneous strong solution assignment on one fixed extension of the ground basis.

Strong relaxed dynamics and couplings are also present in the literature. El Karoui--M\'el\'eard~\cite{ElKarouiMeleard1990} establish white-noise representations and pathwise uniqueness for equations driven by a prescribed martingale measure. They also explain that the random intensity of the driving martingale measure need not determine the resulting state law. Fontbona, Gu\'erin, and M\'el\'eard~\cite{FontbonaGuerinMeleard2010} construct couplings of martingale measures through predictable optimal transport maps between covariance measures. Our objective can be equivalently expressed as selecting particular martingale measures, simultaneously for all relaxed inputs on a fixed extension, consistently with joint strict approximation (Proposition~\ref{prop:relaxed_state_assignment_MM}, Remark~\ref{remark:parametrized_to_MM}).

Relaxed controls have acquired renewed relevance in continuous-time reinforcement learning, where distribution-valued policies describe exploration in an unknown environment through random action selection \cite{JMLR:v21:19-144}. The relationship between distribution-valued policies, their implementation through sampled actions, and the resulting state dynamics is a natural question. Bender and Thuan~\cite{BenderThuan2026} develop grid-sampling limit dynamics for randomized feedback policies in models with diffusion and jumps.
We characterize a common strong assignment for general predictable open-loop controls through joint ground-stable approximation. Our approximating strict controls are not required to be $\bG$-predictable; randomized sample-and-hold implementations in the spirit of~\cite{BenderThuan2026} are permissible (Remark~\ref{remark:sample_and_hold_cylindrical}).

\subsection{Structure of paper}

Section~\ref{sec:2} introduces the control spaces and identifies martingale characteristics of ground-stable strict limits (Proposition~\ref{prop:martingale_characteristics}). These are used in Section~\ref{sec:3} to establish both approximation directions of~(\ref{intro:residual_SDE}) for a single relaxed control (Theorems~\ref{theorem:residual_representation} and~\ref{theorem:residual_approximation}). Section~\ref{sec:4} proves finite-family approximation of~(\ref{intro:parametrized_SDE}) (Theorem~\ref{theorem:cylindrical_approximation}), and conversely, shows that a parametrized relaxed state assignment with the pairwise approximation property is the strong solution map of dynamics of the same form (Theorem~\ref{theorem:cylindrical_converse}). Section~\ref{sec:5} recovers pairwise action couplings for a state assignment indexed by ordinary (non-parametrized) relaxed controls (Proposition~\ref{prop:gamma^rho,rho'}), obtains parametrizations (Proposition~\ref{prop:common_label_lift}), and reconstructs dynamics of the form~(\ref{intro:parametrized_SDE}) under its additional assumptions (Theorem~\ref{theorem:ordinary_relaxed_assignment_common_label_dynamics}). Section~\ref{sec:6} identifies the associated martingale measures and coherent choices made by the common driver (Proposition~\ref{prop:relaxed_state_assignment_MM} and Remark~\ref{remark:parametrized_to_MM}). Section~\ref{sec:7} concludes.

\subsection{Notation and conventions}

All filtrations are indexed by $[0,T]$ and satisfy the usual conditions. A filtration specified through generators means the usual augmentation, under the ambient probability measure, of the corresponding raw filtration. For a process $X$ and ambient filtration $\bF$, write $\bF^X$ for its natural filtration. Write $\cM_{\loc,0}(\bF;\bRm)$, $\cM^c_{\loc,0}(\bF;\bRm)$, and $\cM^{2,c}_0(\bF;\bRm)$ for the spaces of \cadlag local, continuous local, and continuous square-integrable $m$-dimensional $\bF$-martingales null at zero, respectively. When $m = 1$, omit $\bR$. 
We frequently work with filtered extensions, and the enlarged probability space and probability measure are typically not relabeled, with ground variables and filtrations pulled back without changing notation. Write $\bH \subset \bK$ if $\sH_t \subset \sK_t$ for all $t \in [0,T]$, and $\bH \hookrightarrow \bK$ if moreover $\bH$ is immersed in $\bK$, which is to say that every (bounded) \cadlag $\bH$-martingale remains a $\bK$-martingale. An $\bH$-Brownian motion remains $\bK$-Brownian in this case. 

The standard basis of $\bRd$ is denoted $\epsilon_1,\dots,\epsilon_d$. Write $|\cdot|$ for the Euclidean norm of vectors and the Frobenius norm of matrices. Let $\col(x,y)\in \bR^{m+m'}$ be vertical concatenation of vectors $x \in \bRm$ and $y \in \bR^{m'}$, and write $\col_{i \in I}(x_i)$ for the analogous vertical concatenation for finite indexing sets $I$.
Index sets are always nonempty.
Let $\bS^m_+$ be the positive semidefinite symmetric $m \times m$ matrices, and $R^{1/2}$ the principal square root of $R \in \bS^m_+$. For an $m \times \ell$ matrix $D$, let $D^\dagger$ be its Moore--Penrose inverse; $D \mapsto D^\dagger$ is a Borel function on $\bR^{m \times \ell}$.

For separable real Hilbert spaces $H,H'$, let $\cL(H,H')$ and $\cL_2(H,H')$ denote the bounded and Hilbert--Schmidt operators, respectively. Abbreviate $\cL(H) \coloneq \cL(H,H)$. When the underlying Hilbert space is clear from the context, we write its inner product and norm as $\langle \cdot, \cdot \rangle$ and $\|\cdot\|$. For $T \in \cL(H,H')$, write $T^*$, $\ker T$, and $\ran T$ for its adjoint, kernel, and range, and write $P_M$ for the orthogonal projection onto a closed subspace $M \subset H$. An $\cL(H,H')$-valued process is \textit{predictable} if it is predictable for the strong operator topology.

For nonnegative quantities $f,g$, write $f \lesssim g$ if $f \leqslant C g$ for a finite deterministic constant~$C$, which may change from line to line.

\section{Controls, parametrizations, and strict-state limits}
\label{sec:2}

We introduce the control spaces and identify martingale characteristics of ground-stable limits of strict control-state pairs.

\subsection{Ground setting}
\label{sec:strict_controls_and_systems}

Fix a ground stochastic basis $(\Omega, \sG_T,\bG = (\sG_t)_{0 \leqslant t \leqslant T},P)$ over a finite time horizon $T > 0$ carrying a $d$-dimensional Brownian motion $W$. Put $\oOmega \coloneq \Omega \times [0,T]$, $\oP \coloneq P \otimes \di t$, and for any filtration $\bH$, denote by $\sP^{\bH}$ its predictable $\sigma$-algebra on $\oOmega$. Take the space of \textit{actions} to be a pointed Polish space $(\mA,\ma_0,d)$ containing at least two points.

Throughout the section, $\bK$ denotes a filtration, possibly defined on an extension of the ground space that is not relabeled, for which $\bG \hookrightarrow \bK$. 

\begin{definition}
    A $\bK$-predictable $u : \oOmega \rightarrow \mA$ is a \textit{$\bK$-strict control}.
\end{definition}

When no filtration is specified, strict controls are understood to be on the ground space and $\bG$-predictable. The larger $\bK$ supports randomized action selection; see Example~\ref{example:sample_and_hold_labeling_scheme}.

Fix $p \in [1,\infty)$ and $q \in [2,\infty)$ throughout. For $\bK$-strict controls $u,u'$, define
$$
    d_{p,q}(u,u') \coloneq E \bigg[\bigg(\int_0^T d(u_s,u'_s)^q \dint s\bigg)^{p/q}\bigg]^{1/p},
$$
and the corresponding space
$$
    \cU^{p,q}(\bK) \coloneq L^{p,q}(\oOmega,\sP^{\bK},\oP;\mA)
$$
of equivalence classes, up to equality $\oP$-a.e., of $\bK$-strict controls $u$ with $d_{p,q}(u,\ma_0) < \infty$. We require $q \geqslant 2$ to bound the $L^2$-terms entering the diffusion estimates. Since $\mA$ is Polish, $(\cU^{p,q}(\bK), d_{p,q})$ is a complete metric space. If $p = q$, the usual space $L^p(\oOmega,\sP^{\bK},\oP;\mA)$ is recovered. Whenever $p'\leqslant p$ and $q' \leqslant q$, $\cU^{p,q}(\bK) \subset \cU^{p',q'}(\bK)$; write $\cU(\bK) \coloneq \cU^{1,2}(\bK)$ for the total space. For $\bK = \bG$, we simply put $\cU^{p,q} \coloneq \cU^{p,q}(\bG)$ and $\cU \coloneq \cU(\bG)$.

The following notion specifies the state dimension, initial condition, and drift and diffusion coefficients of the It\^o equation determining the state corresponding to a strict control.

\begin{definition}
    \label{def:system}
    A \textit{system} $\mS = (m_{\mS}, x_0^{\mS}, b_{\mS}, \sigma_{\mS})$ consists of $m_{\mS} \in \bN$, $x_0^{\mS} \in L^1(\sG_0;\bR^{m_\mS})$, and jointly measurable $b_{\mS} : [0,T] \times \bR^{m_{\mS}} \times \mA \rightarrow \bR^{m_{\mS}}$, $\sigma_{\mS} : [0,T] \times \bR^{m_{\mS}} \times \mA \rightarrow \bR^{m_{\mS}\times d}$ with
    \begin{longlist}
        \item $|b_{\mS}(t,x,a)| + |\sigma_{\mS}(t,x,a)| \lesssim 1 + |x| + d(a,\ma_0)$,
        \item $|b_{\mS}(t,x,a) - b_{\mS}(t,x',a')| + |\sigma_{\mS}(t,x,a) - \sigma_{\mS}(t,x',a')| \lesssim |x-x'| + d(a,a')$,
    \end{longlist}
    for all $t \in [0,T]$, $x,x' \in \bR^{m_{\mS}}$, $a,a' \in \mA$. Let $\cS$ denote the class of all systems.
\end{definition}

For $(u,\mS) \in \cU(\bK) \times \cS$, consider the strict-state equation
\begin{equation}
    \label{eq:SDE_strict}
    X^{u;\mS}_t = x_0^{\mS} + \int_0^t b_{\mS}(s,X^{u;\mS}_s, u_s) \dint s +  \int_0^t \sigma_{\mS}(s,X^{u;\mS}_s, u_s) \dint W_s.
\end{equation}
When the underlying system is fixed, we usually abbreviate $X^u \coloneq X^{u;\mS}$. Throughout, a \textit{$\bK$-strong solution} of an SDE means a continuous $\bK$-adapted process, up to indistinguishability, for which all Lebesgue and stochastic integrals are well defined and the equation holds.

\begin{proposition}
    \label{prop:strict_well_posed}
    For $(u,\mS) \in \cU(\bK) \times \cS$, (\ref{eq:SDE_strict}) admits a unique $\bK$-strong solution $X^{u;\mS}$; if $u \in \cU^{p,2}(\bK)$ and $x_0^\mS$ is $p$-integrable, $E[\|X^{u;\mS}\|_\infty^p] \lesssim 1 + E[|x_0^\mS|^p] + d_{p,2}(u,\ma_0)^p$. For $u,u' \in \cU^{p,2}(\bK)$, $E[\|X^{u;\mS} - X^{u';\mS} \|_\infty^p] \lesssim d_{p,2}(u,u')^p$.
\end{proposition}

\begin{proof}
    Strong well-posedness follows by the usual Picard iteration argument. The estimates are consequences of the Burkholder--Davis--Gundy and Gr\"onwall inequalities. 
\end{proof}

The $\bK$-strong solution $X^{u;\mS}$ of~(\ref{eq:SDE_strict}) is the \textit{state} corresponding to the \textit{input} $(u,\mS)$. 

For a finite index set $I$, $\mA^I$ is a pointed Polish space with reference $(\ma_0)_{i \in I}$ and metric $d_I(a, a') \coloneq (\sum_{i \in I} d(a_i, a_i')^2 )^{1/2}$. The \textit{stacked system} of $(\mS^i)_{i \in I} \subset \cS$ is the system
\begin{equation}
    \label{eq:system_stacked}
    \mS^{\stk} \coloneq (m_{\stk}, x_0^{\stk}, b_{\stk},\sigma_{\stk}),
\end{equation}
on $\mA^I$ given by $m_{\stk} \coloneq \sum_{i\in I} m_{\mS^i}$, $x_0^{\stk} \coloneq \col_{i \in I}x^{\mS^i}_0$, $b_{\stk}(t,x, a) \coloneq \col_{i \in I} b_{\mS^i}(t,x_i,a_i)$, and $\sigma_{\stk}(t,x, a) \coloneq \col_{i \in I} \sigma_{\mS^i}(t,x_i,a_i)$. Indeed, the coefficients satisfy Definition~\ref{def:system} for~$\mA^I$. Any $U = (u^i)_{i \in I} \in \cU^{p,q}(\bK)^I$ is identified as an element of $\cU^{p,q}(\bK;\mA^I) \coloneq L^{p,q}(\oOmega,\sP^{\bK},\oP; \mA^I)$; the $(p,q)$-integrability holds because $I$ is finite. Write $X^{U;\mS^{\stk}}$ for the unique strong solution of~(\ref{eq:SDE_strict}) for $(U,\mS^{\stk})$ over $\mA^I$. Pathwise uniqueness in Proposition~\ref{prop:strict_well_posed}, componentwise, gives
\begin{equation}
    \label{eq:system_stacked_state}
    X^{U;\mS^{\stk}} = \col_{i \in I} X^{u^i;\mS^i}.
\end{equation}
The \textit{direct sum} $\oplus_{i \in I} \mS^i \in \cS$ is the system on $\mA$ obtained by restricting $\mS^{\stk}$ to the diagonal $\{ (a)_{i\in I} \in \mA^I : a \in \mA\}$. For all $u \in \cU(\bK)$,
\begin{equation}
    \label{eq:system_direct_sum_state}
    X^{u;\oplus_{i\in I} \mS^i} = \col_{i \in I} X^{u;\mS^i}.
\end{equation}

\subsection{Relaxed controls and Young measures}

Fix a filtered extension $(\Omega, \sF_T, \bF, P)$ of the ground basis---we do not require it to extend $(\Omega,\sK_T,\bK,P)$. Let $\cP(\mA)$ be the pointed Polish space of Borel probability measures, carrying the narrow topology and the corresponding Borel $\sigma$-algebra, with reference~$\delta_{\ma_0}$. Write $\cP_q(\mA)$ for those $\mu \in \cP(\mA)$ with finite $q$-moments, which is complete under the Wasserstein $q$-distance $\cW_q$. Its Borel $\sigma$-algebra coincides with the trace $\sigma$-algebra induced by $\cP(\mA)$, which follows by applying the Lusin--Suslin theorem \cite[Theorem~6.8.6]{bogachev2007measure} to the continuous injection $\cP_q(\mA) \hookrightarrow \cP(\mA)$.

\begin{definition}
    An $\bF$-predictable $\rho : \oOmega \rightarrow \cP(\mA)$ is an \textit{$\bF$-relaxed control}.
\end{definition}

A map $\rho : \oOmega \rightarrow \cP(\mA)$ is $\bF$-predictable if and only if $(\omega,t) \mapsto \int_{\mA}f(a) \, \rho_t(\omega,\di a)$ is for every $f \in C_b(\mA)$; see~\cite[Lemma~A5.5]{kallenberg2021foundations}. Thus equivalently, an $\bF$-relaxed control is a probability kernel from $(\oOmega, \sP^\bF)$ to $\mA$. Strict controls $u \in \cU$ are identified as $\bG$- (hence $\bF$-) relaxed controls under the Dirac embedding $u_t \mapsto \delta_{u_t}$. Analogously to strict controls, write
$$
    \cQ^{p,q}(\bF) \coloneq L^{p,q}\big(\oOmega, \sP^{\bF}, \oP; \cP_q(\mA)\big),
$$
with metric still denoted $d_{p,q}$, and set $\cQ(\bF) \coloneq \cQ^{1,2}(\bF)$. When the ambient filtration is clear from the context, we usually omit $\bF$ from the notation and simply say \say{relaxed control.} 

For a pointed Polish space $(E, \me_0, d_E)$, let $\cM_+([0,T] \times E)$ be the finite positive Borel measures on $[0,T] \times E$ carrying the narrow topology, which is metrized by
$$
    d_{\BL}(\mu,\mu') \coloneq \sup\bigg\{ \bigg| \int_{[0,T] \times E} f \dint (\mu - \mu') \bigg| : \|f\|_\infty \leqslant 1, \, \Lip(f) \leqslant 1 \bigg\},
$$
where $\Lip(f)$ is the smallest Lipschitz constant of $f$~\cite[Theorem~8.3.2]{bogachev2007measure}. Let
$$
    \cY_q^T(E) \coloneq \bigg\{\mu \in \cM_+([0,T]\times E) : \mu(\di t, E) = \di t, \, \int_{[0,T] \times E} d_E(e,\me_0)^q \, \mu(\di t, \di e) < \infty\bigg\}
$$
consist of the $\mu$ with time marginal $\di t$ and finite $q$-moments. Every $\mu \in \cY_q^T(E)$ has mass~$T$. Equip $[0,T] \times E$ with $d_{[0,T] \times E}((t,e), (t',e')) \coloneq (|t-t'|^2 + d_E(e,e')^2)^{1/2}$, and $\cY_q^T(E)$ with the mass-$T$ Wasserstein $q$-distance $(\mu,\mu') \mapsto T^{1/q}\cW_q(\mu/T, \mu'/T)$, still denoted $\cW_q$. The time marginal $\di t$ is preserved under Wasserstein convergence, thus $(\cY^T_q(E), \cW_q)$ is Polish. 

For $u \in \cU^{p,q}(\bK)$ and $\rho \in \cQ^{p,q}$, let $\eta^u, \eta^\rho : \Omega \rightarrow \cY^T_q(\mA)$ be their (random) Young measures
$$
    \eta^u(\di t, \di a) \coloneq \dint t \, \delta_{u_t}(\di a), \qquad \eta^\rho(\di t, \di a) \coloneq \dint t \, \rho_t(\di a).
$$
To express the information generated by a random Young measure, we use the following process. For a pointed Polish space $(E,\me_0,d_E)$ and $\mu \in \cY^T_2(E)$, write 
$$
    \mu_{\leqslant t} \coloneq \mu|_{[0,t] \times E} \in \cM_+([0,T] \times E), \qquad \mu_{\leqslant} \coloneq (\mu_{\leqslant t})_{0 \leqslant t\leqslant T} \in C\big( [0,T], \cM_+([0,T] \times E) \big),
$$
where $\mu|_{[0,t] \times E}(A) \coloneq \mu(A \cap ([0,t] \times E))$ for all $A \in \sB([0,T] \times E)$.

\begin{lemma}
    \label{lemma:Young_history_continuity}
    The map $\cY^T_2(E) \ni \mu \mapsto \mu_{\leqslant} \in C([0,T]; \cM_+([0,T] \times E))$ is continuous; it is well defined because $d_{\BL}(\mu_{\leqslant t}, \mu_{\leqslant s}) = |t-s|$.
\end{lemma}

\begin{proof}
    For $s \leqslant t$, $d_{\BL}(\mu_{\leqslant t}, \mu_{\leqslant s}) = |t-s|$ as $\mu_{\leqslant t} - \mu_{\leqslant s} = \mu|_{(s,t] \times E} \in \cM_+([0,T] \times E)$ has mass $t-s$. If $\mu^n \rightarrow \mu$ in $\cY^T_2(E)$,  $\mu^n \rightarrow \mu$ narrowly, hence $\mu^n_{\leqslant t} \rightarrow \mu_{\leqslant t}$ narrowly for every~$t$ since $\mu(\{t\} \times E) = 0$; $\sup_{n \geqslant 1} d_{\BL}(\mu^n_{\leqslant t}, \mu^n_{\leqslant s}) = |t-s|$ upgrades to uniform convergence.
\end{proof}

\begin{lemma}
    \label{lemma:predictable_disintegration}
    Let $\eta : \Omega \rightarrow \cY_2^T(E)$ be a random Young measure and assume $\eta_{\leqslant}$ is adapted to some filtration $\bH$. There is a unique, up to $\oP$-a.e.~equality, $\bH$-predictable $\rho : \oOmega \rightarrow \cP(E)$ such that $\eta(\omega;\di t, \di e) = \di t \, \rho_t(\omega;\di e)$ a.s. Moreover, $\rho$ may be chosen $\cP_2(E)$-valued. 
\end{lemma}

\begin{proof}
    For every $A \in \sB(E)$, the process $t \mapsto \eta([0,t] \times A)$ is $\bH$-adapted and continuous, hence $\bH$-predictable. Thus $\eta$ is a predictable random measure and equals its compensator. The predictable disintegration in the proof of~\cite[Theorem~10.21]{kallenberg2021foundations}, with marginal $\di t$, gives $\eta(\di t, \di e) = \di t \, \rho_t(\di e)$ for an $\bH$-predictable $\rho$. Uniqueness of disintegration gives the required uniqueness. At last, $\eta$ taking values in $\cY^T_2(E)$ implies $\rho$ has finite second moments $\oP$-a.e. Replacing $\rho_t$ by $\delta_{\me_0}$ on the predictable exceptional set yields the claimed choice.
\end{proof}

For every $\rho \in \cQ^{p,q}(\bG)$, $\bG$-predictable simple approximation with finite-action support and classical chattering as in Warga~\cite[Theorem~2.4]{Warga1967} give $u^n \in \cU^{p,q}$ with $\cW_q(\eta^{u^n}, \eta^\rho) \rightarrow 0$ in $L^p(\Omega)$. If  $\rho \in \cQ^{p,q}(\bF)$ and $\bG \hookrightarrow \bF$, Section~3.3 supplies $u^n \in \cU^{p,q}$ whose Young measures converge to $\eta^\rho$, $\sG_T$-stably in law (in the sense of~\cite{AldousEagleson1978}). In general, however, different approximating sequences of $\rho$ may correspond to different state limit laws; see also Remark~\ref{remark:limiting_relaxed_state_determined}. 

When jointly approximating several relaxed controls, the situation is even more delicate. A pair $(u,u')$ of strict controls on the same stochastic basis determines the joint Young measure
$$
    \eta^{u,u'}(\di t, \di a, \di a') \coloneq \di t\, \delta_{(u_t,u'_t)}(\di a, \di a')
$$
on $[0,T] \times \mA^2$.
Even for deterministic $\rho,\rho' \in \cQ^{p,q}$, different sequences of pairs of strict controls $(u,u')$ for which $\eta^u,\eta^{u'}$ approximate $\eta^\rho,\eta^{\rho'}$ may yield different joint Young measure limits, which may correspond to different joint state laws. 

\subsection{Parametrized relaxed controls and labeling schemes}
\label{sec:parametrized_relaxed_controls}
To intrinsically specify action couplings, we parametrize action distributions by the \textit{label} $r \in [0,1]$. 

Let $\sB \coloneq \sB[0,1]$ be the Borel $\sigma$-algebra on $[0,1]$.

\begin{definition}
    A $\sP^{\bF} \otimes \sB$-measurable $\pi : \oOmega \times [0,1] \rightarrow \mA$ is a \textit{parametrized relaxed control}.
\end{definition}

Each $\pi$ \textit{parametrizes} the relaxed control $\rho_t^\pi(\omega) \coloneq (\pi_t(\omega,\cdot\,))_\# \dint r$; predictability holds by $\int_\mA f(a) \,\rho^\pi_t(\di a) = \int_0^1 f(\pi_t(r)) \dint r$ for $f \in C_b(\mA)$. Conversely, every relaxed control admits a (generally non-unique) $\sP^\bF \otimes \sB$-measurable \textit{parametrization} \cite[Lemma~4.22]{kallenberg2021foundations}.

A finite family $\Pi = (\pi^i)_{i \in I}$ determines an action coupling $\lambda^\Pi$ of  $(\rho^{\pi^i})_{i \in I}$, where
\begin{equation}
    \label{eq:joint_coupling}
    \lambda^\Pi_t \coloneq (\Pi_t)_\# \dint r, \qquad \Pi_t(r) \coloneq \big(\pi_t^i(r)\big)_{i\in I}.
\end{equation}
These couplings $\lambda^\Pi$ are consistent under subfamily projection.

We call the map $\pi_t(\omega) : [0,1] \rightarrow \mA$ a \textit{profile}. The $q$-integrable \textit{profile space}~is
$$
    \mA_q \coloneq L^q([0,1];\mA),
$$
where $([0,1], \sB, \di r)$ is understood. Then $(\mA_q,\ma_0, d_{\mA_q})$ is a pointed Polish space. Write
$$
    \cR^{p,q} \coloneq L^{p,q}(\oOmega, \sP^{\bF}, \oP; \mA_q)
$$
with metric $d_{p,q}$, and set $\cR \coloneq \cR^{1,2}$. As usual, we write $\cR^{p,q}(\bF)$ and $\cR(\bF)$ when the ambient filtration needs clarification. For profiles $\phi,\phi' \in \mA_q$, $(\phi, \phi')_\# \dint r \in \Gamma(\phi_\# \dint r, \phi'_\# \dint r)$, hence $\cW_q(\phi_\# \dint r, \phi'_\# \dint r) \leqslant d_{\mA_q}(\phi,\phi')$. So $\cR^{p,q} \ni \pi \mapsto \rho^\pi \in \cQ^{p,q}$ is well defined, 1-Lipschitz, and
\begin{equation}
    \label{eq:relaxed_parametrized_size_comparison}
    d_{p,q}(\rho^\pi,\rho^{\pi'}) \leqslant d_{p,q}(\pi,\pi'), \qquad d_{p,q}(\pi,\ma_0) = d_{p,q}(\rho^\pi,\delta_{\ma_0}).
\end{equation}

Every $u \in \cU$ is canonically realized as the parametrized relaxed control $\pi^u_t(\omega,r) \coloneq u_t(\omega)$. For $u \in \cU^{p,q}(\bK)$, a $\bK$-predictable \textit{label process} $\vartheta : \oOmega \rightarrow [0,1]$, and $\pi \in \cR^{p,q}$, consider the random Young measures $\nu^{u,\vartheta}, \nu^\pi : \Omega \rightarrow \cY^T_q([0,1] \times \mA)$ given by
\begin{equation}
    \label{eq:Young_measures_pi}
    \nu^{u,\vartheta}(\di t, \di r, \di a) \coloneq \dint t \, \delta_{\vartheta_t}(\di r) \, \delta_{u_t}(\di a), \qquad \nu^\pi(\di t, \di r, \di a) \coloneq \dint t \dint r \, \delta_{\pi_t(r)}(\di a).
\end{equation}
Their projections into $\cY_q^T(\mA)$ are $\eta^u$ and $\eta^\pi \coloneq \eta^{\rho^\pi}$. Note $\nu^u(\di t, \di r, \di a) \coloneq \di t \dint r \, \delta_{u_t}(\di a)$ pairs the same~$u_t$ with every $r \in [0,1]$, and this label-action independence is preserved by narrow limits. By contrast, $\nu^{u,\vartheta}$ retains the dependence between the label and the action, like in $\nu^\pi$.

For $\rho, \rho' \in \cQ^{p,q}$, $\|\cW_q(\eta^\rho, \eta^{\rho'})\|_{L^p(\Omega)} \leqslant d_{p,q}(\rho,\rho')$ follows by finding a measurable optimal coupling. By matching times and labels, pairing $u_t,u_t'$, respectively $\pi_t(r),\pi'_t(r)$, and taking $L^p(\Omega)$-norms, we see that for $u,u' \in \cU^{p,q}(\bK)$, $\bK$-predictable $\vartheta : \oOmega \rightarrow [0,1]$, and $\pi,\pi' \in \cR^{p,q}$,
\begin{equation}
    \label{eq:Young_measure_map_continuous}
    \|\cW_q(\nu^{u,\vartheta}, \nu^{u',\vartheta})\|_{L^p(\Omega)} \leqslant d_{p,q}(u,u'),  \qquad \|\cW_q(\nu^\pi, \nu^{\pi'})\|_{L^p(\Omega)} \leqslant d_{p,q}(\pi,\pi').
\end{equation}
Thus $d_{p,q}$-convergence of controls yields $\cW_q$-convergence in $L^p(\Omega)$ of their Young measures.

For $\nu^{u,\vartheta}$ to approximate $\nu^\pi$, the time-label marginal $\di t \, \delta_{\vartheta_t}(\di r)$ should approximate $\di t \dint r$. 

\begin{definition}
    \label{def:labeling_scheme}
    A \textit{labeling scheme} $(\bK^n, \vartheta^n)_{n \geqslant 1}$ consists of pairs $(\bK^n,\vartheta^n)$, each defined on an extension of the ground space with $\bG \hookrightarrow \bK^n$, such that $\vartheta^n : \oOmega \rightarrow [0,1]$ is $\bK^n$-predictable and $d_{\BL}(\di t \, \delta_{\vartheta^n_t}(\di r),\di t \dint r) \rightarrow 0$ in probability.
\end{definition}

\begin{example}
    \label{example:standard_labeling_scheme}
    Define $\theta_n : [0,T] \rightarrow [0,1]$ by $\theta_n(t) \coloneq nt/T - \lfloor nt/T \rfloor$, the fractional part of~$nt/T$. Whenever $\bG \hookrightarrow \bK^n$, $(\bK^n,\theta_n)_{n \geqslant 1}$ is a labeling scheme; in particular, $(\bG,\theta_n)_{n \geqslant 1}$ is a labeling scheme. Indeed, $\int_0^T f(s,\theta_n(s)) \dint s = T/n \sum_{k=0}^{n-1} \int_0^1 f((k+r)T/n,r) \dint r$ for any 1-Lipschitz $f : [0,T] \times [0,1] \rightarrow \bR$. Comparing over each $[kT/n, (k+1)T/n]$ yields $|\int_0^T f(s,\theta_n(s)) \dint s - \int_0^T \int_0^1 f(s,r) \dint r \dint s| \leqslant T^2/n$, so $d_{\BL}(\di t\, \delta_{\theta_n(t)}(\di r), \di t \dint r) \leqslant T^2/n \rightarrow 0$.
\end{example}

Thus by choosing the labeling scheme $(\bG,\theta_n)_{n \geqslant 1}$, the subsequent approximation results hold with strict controls on the ground space. On the other hand, it may be desirable to approximate using randomized action selection; this is also supported in Definition~\ref{def:labeling_scheme}.

\begin{example}
    \label{example:sample_and_hold_labeling_scheme}
    Let $\Delta_n = \{0 = t_0^n < \dots < t_{J_n}^n = T\}$ be a finite partition of $[0,T]$, and write $|\Delta_n| \coloneq \max_{0 \leqslant j < J_n} |t_{j+1}^n - t_j^n|$ for its mesh. Pass to a product extension of the ground space carrying independent $\mU^n_j \sim \Unif[0,1], n \geqslant 1, 0 \leqslant j < J_n$, all independent of $\sG_T$. Let
    $$
        \sK^n_t \coloneq \sG_t \vee \sigma(\mU^n_j : t^n_j \leqslant t), \qquad \vartheta^n_t \coloneq \mU^n_j \quad \text{for } t \in (t_j^n,t_{j+1}^n],
    $$
    and $\vartheta^n_0 = 0$. Thus at $t_j^n$, the fresh label $\mU^n_j$ is revealed and held by $\vartheta^n$ along $(t_j^n,t_{j+1}^n]$. The filtration $\bK^n$ contains the cumulative label information together with the ground information. When $|\Delta_n| \rightarrow 0$, we call $(\bK^n,\vartheta^n)_{n \geqslant 1}$ a \textit{sample-and-hold} labeling scheme. Indeed, $\bG \hookrightarrow \bK^n$ by independence, $\vartheta^n : \oOmega \rightarrow [0,1]$ is $\bK^n$-predictable, and for $f \in C([0,T] \times [0,1])$,
    \begin{align*}
        E\bigg[\bigg|\int_0^T \int_0^1 f(s,r) \, \delta_{\vartheta^n_s}(\di r) \dint s - \int_0^T \int_0^1 f(s,r) \dint r \dint s \bigg|^2\bigg] &= \sum_{j=0}^{J_n-1} \Var\bigg[ \int_{t_j^n}^{t_{j+1}^n} f(s,\mU^n_j) \dint s\bigg] \\&\leqslant
        |\Delta_n|\int_0^T \int_0^1 |f(s,r)|^2 \dint r \dint s
    \end{align*}
    converges to 0 as $n \rightarrow \infty$. The family $\cF = \{f \in C([0,T] \times [0,1]) : \|f\|_\infty \leqslant 1, \, \Lip(f) \leqslant 1 \}$ is compact by Arzel\`a--Ascoli. So for any $\varepsilon > 0$, there is a finite $\varepsilon$-net $f_1,\dots,f_M$ of $\cF$, and the above shows that $\limsup_{n \rightarrow \infty}\|d_{\BL}(\di t \, \delta_{\vartheta_t^n}(\di r), \di t \dint r)\|_{L^2(\Omega)} \leqslant 2T\varepsilon$. Thus $(\bK^n, \vartheta^n)_{n \geqslant 1}$ is a labeling scheme, as claimed. We obtain \textit{sample-and-hold controls} when a profile is selected before each draw and held, at the sampled action, until the next grid time; see Example~\ref{example:sample_and_hold_control}.
\end{example}

\subsection{Stable convergence and martingale characteristics}

States corresponding to $\rho$ and $\pi$ are also defined on the extension $(\Omega, \sF_T,\bF,P)$. For consistency with~(\ref{eq:SDE_strict}), we aim to realize the control-state pairs as limits of strict control-state pairs, $\sG_T$-stably in law. Recall that for random elements $X_n,X$ taking values in a Polish space $E$, possibly defined on different extensions of the ground space, $X_n$ converges \textit{$\sG_T$-stably in law} to $X$, written $X_n \Rightarrow_{\sG_T} X$, if
$$
    E[Z f(X_n)] \longrightarrow E[Z f(X)]
$$
for all $Z \in L^\infty(\sG_T)$, $f \in C_b(E)$. Equivalently, for any $\sG_T$-measurable Polish-valued $Y$,
$$
    (Y,X_n) \Longrightarrow (Y,X);
$$
see~\cite[Definition~VIII.5.28, Proposition~VIII.5.33(ii)]{JacodShiryaev2003}. In particular, we may take $Y = W$. If $X_n \Rightarrow_{\sG_T} X$ and $\Law(X \giv \sG_T) = \Law(X' \giv \sG_T)$ (as kernels on the ground space), then $X_n \Rightarrow_{\sG_T} X'$. Random elements on the same probability space converging in probability converge $\sG_T$-stably in law. A tight sequence of Polish-valued random elements, possibly defined on different extensions of the ground space, has a $\sG_T$-stably convergent subsequence, with its limit realized on an extension of the ground space~\cite[Theorem~3.8 and Proposition~3.18]{JacodMemin1981}.

We record the following facts to be used later. For a Borel function $\Phi : E \rightarrow E'$, let $D_\Phi$ denote the Borel set of discontinuities of $\Phi$. 

\begin{lemma}
    \label{lemma:stable_convergence_facts}
    Let $E,E'$ be Polish. All variables below are $E$-valued random elements. 
    \begin{longlist}
        \item Let $\Phi : E \rightarrow E'$ be Borel. If $X_n \Rightarrow_{\sG_T} X$, $P[X \in D_\Phi] = 0$, then $\Phi(X_n) \Rightarrow_{\sG_T} \Phi(X)$.
        \item If $X_n \Rightarrow_{\sG_T} X$, $d_E(X_n,X_n') \rightarrow 0$ in probability,
        then $X_n' \Rightarrow_{\sG_T} X$.
        \item If $X_{k,n} \Rightarrow_{\sG_T} X_k$, $X_k \Rightarrow_{\sG_T} X$, there are nondecreasing $k_n \rightarrow \infty$ with $X_{k_n,n} \Rightarrow_{\sG_T} X$.
    \end{longlist}
\end{lemma}

\begin{proof}
    Let $Y$ be $\sG_T$-measurable and Polish-valued.
    For~(i), apply~\cite[Theorem~5.27]{kallenberg2021foundations} to $(Y,X_n)$ and the map $(y,x) \mapsto (y,\Phi(x))$. For~(ii), adjoin $Y$ and apply~\cite[Theorem~5.29]{kallenberg2021foundations}. For~(iii), enumerate $\{X_{k,n},X_k, X : k,n \geqslant 1\}$ as $(Y_n)_{n \geqslant 1}$, and let $\sH \subset \sG_T$ be the $\sigma$-algebra generated by $\Law(Y_n \giv \sG_T), n \geqslant 1$. As $E$ is Polish, $\sH$ is countably generated, and the topology of convergence $\sH$-stably in law is metrizable~\cite[Proposition~2.3]{Balder2001ws}, say with $d_{\st}$. Choose strictly increasing integers $N_k$ such that $d_{\st}(X_{k,n}, X_k) < 1/k$ for all $n \geqslant N_k$. Define nondecreasing $k_n \rightarrow \infty$ by $k_n = k$ for $N_k \leqslant n < N_{k+1}$ and $k_n = 1$ for $n < N_1$. Then $d_{\st}(X_{k_n,n}, X) \leqslant d_{\st}(X_{k_n,n}, X_{k_n}) + d_{\st}(X_{k_n},X) \leqslant 1/k_n + d_{\st}(X_{k_n}, X) \rightarrow 0$. Thus $X_{k_n,n} \Rightarrow_{\sH} X$, and~\cite[Proposition~3.11]{JacodMemin1981} shows $X_{k_n,n} \Rightarrow_{\sG_T} X$ by construction of~$\sH$.
\end{proof}

We now identify martingale characteristics of a $\sG_T$-stable limit of a sequence of strict control-state pairs, in the filtration generated by $\bG$ and the limiting history.

For $\mS \in \cS$, $t \in [0,T], x \in \bR^{m_\mS}, \mu \in \cP_2(\mA)$, define the \textit{averaged coefficients}
\begin{equation}
    \label{eq:ob_osigma}
    \ob_\mS(t,x,\mu) \coloneq \int_{\mA} b_\mS(t,x,a) \, \mu(\di a),\qquad \osigma_\mS(t,x,\mu) \coloneq \int_{\mA} \sigma_\mS(t,x,a) \, \mu(\di a),
\end{equation}
and, with $\sigma^\circ_\mS(t,x,a,\mu) \coloneq \sigma_\mS(t,x,a) - \osigma_\mS(t,x,\mu)$, the \textit{residual covariance matrix}
\begin{align}
    R_\mS(t,x,\mu) &\coloneq \int_\mA \sigma_\mS(t,x,a) \sigma_\mS(t,x,a)^\top \,\mu(\di a) - \osigma_\mS(t,x,\mu) \osigma_\mS(t,x,\mu)^\top \label{eq:R}
    \\&=
    \int_{\mA} \sigma^\circ_\mS(t,x,a,\mu) \sigma^\circ_\mS(t,x,a,\mu)^\top \, \mu(\di a). \label{eq:R_rewritten_centered}
\end{align}
Note that $R_\mS(t,x,\mu) \in \bS^{m_\mS}_+$ and $\tr R_\mS(t,x,\mu) \leqslant \int_\mA |\sigma^\circ_\mS(t,x,a,\delta_{\ma_0})|^2 \,\mu(\di a) \lesssim \cW_2(\mu,\delta_{\ma_0})^2$.

\begin{proposition}
    \label{prop:martingale_characteristics}
    Fix $\mS \in \cS$ and suppress it from the notation. For each $n$, let $\bK^n$ be a filtration on an extension of the ground stochastic basis with $\bG \hookrightarrow \bK^n$, and let $u^n \in \cU(\bK^n)$. Suppose that, on some extension of the ground space,
    \begin{equation}
        \label{eq:martingale_characteristics_convergence}
        (\eta^{u^n}, X^{u^n}) \Longrightarrow_{\sG_T} (\eta,Y)
    \end{equation}
    in $\cY^T_2(\mA) \times C([0,T];\bRm)$. Then $\bG \hookrightarrow \bF^{\eta,Y}$, where $\sF^{\eta,Y}_t \coloneq \sG_t \vee \sigma(\eta_{\leqslant t}, Y_{\cdot \wedge t})$. Write $\eta(\di t, \di a) = \di t \, \rho_t(\di a)$ for the $\bF^{\eta,Y}$-predictable disintegration of $\eta$ from Lemma~\ref{lemma:predictable_disintegration}, and let
    $$
        L_t \coloneq Y_t - x_0 - \int_0^t \ob(s,Y_s,\rho_s) \dint s, \qquad N_t \coloneq L_t - \int_0^t \osigma(s,Y_s,\rho_s) \dint W_s.
    $$
    Then $L \in \cM^c_{\loc,0}(\bF^{\eta,Y};\bR^m)$ with $\di [L]_t = (\int_{\mA} \sigma(t,Y_t,a) \sigma(t,Y_t,a)^\top \rho_t(\di a)) \dint t$, and for any $Z \in \cM_{\loc,0}(\bG)$, $\di [L,Z]_t = \osigma(t,Y_t,\rho_t) \dint [W,Z]_t$. Consequently, $N \in \cM^c_{\loc,0}(\bF^{\eta,Y};\bRm)$ and 
    \begin{equation}
        \label{eq:N_characteristics}
        \di[N]_t = R(t,Y_t,\rho_t) \dint t, \qquad [N,Z] \equiv 0.
    \end{equation}
\end{proposition}

\begin{proof}
    \textit{Step 1. Immersion.} As $\bG \hookrightarrow \bK^n$, $E[\xi F(\eta^{u^n}_{\leqslant s}, X^{u^n}_{\cdot \wedge s})(Z_t - Z_s)] = 0$ for any $s < t$, $\xi \in L^\infty(\sG_s)$, bounded continuous real $F$, and bounded real $\bG$-martingale $Z$. By Lemma~\ref{lemma:Young_history_continuity} and (\ref{eq:martingale_characteristics_convergence}), $E[\xi F(\eta_{\leqslant s}, Y_{\cdot \wedge s}) (Z_t - Z_s)] = 0$. A monotone class argument yields the martingale identity of $Z$ in the raw filtration for $\bF^{\eta,Y}$; as $Z$ is bounded and \cadlag, this extends by completion and right limits to $\bF^{\eta,Y}$. Hence $\bG \hookrightarrow \bF^{\eta,Y}$.
    
    \textit{Step 2. Continuity.} For $f = b,\sigma,\sigma\sigma^\top$, $K < \infty$, there is $C_K < \infty$ such that for $|y|, |y'| \leqslant K$, $|f(t,y,a)| \leqslant C_K(1 + d(a,\ma_0)^2)$ and $|f(t,y,a) - f(t,y',a)| \leqslant C_K( 1+ d(a,\ma_0))|y-y'|$. Let
    $$
        \Phi_f : (\mu,x) \longmapsto \bigg(\int_{[0,t] \times \mA} f(s,x_s,a) \,\mu(\di s, \di a)\bigg)_{0 \leqslant t\leqslant T}    
    $$ 
    be defined on $\cY^T_2(\mA) \times C([0,T];\bR^m)$. We claim $\Phi_f$ is continuous when the codomain is equipped with the topology of uniform convergence. Let $\cW_2(\mu_n,\mu) \rightarrow 0$, $x_n \rightarrow x$ uniformly, and $K \coloneq \sup_{n \geqslant 1} \|x_n\|_\infty \vee \|x\|_\infty < \infty$. By the Lipschitz bound on $f$ and $\cW_2$-convergence, $\|\Phi_f(\mu_n, x_n) - \Phi_f(\mu_n,x)\|_\infty \leqslant C_K \|x_n - x\|_\infty \int_{[0,T] \times \mA} ( 1 + d(a,\ma_0) ) \, \mu_n(\di s, \di a) \rightarrow 0$. It is left to show $\|\Phi_f(\mu_n,x) - \Phi_f(\mu,x)\|_\infty \rightarrow 0$. For $R > 0$, find continuous $\chi_R : \mA \rightarrow [0,1]$ with $\chi_R(a) = 1$ for $d(a,\ma_0) \leqslant R$ and $\chi_R(a) = 0$ for $d(a,\ma_0) \geqslant R+1$. Set $f_R(t,y,a) \coloneq f(t,y,a) \chi_R(a)$. By quadratic growth and $\cW_2$-convergence, $\lim_{R \rightarrow \infty} \|\Phi_{f - f_R}(\mu,x)\|_\infty = 0$ and $\lim_{R \rightarrow \infty} \sup_{n \geqslant 1} \|\Phi_{f-f_R}(\mu_n,x)\|_\infty = 0$. Moreover, the map $(s,a) \mapsto f_R(s,x_s,a)$ is bounded, Borel in $s$, and continuous in $a$. As the measures have the same time-marginal~$\di t$, \cite[Corollary~2.9]{JacodMemin1981} gives $\Phi_{f_R}(\mu_n,x)_t \rightarrow \Phi_{f_R}(\mu,x)_t$ for each $t$; uniform convergence follows from $|\Phi_{f_R}(\mu_n,x)_t - \Phi_{f_R}(\mu_n,x)_s| \leqslant C_{K,R}|t-s|$. This proves the claimed continuity of $\Phi_f$.

    \textit{Step 3. Limiting characteristics.} Set $A^n \coloneq \Phi_b(\eta^{u^n}, X^{u^n})$, $D^n \coloneq \Phi_\sigma(\eta^{u^n}, X^{u^n})$, and $C^n \coloneq \Phi_{\sigma\sigma^\top}(\eta^{u^n},X^{u^n})$. Define $A,D,C$ analogously from $(\eta,Y)$. Lemma~\ref{lemma:stable_convergence_facts}(i), (\ref{eq:martingale_characteristics_convergence}), continuity of the $\Phi_f$, and Lemma~\ref{lemma:Young_history_continuity} yield
    \begin{equation}
        \label{eq:martingale_characteristics_proved_convergence}
        (x_0,\eta^{u^n}_{\leqslant}, X^{u^n}, A^n, D^n, C^n) \Longrightarrow_{\sG_T} (x_0,\eta_{\leqslant}, Y, A,D,C).
    \end{equation}
    In particular, $Y_0 = x_0$ a.s. 
    Fix $Z \in \cM^{2,c}_0(\bG)$. The Kunita--Watanabe decomposition~\cite[Theorem~III.4.11]{JacodShiryaev2003} of $Z$ with respect to~$W$ yields
    $$
        Z = \kappa^\top \bigcdot W + Z^\perp,
    $$
    where $\kappa$ is $\bG$-predictable, $E[\int_0^T |\kappa_s|^2 \dint s] < \infty$, and $Z^\perp \in \cM^{2,c}_0(\bG)$ satisfies $[Z^\perp,W] \equiv 0$. For each $n$, set $L^n_t \coloneq X^{u^n}_t - x_0 - A^n_t = \int_0^t \sigma(s,X^{u^n}_s,u^n_s) \dint W_s$, $M^n \coloneq \col(L^n,W,Z^\perp)$, and $M \coloneq \col(L,W,Z^\perp)$. As $\bG \hookrightarrow \bK^n$, $M^n$ is a continuous local $\bK^n$-martingale null at zero, and
    $$
        [M^n]_t = \Gamma^n_t \coloneq \begin{pmatrix}
            C^n_t & D^n_t & 0 \\
            (D^n_t)^\top & t  I_d & 0 \\
            0 & 0 & [Z^\perp]_t
        \end{pmatrix}.
    $$
    Let $\Gamma$ be the analogous matrix with $C,D$ in place of $C^n,D^n$. For $s < T, G \in \sG_s$, set $V^{G,s}_t(\omega) \coloneq 1_G(\omega) 1_{[s,T]}(t)$. Then $V^{G,s}$ is $\bG$-adapted and \cadlag. Lemma~\ref{lemma:stable_convergence_facts}(i) and (\ref{eq:martingale_characteristics_convergence}) give $(V^{G,s},x_0, \eta^{u^n}_{\leqslant}, X^{u^n}, M^n, \Gamma^n) \Rightarrow (V^{G,s},x_0, \eta_{\leqslant}, Y, M, \Gamma)$. By~\cite[Theorem~IX.2.4 and Remark~IX.2.21]{JacodShiryaev2003} with auxiliary coordinates $(V^{G,s}, x_0, \eta^{u^n}_{\leqslant}, X^{u^n})$, $M$ is a continuous local martingale null at zero in the filtration $\bH^{G,s} \subset \bF^{\eta,Y}$ generated by the limiting coordinates, and $[M] = \Gamma$. The stopping times $\tau_k \coloneq \inf\{ t \geqslant 0 : |M_t| \geqslant k\} \wedge T$ depend only on $M$, and $M^{\tau_k}$ is a bounded $\bH^{G,s}$-martingale for all $s< T$ and $G \in \sG_s$. Hence $E[1_{G \cap F}(M_{t \wedge \tau_k} - M_{s \wedge \tau_k})] = 0$ for $t > s$ and $F \in \sigma(\eta_{\leqslant s}, Y_{\cdot \wedge s}) \subset \sH_s^{G,s}$. Since $\{G \cap F : G \in \sG_s, F \in \sigma(\eta_{\leqslant s}, Y_{\cdot \wedge s})\}$ is a $\pi$-system generating the raw $\sigma$-algebra defining $\sF^{\eta,Y}_s$, a monotone-class argument, followed by completion and right limits, shows that $M$ is a continuous local $\bF^{\eta,Y}$-martingale. Hence $L \in \cM_{\loc,0}^c(\bF^{\eta,Y};\bRm)$. As $[M] = \Gamma$, $\di [L]_t = (\int_\mA \sigma(t,Y_t,a) \sigma(t,Y_t,a)^\top \, \rho_t(\di a)) \dint t$, $\di [L,W]_t = \osigma(t,Y_t,\rho_t) \dint t$, and $[L,Z^\perp] \equiv 0$. Thus
    $$
        \di [L,Z]_t = \di [L, \kappa^\top \bigcdot W]_t 
        = \osigma(t,Y_t,\rho_t) \kappa_t \dint t = \osigma(t,Y_t,\rho_t) \dint [W,Z]_t.
    $$
    This gives the required formula for $[L,Z]$ for $Z \in \cM_0^{2,c}(\bG)$, and hence for $Z \in \cM_{\loc,0}^c(\bG)$ by localization. For $Z \in \cM_{\loc,0}(\bG)$, let $Z = Z^c + Z^d$ be its continuous-purely-discontinuous decomposition~\cite[Theorem~I.4.18]{JacodShiryaev2003}. By immersion, both components remain local $\bF^{\eta,Y}$-martingales, and $\di [L,Z]_t = \di [L,Z^c]_t = \osigma(t,Y_t,\rho_t) \dint[W,Z^c]_t = \osigma(t,Y_t,\rho_t) \dint [W,Z]_t$ since $[L,Z^d] = [W,Z^d] = 0$. This proves all claims for $L$; the claims for $N$ readily follow.
\end{proof}

We call $N = L - \osigma \bigcdot W$ the \textit{residual martingale} of $Y$. In the next section, the characteristics of $N$ in (\ref{eq:N_characteristics}) lead to representations of strict state limits by the \textit{residual dynamics}. We conclude with a tightness criterion that will be useful later.

\begin{lemma}
    \label{lemma:tight_Young_gives_tight_states}
    Fix $\mS \in \cS$ and suppress it from the notation. For each $n$, let $\bK^n$ be a filtration on an extension of the ground stochastic basis with $\bG \hookrightarrow \bK^n$, and let $u^n \in \cU(\bK^n)$. If the laws of $\eta^{u^n}$ are tight in $\cY_2^T(\mA)$, then the laws of $X^{u^n}$ are tight in $C([0,T];\bRm)$.
\end{lemma}

\begin{proof}
        Write $X^{u^n} = x_0 + D^n + M^n$, where $D^n,M^n$ are the drift and continuous local martingale parts, respectively. Set $A^n \coloneq \int_0^\cdot 1 + d(u^n_s,\ma_0)^2 \dint s = \Phi_f(\eta^{u^n},0)$ in the notation of Proposition~\ref{prop:martingale_characteristics}, with $f(t,x,a) = 1 + d(a,\ma_0)^2$. As the $\eta^{u^n}$ are tight in $\cY_2^T(\mA)$ and $\Phi_f$ is continuous by the proof of Proposition~\ref{prop:martingale_characteristics}, the $A^n$ are tight in $C[0,T]$. For $K > 0$, set $\tau^n_K \coloneq \inf\{t \geqslant 0 : A^n_t \geqslant K\} \wedge T$. Replacing $u^n$ by $\ma_0$ after $\tau^n_K$, Proposition~\ref{prop:strict_well_posed} and pathwise uniqueness give $\sup_n E[\|X^{u^n}_{\cdot \wedge \tau^n_K}\|_\infty]  \lesssim (1 + E[|x_0|] + \sqrt{K})$. So as $(A^n)_{n \geqslant 1}$ is tight and $P[\|X^{u^n}\|_\infty > R] \leqslant P[A^n_T \geqslant K] + C(1 + E[|x_0|] + \sqrt{K})/R$, $(\|X^{u^n}\|_\infty)_{n \geqslant 1}$ is bounded in probability. Now $\TV_{[s,t]}(D^n) + \tr([M^n]_t - [M^n]_s) \lesssim (1+\|X^{u^n}\|_\infty^2)(t-s) + A^n_t - A^n_s$ for $0 \leqslant s \leqslant t \leqslant T$ by linear growth of $b,\sigma$, where $\TV_{[s,t]}(D^n)$ is the total variation of $D^n$ over $[s,t]$. Thus $(D^n)_{n \geqslant 1}$ and $(\tr[M^n])_{n \geqslant 1}$ are tight as their supremum norms are bounded in probability and their moduli of continuity vanish uniformly in probability. Consequently, $(M^n)_{n \geqslant 1}$ is tight~\cite[Theorem~VI.4.13]{JacodShiryaev2003}, and therefore $(X^{u^n})_{n \geqslant 1}$ is tight.
\end{proof}

\section{Finite-dimensional residual dynamics}
\label{sec:3}
For a prescribed ordinary relaxed control, we provide an equation determining a corresponding state for which the control-state pair can be approximated by~(\ref{eq:SDE_strict}), $\sG_T$-stably in law. Conversely, we show that any state for which this approximation holds necessarily solves an equation of the same type on a further Wiener extension. Throughout this section, fix $\mS \in \cS$ and suppress it from the notation.

\subsection{The residual equation}
\label{sec:3.1}
Work on a filtered extension $(\Omega,\sF_T,\bF,P)$ of the ground basis. In this subsection, assume $\bG \hookrightarrow \bF$ and that the extension carries, for some $K \in \bN$, a $K$-dimensional $\bF$-Brownian motion $B$ independent of~$\sG_T$. In particular, $(W,B)$ is $\bF$-Brownian. 

In view of~(\ref{eq:N_characteristics}), we wish to represent the residual martingale $N$ through a square root of the residual covariance matrix $R$ and the orthogonal Brownian motion~$B$.

\begin{definition}
    \label{def:K_factor}
     A \textit{$K$-factor} $Q$ for $\rho \in \cQ$ is a $\sP^{\bF}\otimes \sB(\bRm)$-measurable map 
    $$
        Q : \oOmega \times \bRm \rightarrow \bR^{m \times K}
    $$
    such that, $\oP$-a.e., $Q_t(x) Q_t(x)^\top = R(t,x,\rho_t)$ and $|Q_t(x) - Q_t(x')| \lesssim |x-x'|$ for $x,x' \in \bRm$.
\end{definition}

For a $K$-factor $Q$ for $\rho \in \cQ$, the \textit{residual equation} corresponding to $(\rho,Q)$ is
\begin{equation}
    \label{eq:SDE_residual_factor}
    X^{\rho,Q}_t = x_0 + \int_0^t \ob(s,X_s^{\rho,Q}, \rho_s) \dint s + \int_0^t \osigma(s,X_s^{\rho,Q}, \rho_s) \dint W_s + \int_0^t Q_s(X_s^{\rho,Q}) \dint B_s.
\end{equation}
The canonical choice (verified in Proposition~\ref{prop:residual_well_posed} below) is $K = m$ and
\begin{equation}
    \label{eq:principal_factor}
    Q^\rho_t(x) \coloneq R(t,x,\rho_t)^{1/2}, \qquad X^\rho \coloneq X^{\rho,Q^\rho},
\end{equation}
where the $1/2$ is the principal square root. For $K > m$, $[R(t,x,\rho_t)^{1/2} \ 0_{m \times (K-m)}]$ is a $K$-factor, and for $K < m$, existence of a $K$-factor requires, in particular, $\rank R(t,x,\rho_t) \leqslant K$.

\begin{proposition}
    \label{prop:residual_well_posed}
    For each $K$-factor $Q$ for $\rho \in \cQ$, (\ref{eq:SDE_residual_factor}) has a unique $\bF$-strong solution $X^{\rho,Q}$; if $\rho \in \cQ^{p,2}$ and $x_0$ is $p$-integrable, $E[\|X^{\rho,Q}\|_\infty^p] < \infty$. In the case $\rho^u_t = \delta_{u_t}$ for some $u \in \cU$, (\ref{eq:SDE_residual_factor}) reduces to~(\ref{eq:SDE_strict}) and $X^{\rho^u,Q} = X^u$. The $Q^\rho$ in (\ref{eq:principal_factor}) is an $m$-factor, and
    \begin{equation}
        \label{eq:canonical_factor_regularity}
        |R(t,x,\mu)^{1/2} - R(t,x',\mu')^{1/2}| \lesssim |x-x'| + \cW_2(\mu,\mu').
    \end{equation}
    For $\rho,\rho' \in \cQ^{p,2}$, if $X^\rho,X^{\rho'}$ use the same residual driver $B$, $E[\|X^{\rho} - X^{\rho'}\|_\infty^p] \lesssim d_{p,2}(\rho,\rho')^p$.
\end{proposition}

\begin{proof}
    The averaged coefficients satisfy $|\ob(t,x,\mu)| + |\osigma(t,x,\mu)| \lesssim 1 + |x| + \cW_2(\mu,\delta_{\ma_0})$ and $|\ob(t,x,\mu) - \ob(t,x',\mu')| + |\osigma(t,x,\mu) - \osigma(t,x',\mu')| \lesssim |x-x'| + \cW_2(\mu,\mu')$. Moreover, every $K$-factor $Q$ satisfies $|Q_t(x)|^2 = \tr R(t,x,\rho_t) \lesssim \cW_2(\rho_t,\delta_{\ma_0})^2$, so together with the state-Lipschitz condition in Definition~\ref{def:K_factor}, the estimate of Proposition~\ref{prop:strict_well_posed} gives strong well-posedness and the claimed moment bound when $\rho \in \cQ^{p,2}$ and $x_0$ is $p$-integrable. Since  $\ob(t,x,\delta_{u_t}) = b(t,x,u_t)$, $\osigma(t,x,\delta_{u_t}) = \sigma(t,x,u_t)$, and $R(t,x,\delta_{u_t}) = 0$, (\ref{eq:SDE_residual_factor}) reduces to (\ref{eq:SDE_strict}) for $\rho^u_t = \delta_{u_t}$, and so $X^{\rho^u,Q} = X^u$.
    To establish~(\ref{eq:canonical_factor_regularity}), fix $t,x,x',\mu,\mu'$, and a coupling $\gamma \in \Gamma(\mu,\mu')$. Set $H_\gamma \coloneq L^2(\mA^2,\gamma;\bRd)$, and let $G_1,G_2 \in \cL_2(H_\gamma,\bRm)$ be given by
    $$
        G_1^* y(a,a') = \sigma^\circ(t,x,a,\mu)^\top y, \qquad G_2^* y(a,a') = \sigma^\circ(t,x',a',\mu')^\top y.
    $$
    Then $\|G_1 - G_2\|_{\HS}^2 = \int_{\mA^2} |\sigma^\circ(t,x,a,\mu) - \sigma^\circ(t,x',a',\mu')|^2 \, \gamma(\di a, \di a')$, hence
    $$
        \|G_1 - G_2\|_{\HS}^2 \leqslant \int_{\mA^2} |\sigma(t,x,a) - \sigma(t,x',a')|^2  \gamma(\di a, \di a') \lesssim |x-x'|^2 + \int_{\mA^2} d(a,a')^2 \gamma(\di a, \di a').
    $$
    As $|(G_1 G_1^*)^{1/2} - (G_2 G_2^*)^{1/2}| \leqslant \sqrt{2}\|G_1 - G_2\|_{\HS}$ \cite[Theorem~1]{araki1981inequality} and $G_1G_1^* = R(t,x,\mu)$, $G_2 G_2^* = R(t,x',\mu')$, taking the infimum over $\gamma \in \Gamma(\mu,\mu')$ yields~(\ref{eq:canonical_factor_regularity}). Thus $Q^\rho$ is an $m$-factor for $\rho$, and the stability estimate follows as in Proposition~\ref{prop:strict_well_posed}.
\end{proof}

We call the strong solution $X^{\rho,Q}$ of~(\ref{eq:SDE_residual_factor}) the \textit{state} corresponding to $(\rho,Q)$ (and $B$).

\subsection{Residual representation of a stable limit}

Return to the setting of Proposition~\ref{prop:martingale_characteristics}. We express the limiting state $Y$ in (\ref{eq:martingale_characteristics_convergence}) as the solution to an equation of the form~(\ref{eq:SDE_residual_factor}) on a further Wiener extension. To that end, we use the following Brownian independence and representation results.

\begin{lemma}
    \label{lemma:conditional_independence_BM}
    Let $\bG\hookrightarrow \bF$ and $B$ a $K$-dimensional $\bF$-Brownian motion. Then $B$ is independent of $\sG_T$ if and only if $[B,Z] \equiv 0$ for every $Z \in \cM_{\loc,0}(\bG)$. In this case, $\bG \vee \bF^B \hookrightarrow \bF$.
\end{lemma}

\begin{proof}
    Assume $B$ is independent of $\sG_T$. For $\tau^n \coloneq \inf\{t \geqslant 0 : |B_t| \geqslant n\} \wedge T$ and a bounded real $\bG$-martingale $Z$, $B^{\tau^n}$ and $ ZB^{\tau^n}$ are $\bG \vee \bF^B$-martingales. Hence $[B^{\tau^n}, Z] \equiv 0$, thus $[B,Z] \equiv 0$. Localizing the continuous part of $Z$ extends to $Z \in \cM_{\loc,0}(\bG)$.
    
    Conversely, suppose $[B,Z] \equiv 0$ for all bounded $Z$. Let $\xi \in L^\infty(\sG_T)$, $\zeta \in L^\infty(\sF^B_T)$, and set $Z_t \coloneq E[\xi \giv \sG_t]$, $M_t \coloneq E[\zeta \giv \sF^B_t]$. The martingale representation theorem provides an $\bF^B$-predictable $\psi$ with $M_t = E[\zeta] + \int_0^t \psi_s^\top \dint B_s$; $M$ is an $\bF$-martingale as $B$ is $\bF$-Brownian. Immersion asserts that~$Z$ is an $\bF$-martingale, and $[Z,M] \equiv 0$ by assumption. Thus $ZM$ is a bounded $\bF$-martingale, hence $E[\xi \zeta \giv \sF_t] = Z_t M_t$. Taking expectations yields $E[\xi \zeta] = E[\xi]E[\zeta]$, which shows that $B$ is independent of $\sG_T$. Moreover, as $Z_t M_t$ is $\sG_t \vee \sF^B_t$-measurable,  $E[\xi \zeta \giv \sF_t] = E[\xi \zeta \giv \sG_t \vee \sF_t^B]$ by the tower law. A monotone-class argument shows $E[V \giv \sF_t] = E[V \giv \sG_t \vee \sF^B_t]$ for $V \in L^\infty(\sG_T \vee \sF^B_T)$, which proves $\bG \vee \bF^B \hookrightarrow \bF$.
\end{proof}

\begin{lemma}
    \label{lemma:orthogonal_Brownian_representation}
    Let $\bG \hookrightarrow \bF$ and $N \in \cM^c_{\loc,0}(\bF;\bRm)$ with $\di [N]_t = V_t V_t^\top \dint t$ and $[N,Z] \equiv 0$ for an $\bF$-predictable $V : \oOmega \rightarrow \bR^{m \times K}$ and every $Z \in \cM_{\loc,0}(\bG)$. On a further filtered Wiener extension in which $\bF$ is immersed, there is a $K$-dimensional Brownian motion~$\tB$, with respect to the enlarged filtration, independent of $\sG_T$, such that $N = V \bigcdot \tB$.
\end{lemma}

\begin{proof}
    On a further filtered Wiener extension in which $\bF$ is immersed, let $B^0$ be a $K$-dimensional Brownian motion independent of $\sF_T$. Set $\tB \coloneq V^\dagger \bigcdot N + (I_K - V^\dagger V) \bigcdot B^0$. Since $V^\dagger V$ is an orthogonal projection, $V^\dagger VV^\top (V^\dagger)^\top = V^\dagger V$. Hence $\tB$ is well defined, $\di [V^\dagger \bigcdot N]_t = V_t^\dagger V_t \dint t$, and $\di[(I_K - V^\dagger V) \bigcdot B^0]_t = (I_K - V^\dagger_t V_t) \dint t$. Thus $\di [\tB]_t = I_K \dint t$ since $[N,B^0] \equiv 0$, so $\tB$ is a $K$-dimensional Brownian motion independent of $\sG_T$ by Lemma~\ref{lemma:conditional_independence_BM}. Now $V \bigcdot \tB$ is well defined as $\int_0^T |V_s|^2 \dint s = \tr [N]_T < \infty$ a.s. Moreover, $V(I_K - V^\dagger V) = 0$ and $(I_m - VV^\dagger) \bigcdot N = 0$ as $[(I_m - VV^\dagger) \bigcdot N] \equiv 0$. Thus $V \bigcdot \tB = VV^\dagger \bigcdot N = I_m \bigcdot N = N$.
\end{proof}

\begin{theorem}
    \label{theorem:residual_representation}
    Retain the assumptions of Proposition~\ref{prop:martingale_characteristics}. There is a further filtered Wiener extension in which $\bF^{\eta,Y}$ is immersed, carrying an $m$-dimensional Brownian motion~$\tB$ for the enlarged filtration, independent of $\sG_T$, such that
    $$
        Y_t = x_0 + \int_0^t \ob(s,Y_s,\rho_s) \dint s + \int_0^t \osigma(s,Y_s,\rho_s) \dint W_s + \int_0^t R(s,Y_s,\rho_s)^{1/2} \dint \tB_s.
    $$
\end{theorem}

\begin{proof}
    Apply Lemma~\ref{lemma:orthogonal_Brownian_representation} to (\ref{eq:N_characteristics}) with $V_t = R(t,Y_t,\rho_t)^{1/2}$.
\end{proof}

\begin{remark}
    \label{remark:limiting_relaxed_state_determined}
    Retain the assumptions of Proposition~\ref{prop:martingale_characteristics}. If $\rho \in \cQ(\bG)$, Theorem~\ref{theorem:residual_representation} and Proposition~\ref{prop:residual_well_posed} uniquely determine $\Law(Y\giv \sG_T)$, independently of the approximating sequence. For a general $\rho \in \cQ(\bF)$, independence of $\tB$ from $\sG_T$ does not specify the joint law of $\rho$ and~$\tB$, and one can check that $\Law(Y)$ need not be uniquely determined.
\end{remark}

\subsection{Strict approximation}

We show that the solution $X^{\rho,Q}$ to~(\ref{eq:SDE_residual_factor}) admits a $\sG_T$-stable in law strict approximation along any prescribed labeling scheme and time grids of mesh converging to zero. This is established first for a dense subclass of elementary parametrized relaxed controls, and then extended to general inputs by stability of the principal factor~(\ref{eq:principal_factor}).

Retain the stochastic basis and Brownian-driver assumptions of Section~\ref{sec:3.1}. 
For a partition $\fp = \{I_1,\dots,I_L\}$ of $[0,1]$ into intervals of positive length, let $\cA_\fp \subset \mA_2$ denote the profiles $\phi : [0,1] \rightarrow \mA$, up to $\di r$-equality, that are constant on each $I_\ell$ $\di r$-a.e.; $\phi \in \cA_\fp$ are evaluated using their interval-step representatives. Put
$$
    \cR^{p,q}_{\fp} \coloneq \{\pi \in \cR^{p,q} : \pi_t(\omega) \in \cA_\fp \ \oP\text{-a.e.} \}, \qquad \cR_\fp \coloneq \cR^{1,2}_\fp.
$$
We say $\pi \in \cR_\fp$ is \textit{$\fp$-elementary} if there are $0 = t_0 < \dots < t_J = T$ and, for each $0 \leqslant j < J$, finite-range $\sF_{t_j}$-measurable $\varphi_j : \Omega \rightarrow \cA_\fp$ such that $\pi_t = \varphi_j$ for $t \in (t_j,t_{j+1}]$ $\oP$-a.e.~and $\pi_0 \equiv \ma_0$. A $\pi \in \cR$ is \textit{elementary} if it is $\fp$-elementary for some $\fp$. An elementary $\pi$ assumes only finitely many action values, and therefore belongs to every $\cR^{p,q}$.

\begin{lemma}
    \label{lemma:elementary_dense}
    \begin{longlist}
        \item The elementary parametrized relaxed controls are dense in $\cR^{p,q}$. 
        \item The $\fp$-elementary parametrized relaxed controls are dense in $\cR^{p,q}_\fp$.
        \item If $\bF = \bG \vee \bF^B$, the dense $\pi$ in (i)--(ii) can be chosen with $\varphi_j(\omega) = \phi^{\alpha,\beta}_j \in \cA_\fp$ for $\omega \in G_{j,\alpha} \cap \{(B_{t_1},\dots,B_{t_j}) \in C_{j,\beta}\}$, where $G_{j,1},\dots,G_{j,N_j} \in \sG_{t_j}$ is a partition of $\Omega$, and $C_{j,1},\dots,C_{j,K_j}$ is a Borel partition of $(\bR^K)^j$ satisfying $P[(B_{t_1},\dots,B_{t_j}) \in \partial C_{j,\beta}] = 0$.
    \end{longlist}
\end{lemma}

\begin{proof}
    For (i), choose a countable dense family $(\phi_\ell)_{\ell \geqslant 1}$ of $\mA_q$, where $\phi_1 \equiv \ma_0$ and every $\phi_\ell$ belongs to $\cA_\fp$ for some $\fp$; for (ii), choose such a family dense in $(\cA_\fp, d_{\mA_q})$ instead. Given $\pi^0$ to be approximated, let $\pi^n_t(\omega)$ be the closest element in $(\phi_\ell)_{\ell=1}^n$ to $\pi^0_t(\omega)$, resolving ties by the smallest index. Then $\pi^n \in \cR^{p,q}$, $d_{\mA_q}(\pi^n_t, \pi^0_t) \rightarrow 0$, and $d_{\mA_q}(\pi^n_t,\pi^0_t) \leqslant d_{\mA_q}(\pi^0_t,\ma_0)$; dominated convergence thus yields $d_{p,q}(\pi^n,\pi^0) \rightarrow 0$. We therefore assume without loss of generality that $\pi^0$ takes values in $(\phi_\ell)_{\ell=1}^n$ for some fixed $n$. Approximate each of the finitely many level sets of $\pi^0$ in $\oP$-measure by finite unions of sets of the form $F \times (t,t'], F \in \sF_t$, and disjointify; assign $\ma_0$ on the complement. The resulting $\pi$ belongs to $\cR_\fp$, where $\fp$ is a common refinement for which $(\phi_\ell)_{\ell=1}^n \subset \cA_\fp$; for (ii), all profiles already belong to the prescribed $\cA_\fp$. Refining the time endpoints shows that $\pi$ is of the required elementary form. Thus (i)--(ii) follow by dominated convergence, by sufficiently approximating the level sets.

    For (iii), the algebra generated by $\sG_t$ and the events $\{B^\ell_s \leqslant b\}$, where $0 < s \leqslant t$, $b \in \bR$, $1 \leqslant \ell \leqslant K$, generates $\sF_t$ up to completion. Construct $\pi$ in the same way as in (i)--(ii), but for each $F \times (t,t']$ require $F$ to belong to this generating algebra. Refining the decision times to include the finitely many observation times, and refining the ground events and Brownian events into finite partitions shows $\pi$ is of the required form; the continuity property holds because each $B^\ell_s, s > 0$, has an atomless distribution. 
\end{proof}

Let $\pi \in \cR_\fp$, where $\fp = \{I_1,\dots,I_L\}$; assume $L \geqslant 2$ by refining the partition if necessary. Fix an orthonormal basis $1,e_1,\dots,e_{L-1}$ of the subspace of $L^2[0,1]$ consisting of functions constant on each $I_\ell$. By orthogonality with $1$, $\int_0^1 e_\ell(r) \dint r = 0$ for all $1 \leqslant \ell \leqslant L-1$. Set 
$$
    K_\fp \coloneq d(L-1), \qquad h_\fp(r) \coloneq [e_1(r)I_d \ \dots \ e_{L-1}(r)I_d] \in \bR^{d \times K_\fp}.
$$
Note $\int_0^1 h_\fp(r) \dint r = 0$ and $\int_0^1 h_\fp(r)^\top h_\fp(r) \dint r = I_{K_\fp}$. Define the $K_\fp$-factor $Q^\pi_\fp$ for $\rho^\pi$ by
\begin{equation}
    \label{eq:elementary_factor}
    Q^\pi_{\fp} : \oOmega \times \bR^m \rightarrow \bR^{m \times K_\fp}, \quad Q^\pi_{\fp,t}(x) \coloneq \int_0^1 \sigma\big( t,x,\pi_t(r) \big) h_\fp(r) \dint r \in \bR^{m \times K_\fp}.
\end{equation}
Indeed, $Q^\pi_{\fp}$ is $\sP^\bF \otimes \sB(\bR^m)$-measurable, $Q_{\fp,t}^\pi(x) Q_{\fp,t}^\pi(x)^\top = R(t,x,\rho^\pi_t)$ by~(\ref{eq:R_rewritten_centered}) and Parseval's identity, and $|Q^\pi_{\fp,t}(x) - Q^{\pi}_{\fp,t}(x')| \lesssim |x-x'|$ by Bessel's inequality. 

Let $\bK$ be a filtration on an extension with $\bG \hookrightarrow \bK$, and $\vartheta : \oOmega \rightarrow [0,1]$ a $\bK$-predictable label process. We call $u$ a \textit{$\fp$-elementary label control} relative to $\vartheta$ if, for some finite partition $0 = t_0 < \dots < t_J = T$ of $[0,T]$, $u_t = \varphi_j(\vartheta_t)$ for $t \in (t_j,t_{j+1}]$ $\oP$-a.e.\ and $u_0 = \ma_0$, where $\varphi_j(\omega) = \phi_j^{\alpha,\beta} \in \cA_\fp$ for $\omega \in G_{j,\alpha} \cap \{(Z^\vartheta_{t_1},\dots,Z^\vartheta_{t_j}) \in C_{j,\beta}\}$. Here $Z^{\vartheta,\ell}_t \coloneq \int_0^t 1_{I_\ell}(\vartheta_s) \dint W_s$, $Z^\vartheta_t \coloneq \col_{1 \leqslant \ell \leqslant L} Z^{\vartheta,\ell}_t$, $G_{j,1},\dots,G_{j,N_j} \in \sG_{t_j}$ is a partition of $\Omega$, and $C_{j,1},\dots,C_{j,K_j}$ is a Borel partition of $(\bR^{dL})^j$. We call $u$ \textit{elementary} if it is $\fp$-elementary for some $\fp$. Each elementary~$u$ is $\bK$-predictable and assumes finitely many values in $\mA$, hence belongs to $\cU^{p,q}(\bK)$. 

\begin{example}
    \label{example:sample_and_hold_control}
    Let $(\bK^n, \vartheta^n)_{n \geqslant 1}$ be a sample-and-hold labeling scheme in the sense of Example~\ref{example:sample_and_hold_labeling_scheme}. For $0 \leqslant j < J_n$, denote
    $$
        \sD^n_j \coloneq \sG_{t^n_j} \vee \sigma(\mU^n_\ell : \ell < j)
    $$
    the information available before the fresh label $\mU^n_j$ is revealed. 
    We call $u^n \in \cU^{p,q}(\bK^n)$ a \textit{sample-and-hold control} on $\Delta_n$ if for each $0 \leqslant j < J_n$, there are $\sD^n_j $-measurable (random) profiles $\varphi^n_j : \Omega \rightarrow \mA_q$ for which, using $\sD^n_j \otimes \sB$-measurable representatives of $\varphi_j^n$, 
    $$
        u^n_t = \varphi^n_j(\mU^n_j), \quad t \in (t_j^n, t_{j+1}^n].
    $$
    Thus the profile $\varphi^n_j$ is selected before the fresh label $\mU^n_j$ is revealed, and the resulting sampled action $\varphi_j^n(\mU^n_j)$ is held until the next grid time. 
    
    Suppose the decision times of an elementary label control $u$ relative to $\vartheta^n$ belong to $\Delta_n$. By repeating the selected profiles, write $u_t = \varphi_j(\vartheta_t^n) = \varphi_j(\mU^n_j)$ for $t \in (t_j^n,t_{j+1}^n]$. Then $\varphi_j$ is $\sD^n_j$-measurable as $Z^{\vartheta^n,\ell}_{t^n_j} = \sum_{i=0}^{j-1} 1_{I_\ell}(\mU^n_i)(W_{t^n_{i+1}} - W_{t^n_i})$ is. So $u$ is sample-and-hold on $\Delta_n$.
\end{example}

\begin{proposition}
    \label{prop:elementary_strictification}
    Fix a partition $\fp = \{I_1,\dots,I_L\}$ of $[0,1]$, $L \geqslant 2$, and $\pi \in \cR^{p,q}_\fp(\bF)$. Assume $B$ in~(\ref{eq:SDE_residual_factor}) has dimension $K = K_\fp$. Let $(\bK^n,\vartheta^n)_{n \geqslant 1}$ be a labeling scheme, and define
    $B^n_t \coloneq \int_0^t h_\fp(\vartheta^n_s)^\top \dint W_s$ on the extension carrying $(\bK^n,\vartheta^n)$. Then $B^n \Rightarrow_{\sG_T} B$, and   
    \begin{longlist}
        \item (stability) for any $u^n \in \cU^{p,q}(\bK^n)$ satisfying
        \begin{equation}
            \label{eq:elementary_strictification_input_convergence}
            (B^n,\nu^{u^n,\vartheta^n}) \Longrightarrow_{\sG_T} (B,\nu^\pi)
        \end{equation}
        in $C([0,T];\bR^K_\fp) \times \cY_q^T([0,1] \times \mA)$, we have 
        \begin{equation}
            \label{eq:elementary_strictification_output_convergence}
            (B^n,\nu^{u^n,\vartheta^n}, X^{u^n}) \Longrightarrow_{\sG_T} (B,\nu^\pi,X^{\rho^\pi,Q^\pi_\fp});
        \end{equation}
        \item (elementary approximation) for finite partitions $(\Delta_n)_{n \geqslant 1}$ of $[0,T]$ with $|\Delta_n| \rightarrow 0$, there are $\fp$-elementary $u^n \in \cU^{p,q}(\bK^n)$ relative to $\vartheta^n$, with decision times belonging to $\Delta_n$, satisfying~(\ref{eq:elementary_strictification_input_convergence}) and hence~(\ref{eq:elementary_strictification_output_convergence}).
    \end{longlist}

\end{proposition}

\begin{proof}
    \textit{Step 1. Residual driver approximation.} Abbreviate $h \coloneq h_\fp$. As $h$ is bounded, $E[|B^n_t - B^n_s|^4] \lesssim |t-s|^2$, so $(B^n)_{n \geqslant 1}$ is tight. Definition~\ref{def:labeling_scheme}, the $\di t \dint r$-a.e.~continuity of $h$, and the continuous mapping theorem yield $\di t \, \delta_{h(\vartheta^n_t)}(\di y) \rightarrow \di t \, (h_\#\di r)(\di y)$ in probability in $\cY_2^T(\bR^{d\times K})$. Thus for every subsequence, there is a further subsequence, not relabeled, for which $(\di t \, \delta_{h(\vartheta^n_t)}(\di y), B^n)$ converges $\sG_T$-stably in law, with first marginal $\di t \, ((h)_\#\di r)(\di y)$. Proposition~\ref{prop:martingale_characteristics} applied with strict controls $h(\vartheta^n)$ on the action space $\bR^{d \times K}$ and the system with state dimension~$K$, zero initial condition and drift, and diffusion coefficient $y \mapsto y^\top$, followed by Lemma~\ref{lemma:conditional_independence_BM}, shows the limit is $K$-dimensional Brownian motion independent of $\sG_T$. As the original subsequence was chosen arbitrarily, this gives $B^n \Rightarrow_{\sG_T} B$.

    \textit{Step 2. Stability.} Suppose (\ref{eq:elementary_strictification_input_convergence}) holds. Then $(B^n,\nu^{u^n,\vartheta^n},X^{u^n})_{n \geqslant 1}$ is tight by Lemma~\ref{lemma:tight_Young_gives_tight_states}. Hence for any subsequence, there is a further subsequence, not relabeled, such that 
    \begin{equation}
        \label{eq:elementary_strictification_1}
        (B^n,\nu^{u^n,\vartheta^n}, X^{u^n}) \Longrightarrow_{\sG_T} (\hB,\hnu, \hX),
    \end{equation}
    where $\Law(\hB,\hnu \giv \sG_T) = \Law(B,\nu^\pi \giv \sG_T)$. Since $L^q([0,T];\cA_\fp) \ni \phi \mapsto \nu^\phi \in \cY^T_q([0,1] \times \mA)$ is continuous and injective, its image is Borel and its inverse on that image is Borel by the Lusin--Suslin theorem, so $\hnu$ belongs to that image a.s. Let $\sK_t \coloneq \sG_t \vee \sigma(\hnu_{\leqslant t}, \hX_{\cdot \wedge t}, \hB_{\cdot \wedge t})$. On each $I_\ell$, applying Lemma~\ref{lemma:predictable_disintegration} to $\heta^\ell(\di t, \di a) \coloneq (1/|I_\ell|) \hnu(\di t, I_\ell,\di a)$, which has $\bK$-adapted history, and applying the Lusin--Suslin theorem again to the Dirac embedding, shows that $\hnu = \nu^{\hpi}$ for some $\hpi \in \cR^{p,q}_{\fp}(\bK)$. Moreover, the preceding Borel bijection also shows, together with $\Law(\hB, \hnu \giv \sG_T) = \Law(B,\nu^\pi \giv \sG_T)$, that $\Law(\hB,\hpi \giv \sG_T) = \Law(B,\pi \giv \sG_T)$. Now define $v^n_t = (\vartheta^n_t, u^n_t, h(\vartheta^n_t))$. As $h$ is bounded and $\di t \dint r$-a.e.\ continuous, Lemma~\ref{lemma:stable_convergence_facts}(i) gives
    $$
        \eta^{v^n} \Longrightarrow_{\sG_T} \heta, \qquad \heta(\di t, \di r, \di a, \di y) = \di t \dint r \, \delta_{\hpi_t(r)}(\di a) \, \delta_{h(r)}(\di y),
    $$
    jointly with~(\ref{eq:elementary_strictification_1}). Consider the system on $[0,1] \times \mA \times \bR^{d \times K}$ with state dimension $m + K$, initial condition $\col(x_0, 0_K)$, and coefficients $\tb(t,(x,x'), (r,a,y)) \coloneq \col(b(t,x,a), 0_K)$ and $\tsigma(t,(x,x'), (r,a,y)) \coloneq \col(\sigma(t,x,a), y^\top)$. Form the strict input with $v^n$; the corresponding state is $\col(X^{u^n}, B^n)$. The histories of $\heta$ and $\hnu$ generate the same filtration, since $\hnu$ is the $(t,r,a)$-marginal of $\heta$, and $\heta$ is a deterministic pushforward of $\hnu$. Proposition~\ref{prop:martingale_characteristics} shows that $\bG \hookrightarrow \bK$, $(W,\hB)$ is $\bK$-Brownian, and the residual martingale $\hN$ of $\hX$, for $\rho^{\hpi}$, satisfies $\di [\hN]_t = \hQ_t \hQ_t^\top \dint t$ and $\di [\hN,\hB]_t = \hQ_t \dint t$, where $\hQ_t \coloneq Q^{\hpi}_{\fp,t}(\hX_t)$. Therefore $[\hN - \hQ \bigcdot \hB] \equiv 0$, hence $\hN = \hQ \bigcdot \hB$, so $\hX$ solves~(\ref{eq:SDE_residual_factor}) for $(\rho^{\hpi}, Q^{\hpi}_{\fp})$ and drivers $(W,\hB)$. Since $\Law(x_0,W,\hB,\hpi \giv \sG_T) = \Law(x_0,W,B,\pi \giv \sG_T)$, and the corresponding states solve~(\ref{eq:SDE_residual_factor}) with factors $Q^{\hpi}_\fp$ and $Q^\pi_\fp$, respectively, strong well-posedness by Proposition~\ref{prop:residual_well_posed} yields $\Law(\hB,\hnu,\hX \giv \sG_T) = \Law(B,\nu^\pi, X^{\rho^\pi,Q^\pi_\fp} \giv \sG_T)$. As the subsequence was chosen arbitrarily, this and~(\ref{eq:elementary_strictification_1}) give (\ref{eq:elementary_strictification_output_convergence}). It remains to establish~(ii).

    \textit{Step 3. $\bG \vee \bF^B$-predictable $\pi$.} Suppose $\pi$ is $\bH$-predictable, where $\bH \coloneq \bG \vee \bF^B$. By Lemma~\ref{lemma:conditional_independence_BM}, $\bH \hookrightarrow \bF$. By Lemma~\ref{lemma:elementary_dense}(ii)--(iii), (\ref{eq:Young_measure_map_continuous}), and Lemma~\ref{lemma:stable_convergence_facts}(iii), assume without loss of generality that $\pi : \oOmega \rightarrow \cA_\fp$ is of the form in Lemma~\ref{lemma:elementary_dense}(iii). In that notation, write $0 = t_0 < \dots < t_J= T$ for the decision times and $\varphi = (\varphi_0,\dots,\varphi_{J-1}) \in D^J$ for the randomly selected profiles, where $D \subset \cA_\fp$ is finite and contains $\ma_0$. Round the decision times up to $\tau_j^n \coloneq \min\{t \in \Delta_n : t \geqslant t_j\} \in \Delta_n$. For all $n$ large enough, $0 = \tau_0^n < \dots < \tau_J^n = T$; work with these $n$, and set $u^n \equiv \ma_0$ for all smaller $n$. Select $\varphi^n_j$ by the same rule as $\varphi_j$, replacing each $B_{t_i}, i \leqslant j$, by $B^n_{\tau^n_i}$, and set $u^n_t \coloneq \varphi^n_j(\vartheta_t^n)$ for $t \in (\tau^n_j,\tau^n_{j+1}]$ and $u^n_0 \coloneq \ma_0$. Each $u^n$ is $\fp$-elementary with the properties in (ii), since its ground events belong to $\sG_{\tau^n_j}$ and $B^n_t = \sum_{\ell=1}^L (h|_{I_\ell})^\top Z^{\vartheta^n,\ell}_t$, from which we see that the selection events have the right form. Since $(B^n,(\tau^n_j)_{j=1}^{J-1}) \Rightarrow_{\sG_T} (B,(t_j)_{j=1}^{J-1})$ and $(x,s) \mapsto (x(s_1),\dots,x(s_{J-1}))$ is continuous, the continuity in Lemma~\ref{lemma:elementary_dense}(iii), Lemma~\ref{lemma:stable_convergence_facts}(i), and including the finitely many $\sG_T$-measurable indicators in the $\sG_T$-stable limit give, with $\varphi^n \coloneq (\varphi^n_0,\dots,\varphi^n_{J-1})$,
    \begin{equation}
        \label{eq:elementary_strictification_convergence_1}
        (B^n, \varphi^n) \Longrightarrow_{\sG_T} (B,\varphi).
    \end{equation}
    For $\phi = (\phi_0,\dots,\phi_{J-1}) \in D^J$, set $\pi^{\phi}_t \coloneq \phi_j$ for $t \in (t_j,t_{j+1}]$ and $\pi^{\phi}_0 \coloneq \ma_0$, and set $\nu^{\phi} \coloneq \nu^{\pi^{\phi}}$. As $(t,r) \mapsto (t,r,\pi^\phi_t(r))$ is $\di t \dint r$-a.e.~continuous, Definition~\ref{def:labeling_scheme} and the continuous mapping theorem yield $\max_{\phi \in D^J} \cW_q(\di t \, \delta_{\vartheta^n_t}(\di r) \, \delta_{\pi^\phi_t(\vartheta^n_t)}(\di a),\nu^\phi) \rightarrow 0$ in probability---Wasserstein convergence follows because all these measures have a common compact support. Also, $u^n$ differs from $\pi_t^{\varphi^n}(\vartheta^n_t)$ only on the intervals $(t_j,\tau^n_j]$, whose total length is at most $J|\Delta_n|$. So coupling the rounded and unrounded executions at the same $(t,r)$ bounds their $\cW_q$-distance by $(CJ|\Delta_n|)^{1/q}$; the preceding convergence then gives $\cW_q(\nu^{u^n,\vartheta^n},\nu^{\varphi^n}) \rightarrow 0$ in probability. As $\phi \mapsto \nu^{\phi}$ is continuous on $D^J$ and $\nu^\varphi = \nu^\pi$ a.s., (\ref{eq:elementary_strictification_convergence_1}) and Lemma~\ref{lemma:stable_convergence_facts}(i)--(ii) yield (\ref{eq:elementary_strictification_input_convergence}).

    \textit{Step 4. $\bF$-predictable $\pi$.} By Lemma~\ref{lemma:elementary_dense}(ii), (\ref{eq:Young_measure_map_continuous}), and Lemma~\ref{lemma:stable_convergence_facts}(iii), it remains to treat an $\bF$-predictable $\fp$-elementary $\pi$, with decision times $t_j$ and selected tuple $\varphi$, using the same notation as Step~3. Order $D$ and equip the set of words in $D$ with lexicographic order~$\leq_{\lex}$. For $\phi \in D^J$, set
    $$
        c_\phi \coloneq P[\varphi \leq_{\lex} \phi \giv \sH_T], \qquad \Psi(\omega,r) \coloneq \min_{\lex}\{\phi \in D^J : r \leqslant c_\phi(\omega)\}, \quad r \in [0,1].
    $$
    Write $\phi_{\leqslant j} \coloneq (\phi_0,\dots,\phi_j)$ for $\phi \in D^J$, and likewise $\Psi_{\leqslant j}(\omega,r)$. As $\varphi_{\leqslant j}$ is $\sF_{t_j}$-measurable and $\bH \hookrightarrow \bF$, $\Law(\varphi_{\leqslant j} \giv \sH_T) = \Law(\varphi_{\leqslant j} \giv \sH_{t_j})$. 
    By construction of the lexicographic order, $\{\Psi_{\leqslant j}(\omega,r) \leq_{\lex} \phi_{\leqslant j}\} = \{r \leqslant P[\varphi_{\leqslant j} \leq_{\lex} \phi_{\leqslant j} \giv \sH_{t_j}](\omega)\}$, so $\Psi_j$ is $\sH_{t_j} \otimes \sB$-measurable. Pass to a further uniform product extension carrying $U \sim \Unif[0,1]$ independent of $\sF_T$, and put $\hvarphi(\omega) \coloneq \Psi(\omega,U(\omega))$. Then $\Law(\hvarphi \giv \sH_T) = \Law(\varphi \giv \sH_T)$. Let $\Phi$ be the standard normal distribution function. For $\varepsilon \in (0,t_1)$, put $U_\varepsilon \coloneq \Phi(W^1_\varepsilon/\sqrt{\varepsilon}) \sim \Unif[0,1]$, which is independent of $\sH_0$. Right continuity of $\bH$ gives, for $Z \in L^\infty(\sH_T)$ and $f \in C[0,1]$, that $E[Z f(U_\varepsilon)] - E[Z] \int_0^1 f(r) \dint r = E[(E[Z \giv \sH_\varepsilon] - E[Z \giv \sH_0])f(U_\varepsilon)] \rightarrow 0$. Thus $U_\varepsilon \Rightarrow_{\sH_T} U$. Set $\varphi^\varepsilon(\omega) \coloneq \Psi(\omega,U_\varepsilon(\omega))$. Viewed as a measurable function of $(c,r)$, where $c = (c_\phi)_{\phi \in D^J}$, it is continuous whenever $r \neq c_\phi$ for every $\phi$. As $c$ is $\sH_T$-measurable and  $U \neq c_\phi$ for all~$\phi$ a.s., Lemma~\ref{lemma:stable_convergence_facts}(i) gives $\varphi^\varepsilon \Rightarrow_{\sH_T} \hvarphi$; so $\varphi^\varepsilon \Rightarrow_{\sH_T} \varphi$ as $\Law(\hvarphi \giv \sH_T) = \Law(\varphi \giv \sH_T)$. Set 
    \begin{align*}
        \pi^\varepsilon_t \coloneq \begin{cases*}
            \ma_0, &$t \leqslant \varepsilon$, \\
           \pi^{\varphi^\varepsilon}_t, &$t > \varepsilon$.
        \end{cases*}
    \end{align*}
    Since $\varphi^\varepsilon_j$ is $\sH_{t_j \vee \varepsilon}$-measurable, $\pi^\varepsilon$ is $\bH$-predictable and $\fp$-elementary. Continuity of $\phi \mapsto \nu^\phi$ on $D^J$, $\sH_T$-measurability of $B$, and Lemma~\ref{lemma:stable_convergence_facts}(i) give $(B,\nu^{\varphi^\varepsilon}) \Rightarrow_{\sG_T} (B,\nu^\pi)$. Coupling at the same time and label yields $\cW_q(\nu^{\pi^\varepsilon}, \nu^{\varphi^\varepsilon})^q \leqslant \varepsilon \max_{\phi \in D} d_{\cA_q}(\ma_0,\phi)^q$, hence Lemma~\ref{lemma:stable_convergence_facts}(ii) gives $(B,\nu^{\pi^\varepsilon}) \Rightarrow_{\sG_T} (B,\nu^\pi)$. Now choose $\varepsilon_k \downarrow
     0$. For each $k$, Step~3 provides $\fp$-elementary label controls $u^{k,n} \in \cU^{p,q}(\bK^n)$ relative to $\vartheta^n$, with decision times in $\Delta_n$, and satisfying $(B^n,\nu^{u^{k,n},\vartheta^n}) \Rightarrow_{\sG_T} (B,\nu^{\pi^{\varepsilon_k}})$. Lemma~\ref{lemma:stable_convergence_facts}(iii) completes the proof.
\end{proof}

\begin{theorem}
    \label{theorem:residual_approximation}
    Assume $\bG \hookrightarrow \bF$. Let $Q$ be a $K$-factor for $\rho \in \cQ^{p,q}(\bF)$, and let $X^{\rho,Q}$ be the solution of~(\ref{eq:SDE_residual_factor}). Fix a parametrization $\pi$ of $\rho$. Let $(\bK^n,\vartheta^n)_{n \geqslant 1}$ be a labeling scheme and $(\Delta_n)_{n \geqslant 1}$ finite partitions of $[0,T]$ with $|\Delta_n| \rightarrow 0$. There are elementary label controls $u^n \in \cU^{p,q}(\bK^n)$ relative to $\vartheta^n$, with decision times belonging to $\Delta_n$, such that
    \begin{equation}
        \label{eq:residual_approximation_convergence}
        (\nu^{u^n,\vartheta^n}, X^{u^n}) \Longrightarrow_{\sG_T} (\nu^\pi, X^{\rho,Q})
    \end{equation}
    in $\cY^T_q([0,1] \times \mA) \times C([0,T]; \bRm)$. Consequently,
    \begin{equation}
        \label{eq:residual_approximation_convergence_ordinary}
        (\eta^{u^n}, X^{u^n}) \Longrightarrow_{\sG_T} (\eta^\rho, X^{\rho,Q})
    \end{equation}
    in $\cY^T_q(\mA) \times C([0,T];\bRm)$. 
\end{theorem}

\begin{proof}
    By Lemmas~\ref{lemma:conditional_independence_BM}--\ref{lemma:orthogonal_Brownian_representation}, there is a further filtered Wiener extension in which $\bF$ is immersed, and on it an $m$-dimensional Brownian motion $B^0$, independent of $\sG_T$, such that $X^{\rho,Q}$ solves~(\ref{eq:SDE_residual_factor}) for  $Q^\rho$ in~(\ref{eq:principal_factor}) and drivers $(W,B^0)$. Now $\pi \in \cR^{p,q}$ by~(\ref{eq:relaxed_parametrized_size_comparison}), hence Lemma~\ref{lemma:elementary_dense}(i) provides $\fp^k$-elementary $\pi^k$, with $d_{p,q}(\pi^k,\pi) \rightarrow 0$. Abbreviate $\rho^k \coloneq \rho^{\pi^k}$ and let $X^k \coloneq X^{\rho^k}$ be the solution of (\ref{eq:SDE_residual_factor}) for $(\rho^k, Q^{\rho^k})$ and $(W,B^0)$. By Proposition~\ref{prop:residual_well_posed}, (\ref{eq:relaxed_parametrized_size_comparison}), and (\ref{eq:Young_measure_map_continuous}),
    \begin{equation}
        \label{eq:residual_approximation_1}
        E[\|X^k - X^{\rho,Q}\|_\infty^p] \longrightarrow 0, \qquad \|\cW_q(\nu^{\pi^k}, \nu^\pi)\|_{L^p(\Omega)} \longrightarrow 0.
    \end{equation}
    For each $k$, $N^k_t \coloneq \int_0^t R(s,X^k_s,\rho^k_s)^{1/2} \dint B^0_s$ satisfies $\di [N^k]_t = Q^{\pi^k}_{\fp_k,t}(X^k_t) Q^{\pi^k}_{\fp_k,t}(X^k_t)^\top \dint t$ and $[N^k,Z] \equiv 0$ for all $Z \in \cM_{\loc,0}(\bG)$. Lemma~\ref{lemma:orthogonal_Brownian_representation} yields a further Wiener extension in which $\bF$ is immersed, and on it a $K_{\fp_k}$-dimensional Brownian motion $B^k$, independent of $\sG_T$, such that $X^k$ solves~(\ref{eq:SDE_residual_factor}) for $(\rho^k,Q_{\fp_k}^{\pi^k})$ and $(W,B^k)$. Thus Proposition~\ref{prop:elementary_strictification}(ii) applies, and together with~(\ref{eq:residual_approximation_1}) and Lemma~\ref{lemma:stable_convergence_facts}(iii), (\ref{eq:residual_approximation_convergence}) follows; (\ref{eq:residual_approximation_convergence_ordinary}) is obtained by projection.
\end{proof}

\begin{remark}
    \label{remark:residual_approximation_multi_dimensional}
    Theorem~\ref{theorem:residual_approximation} applies to a fixed finite family if the limiting \textit{joint} relaxed control on the product action space is prescribed. For finite $I$, let $\mS^{\stk}$ be the stacked system of $(\mS^i)_{i \in I} \subset \cS$, and take $\lambda \in \cQ^{p,q}(\bF;\mA^I)$, $\rho^i_t \coloneq (\pr^I_i)_\# \lambda_t$. The stacked residual covariance matrix is $R^\lambda = (R^\lambda_{ij})_{i,j \in I}$, with $(i,j)$-block
    $$
        R_{ij}^\lambda(t,x) = \int_{\mA^I} \sigma^\circ_{\mS^i}(t,x_i,a_i,\rho^i_t) \sigma^\circ_{\mS^j}(t,x_j,a_j,\rho^j_t)^\top \, \lambda_t(\di a),
    $$
    where $x = \col_{i \in I} x_i$. In particular, $R^\lambda_{ii}(t,x) = R_{\mS^i}(t,x_i,\rho^i_t)$.

    Fix a parametrization $\Pi$ of $\lambda$, a $K$-factor $Q$ for $(\lambda,\mS^{\stk})$, and for the prescribed drivers $W,B$, write $X^{\lambda,Q;\mS^{\stk}}$ for the corresponding state. For any labeling scheme $(\bK^n,\vartheta^n)_{n \geqslant 1}$ and finite partitions $(\Delta_n)_{n \geqslant 1}$ of $[0,T]$ with $|\Delta_n| \rightarrow 0$, Theorem~\ref{theorem:residual_approximation} over $\mA^I$ and~(\ref{eq:system_stacked_state}) provide elementary $u^{n,i} \in \cU^{p,q}(\bK^n)$ relative to $\vartheta^n$, with decision times contained in $\Delta_n$, such that
    $$
        \big( \nu^{U^n, \vartheta^n}, \col_{i \in I} X^{u^{n,i};\mS^i} \big) \Longrightarrow_{\sG_T} \big( \nu^\Pi, X^{\lambda,Q;\mS^{\stk}} \big),
    $$
    where $U^n \coloneq (u^{n,i})_{i \in I}$, and the Young measures are in $\cY_q^T([0,1] \times \mA^I)$ and defined as in~(\ref{eq:Young_measures_pi}), for the product action space $\mA^I$. Conversely, for any $U^n \coloneq (u^{n,i})_{i \in I} \subset \cU(\bK^n)$ with $(\eta^{U^n}, \col_{i \in I} X^{u^{n,i};\mS^i}) \Rightarrow_{\sG_T} (\eta^\lambda, Y)$, Theorem~\ref{theorem:residual_representation} applied to $\mS^{\stk}$ and the factor $(R^\lambda)^{1/2}$ gives the corresponding residual representation on a further Wiener extension.
\end{remark}

Such $\lambda$ in Remark~\ref{remark:residual_approximation_multi_dimensional} are specified  by parametrized relaxed controls; see~(\ref{eq:joint_coupling}).

\section{Parametrized cylindrical dynamics}
\label{sec:4} 
On a fixed extension of the ground basis, a state is specified for every $(\pi,\mS) \in \cR \times \cS$ as the unique strong solution to an equation driven by a common cylindrical Brownian motion. Any finite collection of these states can be jointly approximated by~(\ref{eq:SDE_strict}), $\sG_T$-stably in law. A strict-to-relaxed stability criterion is also provided. Conversely, any parametrized relaxed state assignment is identified as the strong solution map of an equation of the same form, under the pairwise approximation property. 

\subsection{The equation}

Let $(\Omega, \sF_T, \bF, P)$ be an extension of the ground basis. A \textit{cylindrical continuous local $\bF$-martingale} $\mM$ in a separable Hilbert space $H$ is a continuous linear map
$$
    \mM : H \longrightarrow \cM_{\loc,0}^c(\bF),
$$
where the codomain is equipped with the topology of uniform convergence in probability, or equivalently, $M^n \rightarrow M$ if $[M^n - M]_T \rightarrow 0$ in probability; see~\cite[Section~3.1]{VeraarYaroslavtsev2016} and also~\cite{VanNeervenVeraarWeis2007}. A predictable $Q : \oOmega \rightarrow \cL(H)$ with self-adjoint, positive semidefinite values is a \textit{covariance density} of $\mM$ if
$$
    \di[\mM(h), \mM(h')]_t = \langle Q_t h, h' \rangle \dint t, \quad h,h' \in H.
$$
When $Q \equiv I_H$, $\mM$ is \textit{cylindrical $\bF$-Brownian motion}~\cite[Chapter~4]{DaPratoZabczyk2014}. We call $\mM$ a \textit{cylindrical continuous square-integrable $\bF$-martingale} if $\mM$ additionally takes values in $\cM_0^{2,c}(\bF)$. 

Let $\mM$ be a cylindrical continuous local $\bF$-martingale in $H$ with covariance density $Q$. A predictable $\Phi : \oOmega \rightarrow \cL_2(H,\bRm)$ is \textit{$\mM$-integrable} when $\int_0^T \|\Phi_s Q^{1/2}_s\|_{\HS}^2 \dint s < \infty$ a.s.; it admits the stochastic integral $\Phi \bigcdot \mM \coloneq \int_0^\cdot \Phi_s \dint \mM_s \in \cM^c_{\loc,0}(\bF;\bRm)$ in this case. For $\mM$-integrable $\Phi : \oOmega \rightarrow \cL_2(H,\bRm)$, $\Psi : \oOmega \rightarrow \cL_2(H,\bR^{m'})$, and $h \in H$, \cite[Proposition~3.6]{Yaroslavtsev2018} yields
\begin{equation}
    \label{eq:cylindrical_covariation}
    \di [\Phi \bigcdot \mM, \Psi \bigcdot \mM]_t = \Phi_t Q_t \Psi^*_t \dint t, \qquad \di [\Phi \bigcdot \mM, \mM(h)]_t = \Phi_t Q_t h \dint t.
\end{equation}
When additionally $E[\int_0^T \|\Phi_s Q_s^{1/2}\|_{\HS}^2 \dint s] < \infty$, it holds $\Phi \bigcdot \mM \in \cM^{2,c}_0(\bF;\bRm)$. The following result identifies an element of $\cM_{\loc,0}^c(\bF;\bRm)$ with $\Phi \bigcdot \mM$ under analogous identities to (\ref{eq:cylindrical_covariation}).
\begin{lemma}
    \label{lemma:cylindrical_integral_identification}
    Let $\mM$ be a cylindrical continuous local $\bF$-martingale in $H$ with predictable covariance density $Q$. If $\Phi : \oOmega \rightarrow \cL_2(H,\bRm)$ is predictable and $L \in \cM_{\loc,0}^c(\bF;\bRm)$ satisfies $\di [L]_t = \Phi_t Q_t \Phi_t^* \dint t$ and $\di [L, \mM(h) ]_t = \Phi_t Q_t h \dint t$ for all $h \in H$, then $L = \Phi \bigcdot \mM$.
\end{lemma}

\begin{proof}
    As $\int_0^T \|\Phi_s Q_s^{1/2}\|_{\HS}^2 \dint s = \tr [L]_T < \infty$ a.s., $I \coloneq \Phi \bigcdot \mM$ is well defined. Extend the identity $\di [L, \mM(h) ]_t = \Phi_t Q_t h \dint t$ for $h \in H$ to $\di [L,\Psi \bigcdot \mM]_t = \Phi_t Q_t \Psi^*_t \dint t$ for $\mM$-integrable $\Psi : \Omega \rightarrow \cL_2(H,\bR^{m'})$. This is readily seen for $\Psi^*_t \epsilon'_j \coloneq \sum_{\ell=1}^N \xi^{\ell,j}_t h_\ell$ for $1 \leqslant j \leqslant m'$, $h_\ell \in H$, and bounded predictable $\xi^{\ell,j} : \oOmega \rightarrow \bR$, where $\epsilon'_1,\dots,\epsilon'_{m'}$ is the standard basis of $\bR^{m'}$. For general $\mM$-integrable $\Psi$, find $\Psi^n$ of the preceding form with $\int_0^T \|(\Psi_s - \Psi^n_s)Q^{1/2}_s\|^2_{\HS} \dint s \rightarrow 0$ in probability~\cite[proof of Proposition~3.6]{Yaroslavtsev2018}; the claimed identity follows from the Kunita--Watanabe inequality. Setting $\Psi = \Phi$ yields $[L] = [I] = [L,I]$, hence $L = I$.
\end{proof}

Consider the separable Hilbert spaces
\begin{equation}
    \label{eq:H,H_0}
    \mH \coloneq L^2([0,1];\bRd) = \bRd \oplus \mH_0, \qquad \mH_0 \coloneq \bigg\{ h \in \mH : \int_0^1 h(r) \dint r = 0 \bigg\},
\end{equation}
where $\bRd$ is identified with the constant functions. Henceforth specialize to $H = \mH$. We say a cylindrical local martingale $\mM$ in $\mH$ \textit{lies above} $W$ if $\mM(c) = c^\top W$ for all $c \in \bRd$. 

For the remainder of this subsection and throughout Section~\ref{sec:4.2}, assume $\bG \hookrightarrow \bF$ and that the extension carries a cylindrical $\bF$-Brownian motion $\mB$ in $\mH$, lying above $W$, with $\mB|_{\mH_0}$ independent of $\sG_T$. Fix an orthonormal basis $(\me_k)_{k \geqslant 1}$ of $\mH_0$, and set $\beta^k \coloneq \mB(\me_k)$. The $(\beta^k)_{k \geqslant 1}$ are independent real Brownian motions, jointly independent of $\sG_T$, and 
\begin{equation}
    \label{eq:cylindrical_BM}
    \mB(c+h) = c^\top W + \sum_{k=1}^\infty \langle h,\me_k \rangle \beta^k, \quad c + h \in \bRd \oplus \mH_0 = \mH,
\end{equation}
where the series converges in $L^2(\Omega;C[0,T])$ by Doob's inequality. Such $\mB$ exists: adjoining independent real Brownian motions $(\beta^k)_{k \geqslant 1}$, jointly independent of $\sG_T$, (\ref{eq:cylindrical_BM}) constructs $\mB$ on a countable Wiener extension. 

Applying Lemma~\ref{lemma:conditional_independence_BM} to $(\mB(\me_1), \dots,\mB(\me_n))$ for each $n \geqslant 1$ shows that $[\mB(h),Z] \equiv 0$ for all $h \in \mH_0$ and $Z \in \cM_{\loc,0}(\bG)$. 

For $\mS \in \cS, t \in [0,T], x \in \bR^{m_\mS}, \phi \in \mA_2$, define $\Sigma_{\mS}(t,x,\phi) \in \cL_2(\mH,\bR^{m_\mS})$ by
$$
    \Sigma_\mS(t,x,\phi) h \coloneq \int_0^1 \sigma_\mS\big(t,x,\phi(r) \big) h(r) \dint r.
$$
Thus $\Sigma_\mS(t,x,\phi)$ and $\sigma_\mS(t,x,\phi(\cdot))$ are identified under the canonical isometric isomorphism $\cL_2(\mH, \bR^{m_\mS}) \cong L^2([0,1];\bR^{m_\mS \times d})$. In particular, $\|\Sigma_\mS(t,x,\phi)\|_{\HS}^2 = \int_0^1 |\sigma_\mS(t,x,\phi(r))|^2 \dint r$ and 
\begin{equation}
    \label{eq:diffusion_operator_adjoint}
    \big(\Sigma_\mS(t,x,\phi)^* y\big)(r) = \sigma_\mS\big(t,x,\phi(r)\big)^\top y.
\end{equation} 
As $(t,x,\phi) \mapsto \Sigma_\mS(t,x,\phi)$ is Borel, $\Sigma_\mS(t,X_t,\pi_t)$ is predictable for continuous adapted $X$ and predictable $\pi$. For $(\pi,\mS) \in \cR \times \cS$, the \textit{parametrized cylindrical equation} is
\begin{equation}
    \label{eq:SDE_cylindrical}
    X^{\pi;\mS}_t = x^{\mS}_0 + \int_0^t \ob_{\mS}(s,X^{\pi;\mS}_s,\rho^\pi_s) \dint s + \int_0^t \Sigma_{\mS}(s,X^{\pi;\mS}_s,\pi_s) \dint \mB_s,
\end{equation}
where $\rho^\pi_s \coloneq (\pi_s)_\# \dint r$ as usual. 

\begin{proposition}
    \label{prop:cylindrical_well_posed}
    For $(\pi,\mS) \in \cR \times \cS$, (\ref{eq:SDE_cylindrical}) admits a unique $\bF$-strong solution $X^{\pi;\mS}$; if $\pi \in \cR^{p,2}$ and $x_0^\mS$ is $p$-integrable, $E[\|X^{\pi;\mS}\|_\infty^p] \lesssim 1 + E[|x_0^\mS|^p] + d_{p,2}(\pi,\ma_0)^p$. For $\pi,\pi' \in \cR^{p,2}$, $E[\|X^{\pi;\mS} - X^{\pi';\mS}\|^p_\infty] \lesssim d_{p,2}(\pi,\pi')^p$. If $\pi^u_t \equiv u_t$, (\ref{eq:SDE_cylindrical}) reduces to~(\ref{eq:SDE_strict}) and $X^{\pi^u;\mS} = X^{u;\mS}$.
\end{proposition}

\begin{proof}
    The assumptions on $\sigma_\mS$ and relation $\cL_2(\mH, \bR^{m_\mS}) \cong L^2([0,1];\bR^{m_\mS \times d})$ yield 
    $\|\Sigma_\mS(t,x,\phi)\|_{\HS} \lesssim 1 + |x| + d_{\mA_2}(\phi,\ma_0)$, $\|\Sigma_\mS(t,x,\phi) - \Sigma_\mS(t,x',\phi')\|_{\HS} \lesssim |x-x'| + d_{\mA_2}(\phi,\phi')$; $\ob_\mS$ satisfies the analogous estimates. Thus, the strong well-posedness and stability properties follow as in Proposition~\ref{prop:strict_well_posed}. The final claim is a consequence of $\mB$ lying above $W$ and $\Sigma_\mS(t,x,\pi^u_t)(c+h) = \sigma_\mS(t,x,u_t)c$ for $c + h \in \bRd \oplus \mH_0$.
\end{proof}

\begin{remark}
    \label{remark:multiple_parametrized_states}
    Regard a finite family $\Pi = (\pi^i)_{i \in I} \subset \cR^{p,q}$ as the parametrized relaxed control $\Pi_t(r) \coloneq (\pi_t^i(r))_{i \in I}$ on $\mA^I$ with induced relaxed control $\lambda^\Pi$ as in~(\ref{eq:joint_coupling}). As $I$ is finite, $\Pi \in \cR^{p,q}(\bF;\mA^I)$. Let $(\mS^i)_{i \in I} \subset \cS$ and write $\mS^{\stk}$ for  the stacked system. Then $\ob_{\stk}(t,x,\lambda^\Pi_t) \coloneq \int_{\mA^I} b_{\stk}(t,x,a) \, \lambda^\Pi_t(\di a) = \col_{i \in I} \ob_{\mS^i}(t,x_i,\rho_t^{\pi^i})$ and $\Sigma_{\stk}(t,x,\Pi_t) h \coloneq \int_0^1 \sigma_{\stk}(t,x,\Pi_t(r)) h(r) \dint r = \col_{i \in I} (\Sigma_{\mS^i}(t,x_i,\pi_t^i)h)$. Thus by Proposition~\ref{prop:cylindrical_well_posed},
    \begin{equation}
        \label{eq:system_stacked_parametrized_state}
        X^{\Pi;\mS^{\stk}} = \col_{i \in I}X^{\pi^i;\mS^i},
    \end{equation}
    where $X^{\Pi;\mS^{\stk}}$ is the unique strong solution of~(\ref{eq:SDE_cylindrical}) for $(\Pi,\mS^{\stk})$ over the action space $\mA^I$. Moreover, by (\ref{eq:cylindrical_covariation}) and (\ref{eq:diffusion_operator_adjoint}), 
    \begin{equation}
        \label{eq:action_coupling_state_characteristics}
        \di [X^{\pi^i;\mS^i},X^{\pi^{i'};\mS^{i'}}]_t = \bigg( \int_{\mA^I} \sigma_{\mS^i}(t,X^{\pi^i;\mS^i}_t,a_i) \sigma_{\mS^{i'}}(t,X^{\pi^{i'};\mS^{i'}}_t, a_{i'})^\top \, \lambda^\Pi_t(\di a) \bigg) \dint t
    \end{equation}
    for $i,i' \in I$, which shows that the covariation between controlled states depends on $\lambda^\Pi$.
\end{remark}

Using $\mH = \bRd \oplus \mH_0$, the diffusion coefficient in (\ref{eq:SDE_cylindrical}) decomposes into ground-driven and residual components. Indeed, with $P_{\mH_0}$ the orthogonal projection of $\mH$ onto $\mH_0$ and $c \in \bRd$,
\begin{equation}
    \label{eq:diffusion_decomposition}
    \Sigma_\mS(t,x,\pi_t)c = \osigma_\mS(t,x,\rho_t^\pi) c, \qquad \Sigma_\mS(t,x,\pi_t) P_{\mH_0} \Sigma_\mS(t,x,\pi_t)^* = R_\mS(t,x,\rho^\pi_t).
\end{equation}
Thus the directions in $\bRd$ reproduce the averaged diffusion against $W$, and the directions in $\mH_0$ account for the residual covariance~(\ref{eq:N_characteristics}). The following result shows that given a state and parametrization, the equations (\ref{eq:SDE_residual_factor}) and (\ref{eq:SDE_cylindrical}) are equivalent up to a Wiener extension.

\begin{proposition}
    \label{prop:residual_cylindrical_correspondence}
    Let $\bG \hookrightarrow \bF$. Fix $\rho \in \cQ$, a parametrization $\pi$ of $\rho$, a $K$-factor $Q$ for $\rho$, and $\mS \in \cS$ suppressed from the notation. Then $X^\pi$ admits a representation of the form (\ref{eq:SDE_residual_factor}) for $(\rho,Q)$ on a further Wiener extension in which $\bF$ is immersed. Conversely, for a prescribed residual driver $B$ independent of $\sG_T$, $X^{\rho,Q}$ admits a representation of the form~(\ref{eq:SDE_cylindrical}) on a further Wiener extension in which $\bF$ is immersed, for a cylindrical Brownian motion $\hmB$ in $\mH$, lying above $W$, with $\hmB|_{\mH_0}$ independent of~$\sG_T$.
\end{proposition}

\begin{proof}
    For the state $X$ under consideration, abbreviate $\Sigma_t \coloneq \Sigma(t,X_t,\pi_t)$, $G_t \coloneq \Sigma_t|_{\mH_0}$, and $R_t \coloneq G_t G_t^* = R(t,X_t,\rho_t)$. These are all predictable.

    Take $X = X^\pi$. Its residual martingale is $N \coloneq G \bigcdot \mB|_{\mH_0}$, and $\di [N]_t = R_t \dint t$ by~(\ref{eq:cylindrical_covariation}). By the same argument as in Lemma~\ref{lemma:cylindrical_integral_identification}, for $Z \in \cM_{\loc,0}(\bG)$, $[\mB(h), Z] = 0$ for every $h \in \mH_0$ extends to $[N,Z] \equiv 0$. Since $R_t = Q_t(X_t) Q_t(X_t)^\top$, Lemma~\ref{lemma:orthogonal_Brownian_representation} provides the required representation.

    Take $X \coloneq X^{\rho,Q}$. Its residual martingale is $N \coloneq Q(X) \bigcdot B$. Again, $\di [N]_t = R_t \dint t$ and $[N,Z] \equiv 0$. Also, $P_{\ran G_t^*} = G^*_t R^\dagger_t G_t$ is the orthogonal projection onto $\ran G_t^*$. On a further Wiener extension in which $\bF$ is immersed, adjoin a cylindrical Brownian motion $\tmB$ in $\mH_0$, independent of $\sF_T$. For $h \in \mH_0$, define
    $$
        \mB^0_t(h) \coloneq \int_0^t (R^\dagger_s G_s h)^\top \dint N_s + \int_0^t \langle (I_{\mH_0}-P_{\ran G^*_s}) h, \di \tmB_s \rangle.
    $$
    The two integrals are strongly orthogonal, and $G^*_t R^\dagger_t R_t R^\dagger_t G_t = P_{\ran G_t^*}$ shows they are well defined with bracket densities $\langle h, P_{\ran G^*_t}h \rangle \dint t$ and $\langle h, (I_{\mH_0} - P_{\ran G^*_t})h \rangle \dint t$, respectively. By polarization, $\di [\mB^0(h), \mB^0(h')]_t = \langle h,h'\rangle \dint t$, so $\mB^0$ is a cylindrical Brownian motion in $\mH_0$ on the extension. Also, $[\mB^0(h),Z] \equiv 0$ for $h \in \mH_0, Z \in \cM_{\loc,0}(\bG)$, so $\mB^0$ is independent of $\sG_T$ by Lemma~\ref{lemma:conditional_independence_BM}. As $\ran G_t = \ran R_t$, $\di [N,\mB^0(h)]_t = R_t R_t^\dagger G_t h \dint t = G_t h \dint t$. So $N = G \bigcdot \mB^0$ by $\di [N]_t = G_t G_t^* \dint t$ and Lemma~\ref{lemma:cylindrical_integral_identification}. Setting $\hmB(c + h) \coloneq c^\top W + \mB^0(h)$ for $c + h \in \bRd \oplus \mH_0$, $\hmB$ is a cylindrical Brownian motion in $\mH$, lying above $W$, with $\hmB|_{\mH_0}$ independent of $\sG_T$. Since $\Sigma_t(c+h) = \osigma(t,X_t,\rho_t) c + G_t h$, $X$ satisfies~(\ref{eq:SDE_cylindrical}) for $\hmB$. This completes the proof.
\end{proof}

\subsection{Strict approximation and stability} 
\label{sec:4.2}

\begin{theorem}
    \label{theorem:cylindrical_approximation}
    Let $(\bK^n,\vartheta^n)_{n \geqslant 1}$ be a labeling scheme and $(\Delta_n)_{n \geqslant 1}$ finite partitions of $[0,T]$ with $|\Delta_n| \rightarrow 0$. For any finite families $\Pi \coloneq (\pi^i)_{i \in I} \subset \cR^{p,q}$ and $(\mS^i)_{i \in I} \subset \cS$, there are elementary label controls $U^n \coloneq (u^{n,i})_{i \in I} \subset \cU^{p,q}(\bK^n)$ relative to $\vartheta^n$, with decision times belonging to $\Delta_n$, such that
    \begin{equation}
        \label{eq:cylindrical_approximation}
        \big(\nu^{U^n, \vartheta^n}, (X^{u^{n,i}; \mS^i})_{i \in I}\big) \Longrightarrow_{\sG_T} \big(\nu^{\Pi}, (X^{\pi^i;\mS^i})_{i \in I}\big)
    \end{equation}
    in $\cY_q^T([0,1] \times \mA^I) \times \prod_{i\in I} C([0,T];\bR^{m_{\mS^i}})$.
\end{theorem}

\begin{proof} 
    By Proposition~\ref{prop:residual_cylindrical_correspondence} with the principal factor for $\lambda^\Pi$, and Theorem~\ref{theorem:residual_approximation} with action space~$\mA^I$, there are $U^n$ of the required form satisfying $(\nu^{U^n,\vartheta^n}, X^{U^n;\mS^{\stk}}) \Rightarrow_{\sG_T} (\nu^\Pi, X^{\Pi;\mS^{\stk}})$. The conclusion follows from~(\ref{eq:system_stacked_state}) and (\ref{eq:system_stacked_parametrized_state}).
\end{proof}

\begin{corollary}
    \label{cor:common_control_parametrized_approximation}
    Let $(\bK^n,\vartheta^n)_{n \geqslant 1}$ be a labeling scheme and $(\Delta_n)_{n \geqslant 1}$ finite partitions of $[0,T]$ with $|\Delta_n| \rightarrow 0$. For any $\pi \in \cR^{p,q}$ and finitely many $(\mS^i)_{i \in I} \subset \cS$, there are elementary label controls $u^n \in \cU^{p,q}(\bK^n)$ relative to $\vartheta^n$, with decision times belonging to $\Delta_n$, such that
    $$
        \big(\nu^{u^n, \vartheta^n}, (X^{u^n; \mS^i})_{i \in I}\big) \Longrightarrow_{\sG_T} \big(\nu^{\pi}, (X^{\pi;\mS^i})_{i \in I}\big)
    $$
    in $\cY_q^T([0,1] \times \mA) \times \prod_{i\in I} C([0,T];\bR^{m_{\mS^i}})$.
\end{corollary}

\begin{proof}
    Use Theorem~\ref{theorem:cylindrical_approximation} with $(\pi,\oplus_{i\in I} \mS^i)$, then (\ref{eq:system_direct_sum_state}) and $X^{\pi;\oplus_{i \in I} \mS^i} = \col_{i \in I} X^{\pi;\mS^i}$.
\end{proof}

\begin{remark}
    \label{remark:sample_and_hold_cylindrical}
    If $(\bK^n,\vartheta^n)_{n \geqslant 1}$ is the sample-and-hold labeling scheme of Example~\ref{example:sample_and_hold_labeling_scheme} and $(\Delta_n)_{n \geqslant 1}$ are its sampling grids, Example~\ref{example:sample_and_hold_control} shows that the approximating elementary label controls of Theorem~\ref{theorem:cylindrical_approximation} are sample-and-hold, using the same labels $\mU^n_j$ across $i \in I$. The same is true of the single-control approximation in Corollary~\ref{cor:common_control_parametrized_approximation}. 
\end{remark}

In classical deterministic control theory, the state corresponding to a strict control $u$ is the unique solution to $X^u_t = x_0 + \int_0^t b(s,X^u_s,u_s) \dint s$. Under suitable assumptions, $\eta^{u^n} \rightarrow \eta^\rho$ implies $\|X^{u^n} - X^\rho\|_\infty \rightarrow 0$, where $X^\rho_t = x_0 + \int_0^t \ob(s,X^\rho_s, \rho_s) \dint s$; see Warga~\cite[Section~5]{Warga1967}. We provide an analogous stability property for~(\ref{eq:SDE_cylindrical}), under joint $\sG_T$-stable convergence in law of the labeled Young measures together with the label-generated residual noise.

For this result, choose each $\me_k$ to be Lipschitz; e.g.~apply Gram--Schmidt to the centered $\bRd$-valued polynomials on $[0,1]$ with coefficients in $\bQ^d$.

\begin{proposition}
    \label{prop:cylindrical_stability}
    Let $(\bK^n,\vartheta^n)_{n \geqslant 1}$ be a labeling scheme. Set $\beta^{n,k}_t \coloneq \int_0^t \me_k(\vartheta^n_s)^\top \dint W_s$.
    Let $\Pi = (\pi^i)_{i \in I} \subset \cR^{p,q}$ be a finite family, and suppose $U^n = (u^{n,i})_{i \in I} \subset \cU^{p,q}(\bK^n)$ satisfies
    \begin{equation}
        \label{eq:cylindrical_stability_input}
        \big(\nu^{U^n,\vartheta^n}, (\beta^{n,k})_{k \geqslant 1}\big) \Longrightarrow_{\sG_T} \big(\nu^\Pi, (\beta^k)_{k \geqslant 1}\big)
    \end{equation}
    in $\cY_q^T([0,1] \times \mA^I) \times (C[0,T])^\bN$. Then for any $(\mS^i)_{i \in I} \subset \cS$,
    \begin{equation}
        \label{eq:cylindrical_stability_output}
        \big(\nu^{U^n,\vartheta^n}, (\beta^{n,k})_{k \geqslant 1}, (X^{u^{n,i};\mS^i})_{i \in I} \big) \Longrightarrow_{\sG_T} \big(\nu^\Pi, (\beta^k)_{k \geqslant 1}, (X^{\pi^i;\mS^i})_{i \in I}\big)
    \end{equation}
    in $\cY_q^T([0,1] \times \mA^I) \times (C[0,T])^\bN \times \prod_{i \in I} C([0,T];\bR^{m_{\mS^i}})$.
\end{proposition}

\begin{proof}
    By Remark~\ref{remark:multiple_parametrized_states}, it suffices to consider a single $\pi \in \cR^{p,q}$ and $\mS \in \cS$; suppress $\mS$ from the notation. Suppose (\ref{eq:cylindrical_stability_input}) holds for some $u^n \in \cU^{p,q}(\bK^n)$. Lemma~\ref{lemma:tight_Young_gives_tight_states} and (\ref{eq:cylindrical_stability_input}) show that for every subsequence, there is a further subsequence, not relabeled, such that on some extension of the ground space,
    $$
        \big(\nu^{u^n,\vartheta^n}, (\beta^{n,k})_{k \geqslant 1}, X^{u^n} \big) \Longrightarrow_{\sG_T} \big(\hnu, (\hbeta^k)_{k \geqslant 1}, \hX \big)
    $$
    with $\Law(\hnu, (\hbeta^k)_{k \geqslant 1} \giv \sG_T) = \Law(\nu^\pi, (\beta^k)_{k \geqslant 1} \giv \sG_T)$ by~(\ref{eq:cylindrical_stability_input}). It remains to show that $\Law(\hnu, (\hbeta^k)_{k \geqslant 1}, \hX \giv \sG_T) = \Law(\nu^\pi, (\beta^k)_{k \geqslant 1}, X^\pi \giv \sG_T)$. For each $\ell \in \bN_0 \cup \{\infty\}$, define $\bH^\ell$ by $\sH^\ell_t \coloneq \sG_t \vee \sigma(\hnu_{\leqslant t}, \hX_{\cdot \wedge t}, (\hbeta^k_{\cdot \wedge t})_{k=1}^\ell)$. The Lusin--Suslin theorem applied to the continuous injection $\iota : L^q([0,T];\mA_q) \ni \varphi \mapsto \nu^\varphi \in \cY_q^T([0,1] \times \mA)$ shows $\hnu$ is in its image a.s. Lemma~\ref{lemma:predictable_disintegration} and measurable disintegration in $r$ give predictable action kernels that are Dirac a.e. The Lusin--Suslin theorem applied to the Dirac embedding yields $\hnu = \nu^{\hpi}$ for some $\hpi \in \cR^{p,q}(\bH^0)$, since $\Law(\hnu) = \Law(\nu^\pi)$. For $\ell \in \bN$, apply Proposition~\ref{prop:martingale_characteristics} for the action space $[0,1] \times \mA$, strict controls $(\vartheta^n,u^n)$, and system with state dimension $m+\ell$, initial condition $\col(x_0,0_\ell)$, and coefficients $\tb(t,(x,x'), (r,a)) \coloneq \col(b(t,x,a), 0_\ell)$ and $\tsigma(t,(x,x'),(r,a)) = \col(\sigma(t,x,a), \me_1(r)^\top, \dots, \me_\ell(r)^\top)$; since the $\me_k$ are Lipschitz, this is indeed a well-defined system. The corresponding strict Young measure and state are $\nu^{u^n,\vartheta^n}$ and $\col(X^{u^n}, \beta^{n,1},\dots,\beta^{n,\ell})$. As $(\me_n)_{n \geqslant 1}$ are centered and orthonormal, Proposition~\ref{prop:martingale_characteristics} and L\'evy's characterization give $\bG \hookrightarrow \bH^\ell$ and $(W,\hbeta^1,\dots,\hbeta^\ell)$ is $\bH^\ell$-Brownian. Moreover, with $\hSigma_t = \Sigma(t,\hX_t,\hpi_t)$ and $\hL_t \coloneq \hX_t - x_0 - \int_0^t \ob(s,\hX_s,\rho^{\hpi}_s) \dint s$, we have $\hL \in \cM_{\loc,0}^c(\bH^\ell;\bRm)$ and
    \begin{equation}
        \label{eq:cylindrical_stability_bracket_identities}
        \di [\hL]_t = \hSigma_t \hSigma_t^* \dint t, \qquad \di [\hL,W^i]_t = \hSigma_t \epsilon_i \dint t, \qquad \di [\hL, \hbeta^k]_t = \hSigma_t \me_k \dint t
    \end{equation}
    for $1 \leqslant i \leqslant d$ and $1 \leqslant k \leqslant \ell$. Using level-hitting localizations, the argument in Step~3 of the proof of Proposition~\ref{prop:martingale_characteristics} shows $\bG \hookrightarrow \bH^\infty$, $(W,\hbeta^1,\dots,\hbeta^\ell)$ is $\bH^\infty$-Brownian, by L\'evy's characterization, and $\hL \in \cM_{\loc,0}^c(\bH^\infty; \bRm)$. With $\hbeta^k$ in place of $\beta^k$, (\ref{eq:cylindrical_BM}) constructs a cylindrical $\bH^\infty$-Brownian motion $\hmB$ in $\mH$, lying above $W$, with $\hmB|_{\mH_0}$ independent of~$\sG_T$ since $\Law(\hnu, (\hbeta^k)_{k \geqslant 1} \giv \sG_T) = \Law(\nu^\pi, (\beta^k)_{k \geqslant 1} \giv \sG_T)$. By Kunita--Watanabe and Cauchy--Schwarz, (\ref{eq:cylindrical_stability_bracket_identities}) extends to $\di [\hL, \hmB(h)]_t = \hSigma_t h \dint t$ for $h \in \mH$. By Lemma~\ref{lemma:cylindrical_integral_identification}, $\hL = \hSigma \bigcdot \hB$, so $\hX$ solves~(\ref{eq:SDE_cylindrical}) for $\hpi$ and $\hmB$. Applying the Borel map $\iota^{-1}$ to the preceding conditional-law identity shows $\Law(x_0, W, \hpi, (\hbeta^k)_{k \geqslant 1} \giv \sG_T) = \Law(x_0, W, \pi, (\beta^k)_{k \geqslant 1} \giv \sG_T)$. Strong well-posedness of~(\ref{eq:SDE_cylindrical}) yields $\Law(\hnu, (\hbeta^k)_{k \geqslant 1}, \hX \giv \sG_T) = \Law(\nu^\pi, (\beta^k)_{k \geqslant 1}, X^\pi \giv \sG_T)$, as required.
\end{proof}

\subsection{Representation from pairwise approximation}
\label{sec:4.3}
Work on some filtered extension $(\Omega, \sF_T,\bF,P)$. Fix filtrations $(\bK^n)_{n \geqslant 1}$, each possibly defined on its own extension of the ground space, with $\bG \hookrightarrow \bK^n$. We show that any state assignment indexed by parametrized relaxed controls, satisfying~(\ref{eq:cylindrical_approximation}) even just for families of size two (and joint unlabeled Young measures), is the strong solution map of an equation of the form~(\ref{eq:SDE_cylindrical}).

For $(u,u') \subset \cU(\bK^n)$ and $(\pi,\pi') \subset \cR$, consider the joint Young measures on $[0,T] \times \mA^2$
$$
    \eta^{u,u'}(\di t, \di a, \di a') \coloneq \di t\, \delta_{(u_t,u_t')}(\di a, \di a'), \qquad \eta^{\pi,\pi'}(\di t, \di a, \di a') \coloneq \di t \, \big((\pi_t,\pi_t')_{\#} \dint r\big)(\di a, \di a').
$$
These retain the induced action-pair couplings, but otherwise forget the common label~$r$.

\begin{definition}
    \label{def:parametrized_relaxed_state_assignment}
    A \textit{parametrized relaxed state assignment} is a mapping 
    $$
        (\pi,\mS) \longmapsto Y^{\pi;\mS}
    $$
    on $\cR \times \cS$, where $Y^{\pi;\mS}$ is a continuous $\bF$-adapted $\bR^{m_{\mS}}$-valued process on $(\Omega,\sF_T,\bF,P)$. The assignment has the \textit{pairwise approximation property} relative to $(\bK^n)_{n \geqslant 1}$ if for any $(\pi^1,\mS^1),(\pi^2,\mS^2) \in \cR \times \cS$, not necessarily distinct, there are $u^{n,1}, u^{n,2} \in \cU(\bK^n)$ such that
    $$
        (\eta^{u^{n,1}, u^{n,2}}, X^{u^{n,1};\mS^1}, X^{u^{n,2};\mS^2}) \Longrightarrow_{\sG_T} (\eta^{\pi^1,\pi^2}, Y^{\pi^1;\mS^1}, Y^{\pi^2;\mS^2})
    $$
    in $\cY_2^T(\mA^2) \times C([0,T];\bR^{m_{\mS^1}}) \times C([0,T];\bR^{m_{\mS^2}})$.
\end{definition}

For an assignment $(\pi,\mS) \mapsto Y^{\pi;\mS}$, set
$$
    L^{\pi;\mS}_t \coloneq Y^{\pi;\mS}_t - x_0^\mS - \int_0^t \ob_\mS(s,Y^{\pi;\mS}_s, \rho^\pi_s) \dint s
$$
and abbreviate $\Sigma^{\pi;\mS}_t \coloneq \Sigma_\mS(t,Y^{\pi;\mS}_t,\pi_t)$ and $\osigma^{\pi;\mS}_t = \osigma_\mS(t,Y^{\pi;\mS}_t, \rho^\pi_t)$. Linear growth of $b_\mS, \sigma_\mS$, continuity of $Y^{\pi;\mS}$, and $\pi \in \cR$ show that $L^{\pi;\mS}_t$ is well defined and $\int_0^T \|\Sigma^{\pi;\mS}_s\|_{\HS}^2 \dint s < \infty$ a.s.

Since the pairwise approximation property tests all of $\cR \times \cS$, an additional parametrized relaxed control can encode any event in $\bF$, thereby promoting the martingale characteristics in Proposition~\ref{prop:martingale_characteristics} to the full ambient filtration.

\begin{lemma}
    \label{lemma:ambient_filtration_characteristics}
    Suppose $(\pi,\mS) \mapsto Y^{\pi;\mS}$ has the pairwise approximation property. Then $\bG \hookrightarrow \bF$. For all $(\pi,\cS),(\pi',\cS') \in \cR \times \cS$ and $Z \in \cM_{\loc,0}(\bG)$, $L^{\pi;\mS} \in \cM^c_{\loc,0}(\bF;\bR^{m_\mS})$ and 
    $$
        \di [L^{\pi;\mS}, L^{\pi';\mS'}]_t = \Sigma^{\pi;\mS}_t (\Sigma^{\pi';\mS'}_t)^* \dint t, \qquad \di [L^{\pi;\mS},Z]_t = \osigma^{\pi;\mS}_t \dint [W,Z]_t.
    $$
\end{lemma}

\begin{proof}
    Proposition~\ref{prop:martingale_characteristics} gives $Y_0^{\pi;\mS} = x_0^\mS$, so $L^{\pi;\mS}_0 = 0$. For $\mP^0 \coloneq (1,0,0,0) \in \cS$, the approximation property gives $Y^{\pi;\mP^0} \equiv 0$. Choose $\ma_1 \in \mA \backslash \{\ma_0\}$ and define, for $s < T, F \in \cF_s$,
    \begin{equation}
        \label{eq:pi^F,s}
        \pi^{F,s}_t(\omega;r) \coloneq \begin{cases*}
            \ma_1, &$\omega \in F, \, t > s$, \\
            \ma_0, &otherwise.
        \end{cases*}
    \end{equation}
    Then $\pi^{F,s} \in \cR$. Fix $(\pi,\mS) \in \cR \times \cS$, and apply the pairwise approximation property to $(\pi,\mS), (\pi^{F,s}, \mP^0)$; Proposition~\ref{prop:martingale_characteristics} applied to the corresponding stacked system on $\mA^2$ shows $L^{\pi;\mS}$ is a continuous local martingale in the provided limiting filtration $\bH^{F,s} \subset \bF$, and $\bG \hookrightarrow \bH^{F,s}$. For every $s' > s$, $\eta^{\pi,\pi^{F,s}}((s,s'] \times \mA \times \{\ma_1\}) = (s'-s) 1_F$, and thus $F \in \sH^{F,s}_{s'}$ for all $s' > s$, so that $F \in\sH_s^{F,s}$ by right-continuity. Considering level-hitting localizations of $L^{\pi;\mS}$ and testing against $F \in \sF_s$ as in Step~3 of the proof of Proposition~\ref{prop:martingale_characteristics} shows that $L^{\pi;\mS} \in \cM_{\loc,0}^c(\bF;\bR^{m_\mS})$. The same argument for bounded $\bG$-martingales gives $\bG \hookrightarrow \bF$. 
    The remaining bracket identities follow from a single application of Proposition~\ref{prop:martingale_characteristics} to the stacking of $\mS, \mS'$ and the joint action kernel $(\pi_t,\pi'_t)_\# \dint r$ on $\mA^2$.
\end{proof}

\begin{theorem}
    \label{theorem:cylindrical_converse}
    A parametrized relaxed state assignment $(\pi,\mS) \mapsto Y^{\pi;\mS}$ has the pairwise approximation property if and only if $\bG \hookrightarrow \bF$, and there is a cylindrical $\bF$-Brownian motion $\mB$ in $\mH$, lying above $W$, with $\mB|_{\mH_0}$ independent of $\sG_T$, such that for all $(\pi,\mS) \in \cR \times \cS$, 
    \begin{equation}
        \label{eq:reconstructed_cylindrical_SDE}
        Y^{\pi;\mS}_t = x_0^{\mS} + \int_0^t \ob_{\mS}(s,Y^{\pi;\mS}_s, \rho^\pi_s) \dint s + \int_0^t \Sigma_{\mS}(s,Y^{\pi;\mS}_s, \pi_s) \dint \mB_s.
    \end{equation}
    Moreover, $\mB$ is uniquely determined by the assignment, up to indistinguishability of $\mB(h)$ for each $h \in \mH$. In this case, the assignment has the strict recovery of Proposition~\ref{prop:cylindrical_well_posed}, and the finite-family approximation conclusion of Theorem~\ref{theorem:cylindrical_approximation}. For $\pi \in \cR$ and any finite family $(\mS^i)_{i \in I} \subset \cS$, $Y^{\pi;\oplus_{i \in I} \mS^i} = \col_{i \in I} Y^{\pi;\mS^i}$.
\end{theorem}

\begin{proof}
    Suppose the assignment has the pairwise approximation property. Lemma~~\ref{lemma:ambient_filtration_characteristics} applies; in particular, $\bG \hookrightarrow \bF$. Let $\cB_0 \subset \sB$ be the countable algebra generated by $[0,r]$, $r \in [0,1] \cap \bQ$. Then $\mH_{\step} \coloneq \spn\{1_C \epsilon_i : C \in \cB_0, 1 \leqslant i \leqslant d\}$ is dense in $\mH$. Fix $\ma_1 \in \mA \backslash\{\ma_0\}$ and choose Lipschitz $f : \mA \rightarrow [0,1]$ with $f(\ma_0) = 0$ and $f(\ma_1) = 1$. For $C \in \cB_0$, define 
    \begin{align*}
        \pi^C_t(r) \coloneq \begin{cases*}
            \ma_1, &$r \in C$, \\
            \ma_0, &$r \notin C$.
        \end{cases*}
    \end{align*}
    For $1 \leqslant i \leqslant d$, let $\mP_i \coloneq (1,0,0,\sigma_{\mP_i}) \in \cS$, $\sigma_{\mP_i}(t,x,a) \coloneq f(a) \epsilon_i^\top$. Noting $f(\pi^C_t(r)) = 1_C(r)$, define $\mB(1_C \epsilon_i) \coloneq Y^{\pi^C;\mP_i} = L^{\pi^C;\mP_i}$. By  Lemma~\ref{lemma:ambient_filtration_characteristics}, $\mB(1_C \epsilon_i) \in \cM_{\loc,0}^c(\bF)$ and
    \begin{equation}
        \label{eq:cylindrical_converse_Brownian_identities}
        \di [\mB(1_C \epsilon_i), \mB(1_{C'} \epsilon_j)]_t = \langle 1_C \epsilon_i, 1_{C'} \epsilon_j \rangle \dint t, \qquad \di [\mB(1_C \epsilon_i), Z]_t = |C|\epsilon_i^\top \dint [W,Z]_t
    \end{equation}
    for all $Z \in \cM_{\loc,0}(\bG)$. The first identity in~(\ref{eq:cylindrical_converse_Brownian_identities}) shows we may extend $\mB$ linearly to $\mH_{\step}$, and $E[\sup_{t \leqslant T} |\mB_t(h)|^2] \leqslant 4T\|h\|^2$ for $h \in \mH_{\step}$. So $\mB$ extends uniquely to a continuous linear map $\mB : \mH \rightarrow \cM_0^{2,c}(\bF)$; Kunita--Watanabe gives $\di [\mB(h), \mB(h')]_t = \langle h, h'\rangle \dint t$ for all $h,h' \in \mH$. Thus $\mB$ is a cylindrical $\bF$-Brownian motion in $\mH$. Taking $C = [0,1]$ and $Z = W^i$ in the second identity in~(\ref{eq:cylindrical_converse_Brownian_identities}) gives $[\mB(\epsilon_i) - W^i] \equiv 0$. Hence $\mB$ lies above~$W$. The second identity in~(\ref{eq:cylindrical_converse_Brownian_identities}) extends by linearity and continuity of~$\mB$ to $\di[\mB(h), Z]_t = (\int_0^1 h(r) \dint r)^\top \dint [W,Z]_t$ for every $h \in \mH$ and $Z \in \cM_{\loc,0}(\bG)$. In particular, $[\mB(h),Z] \equiv 0$ for $h \in \mH_0$, and thus Lemma~\ref{lemma:conditional_independence_BM}, applied to every finite orthonormal family in $\mH_0$, gives independence of $\mB|_{\mH_0}$ from $\sG_T$.

    For $(\pi,\mS) \in \cR \times \cS$, applying Lemma~\ref{lemma:ambient_filtration_characteristics} to $(\pi,\mS), (\pi,\mS)$, and then to $(\pi,\mS), (\pi^C,\mP_i)$ for each $C \in \cB_0, 1 \leqslant i \leqslant d$, and extending linearly, shows
    \begin{equation}
        \label{eq:cylindrical_converse_L_identities}
        \di[L^{\pi;\mS}]_t = \Sigma^{\pi;\mS}_t (\Sigma^{\pi;\mS}_t)^* \dint t, \qquad \di [L^{\pi;\mS},\mB(h)]_t = \Sigma^{\pi;\mS}_t h \dint t
    \end{equation}
    for $h \in \mH_{\step}$. Both sides of the latter identity in~(\ref{eq:cylindrical_converse_L_identities}) are continuous in $h$---the left-hand side by Kunita--Watanabe, and the right-hand side by Cauchy--Schwarz; thus~(\ref{eq:cylindrical_converse_L_identities}) holds for all $h \in \mH$. Hence
    Lemma~\ref{lemma:cylindrical_integral_identification} with $Q = I_\mH$ gives~(\ref{eq:reconstructed_cylindrical_SDE}), and Proposition~\ref{prop:cylindrical_well_posed} identifies $Y^{\pi;\mS}$ as its unique strong solution. Any other driver satisfying~(\ref{eq:reconstructed_cylindrical_SDE}) must respect $\mB(1_C\epsilon_i) = Y^{\pi^C;\mP_i}$, which determines $\mB$ on $\mH$ by linearity and continuity; thus uniqueness of the driver holds.

    Conversely, suppose such $\mB$ exists and the assignment satisfies~(\ref{eq:reconstructed_cylindrical_SDE}). Proposition~\ref{prop:cylindrical_well_posed} identifies the assigned states with its unique strong solutions, and Theorem~\ref{theorem:cylindrical_approximation}, applied with $p=1, q = 2, |I| = 2$, and the labeling scheme $(\bK^n,\theta_n)_{n \geqslant 1}$ of Example~\ref{example:standard_labeling_scheme} gives the pairwise approximation property. Moreover, the full conclusion of Theorem~\ref{theorem:cylindrical_approximation} holds. The final assertion follows by pathwise uniqueness and since $\col_{i \in I} Y^{\pi;\mS^i}$ solves~(\ref{eq:reconstructed_cylindrical_SDE}) for $(\pi,\oplus_{i \in I} \mS^i)$.
\end{proof}

Theorems~\ref{theorem:cylindrical_approximation} and~\ref{theorem:cylindrical_converse} characterize the parametrized relaxed state assignments with the pairwise approximation property as the strong solution maps of~(\ref{eq:SDE_cylindrical}), with each assignment determining its cylindrical Brownian driver uniquely. We next consider state assignments indexed by ordinary relaxed controls.

\section{Ordinary relaxed controls and universal coherence}
\label{sec:5}
We determine the structure of a state assignment indexed by ordinary relaxed controls under a finite-family approximation property. Pairwise approximation recovers action-pair couplings from the state covariations. These couplings determine a common residual covariance structure on the relevant relaxed control class, which we first realize by a cylindrical martingale under a separability-type assumption. Then, finite-family approximation shows that these relaxed controls come from coherent parametrizations, and that the underlying dynamics takes the form~(\ref{eq:SDE_cylindrical}), up to a predictable orthogonal projection of the driving space onto the used noise directions.

\subsection{Recovery of pairwise action couplings}

Throughout the section, state assignments are defined on a fixed arbitrary filtered extension $(\Omega, \sF_T,\bF,P)$ of the ground stochastic basis. We also fix filtrations $(\bK^n)_{n \geqslant 1}$, each possibly defined on its own extension, such that $\bG \hookrightarrow \bK^n$ for every $n$. Sometimes we write $z = (\omega,t) \in \oOmega$.

\begin{definition}
    \label{def:relaxed_state_assignment}
    A \textit{relaxed state assignment} is a mapping 
    $$
        (\rho,\mS) \longmapsto Y^{\rho;\mS}
    $$ 
    on $\cQ \times \cS$, where $Y^{\rho;\mS}$ is a continuous $\bF$-adapted $\bR^{m_\mS}$-valued process on $\Omega$, compatible with direct sums of systems, in the sense that $Y^{\rho;\mS \oplus \mS'} = \col(Y^{\rho;\mS}, Y^{\rho;\mS'})$ for all $\rho \in \cQ$, $\mS,\mS' \in \cS$.
\end{definition}

\begin{definition}
    \label{def:pair_approximation}
    The assignment has the \textit{pairwise marginal approximation property} if, for every $(\rho^1,\mS^1), (\rho^2,\mS^2) \in \cQ \times \cS$, not necessarily distinct, there are $u^{n,1},u^{n,2} \in \cU(\bK^n)$ such that
    $$
        (\eta^{u^{n,1}}, \eta^{u^{n,2}}, X^{u^{n,1}; \mS^1}, X^{u^{n,2}; \mS^2}) \Longrightarrow_{\sG_T} (\eta^{\rho^1}, \eta^{\rho^2}, Y^{\rho^1; \mS^1}, Y^{\rho^2;\mS^2})
    $$
    in $\cY_2^T(\mA)^2 \times C([0,T];\bR^{m_{\mS^1}}) \times C([0,T];\bR^{m_{\mS^2}})$.
\end{definition}

Direct-sum compatibility is explicitly assumed in the present setting, unlike in Section~4, where it is deduced by Theorem~\ref{theorem:cylindrical_converse}. This assumption permits several systems to be tested using one approximating strict control per relaxed control, ensuring the pairwise couplings recovered from Definition~\ref{def:pair_approximation} in Proposition~\ref{prop:gamma^rho,rho'} below depend on the relaxed controls alone, and not also on the choice of system. 

For the remainder of the section, fix a relaxed state assignment $(\rho,\mS) \mapsto Y^{\rho;\mS}$ with the pairwise marginal approximation property. For $(\rho,\mS) \in \cQ \times \cS$, set 
$$
    L^{\rho;\mS}_t = Y^{\rho;\mS}_t - x^\mS_0 - \int_0^t \ob_\mS(s,Y^{\rho;\mS}_s, \rho_s) \dint s.
$$
The drift integral is well defined by linear growth, continuity of $Y^{\rho;\mS}$, and $\rho \in \cQ$. The next lemma recovers a joint-action kernel from joint control-state convergence. We will also use this result in Section~5.3, but in the present subsection, we take $|I| = 2$. For finite $I$ and $(\mu^i)_{i \in I} \subset \cP(\mA)$, let $\Gamma((\mu^i)_{i \in I}) \coloneq \{ \lambda \in \cP(\mA^I) : (\pr^I_i)_\# \lambda = \mu^i \text{ for every } i \in I\}$. 

\begin{lemma}
    \label{lemma:common_kernel}
    Let $I$ be finite and $J$ countable. For $i \in I, j \in J$, let $\rho^i \in \cQ$, $\mS^{i,j} \in \cS$, and suppose there are $u^{n,i} \in \cU(\bK^n), i \in I$, such that
    \begin{equation}
        \label{eq:common_kernel_convergence}
        \big((\eta^{u^{n,i}})_{i \in I}, (X^{u^{n,i};\mS^{i,j}})_{i \in I, j \in J}\big) \Longrightarrow_{\sG_T} \big((\eta^{\rho^i})_{i \in I}, (Y^{\rho^i;\mS^{i,j}})_{i \in I, j \in J}\big)
    \end{equation}
    in $\prod_{i \in I} \cY^T_2(\mA) \times \prod_{i \in I, j \in J} C([0,T]; \bR^{m_{\mS^{i,j}}})$. Set $\sH_t = \sG_t \vee \sigma(\eta^{\rho^i}_{\leqslant t}, \,Y^{\rho^i;\mS^{i,j}}_{\cdot \wedge t} : i \in I, j \in J)$. Then $\bG \hookrightarrow \bH$, each $\rho^i$ has an $\bH$-predictable representative, and $L^{\rho^i;\mS^{i,j}} \in \cM_{\loc,0}^c(\bH;\bR^{m_{\mS^{i,j}}})$. Moreover, there is an $\bH$-predictable probability kernel $\lambda$ on $\mA^I$, such that $\lambda_z \in \Gamma((\rho_z^i)_{i \in I})$ $\oP$-a.e.\ and, for all $(i,j),(i',j') \in I \times J$,
    \begin{equation}
        \label{eq:common_kernel_equality}
        \di [L^{\rho^i;\mS^{i,j}}, L^{\rho^{i'};\mS^{i',j'}}]_t = \bigg( \int_{\mA^I} \sigma_{\mS^{i,j}}(t,Y^{\rho^i;\mS^{i,j}}_t, a_i) \sigma_{\mS^{i',j'}}(t,Y^{\rho^{i'};\mS^{i',j'}}_t, a_{i'})^\top \lambda_t(\di a) \bigg) \dint t.
    \end{equation}
\end{lemma}

\begin{proof}
    \textit{Step 1. Joint limit on an extension.}
    Set $Y^{i,j} \coloneq Y^{\rho^i;\mS^{i,j}}$ and $L^{i,j} \coloneq L^{\rho^i;\mS^{i,j}}$. By Lemma~\ref{lemma:predictable_disintegration}, choose $\rho^i$ predictable for the filtration given by $\sG_t \vee \sigma(\eta^{\rho^i}_{\leqslant t})$. Set $U^n \coloneq (u^{n,i})_{i \in I}$ and $\eta^{n,I}(\di t, \di a) \coloneq \di t \, \delta_{U^n_t}(\di a)$ on $[0,T] \times \mA^I$. Let $\pr^I_i(a) \coloneq a_i$ for $a \in \mA^I$. For $\eta \in \cY^T_2(\mA^I)$, let $\fm_i(\eta) \coloneq (\id_{[0,T]}, \pr^I_i)_\# \eta$ be its $i$th action marginal.
   
    We claim the laws of $\eta^{n,I}$ are tight in $\cY^T_2(\mA^I)$. By~(\ref{eq:common_kernel_convergence}), the $\eta^{u^{n,i}}$ are tight in $\cY^T_2(\mA)$, so there are compact $K_i \subset \cY^T_2(\mA)$ with $\sup_{n \geqslant 1} P[\eta^{u^{n,i}} \notin K_i] < \varepsilon/|I|$ for fixed $\varepsilon > 0$; thus $\sup_{n \geqslant 1} P[\eta^{n,I} \notin K] < \varepsilon$ for $K \coloneq \{\eta \in \cY^T_2(\mA^I) : \fm_i(\eta) \in K_i, \, i \in I\}$. It suffices to show $K \subset \cY^T_2(\mA^I)$ is compact. Any sequence in $K$ has a subsequence with convergent action marginals in $\cY^T_2(\mA)$. Hence, the subsequence is tight on $[0,T] \times \mA^I$. Prokhorov's theorem provides a narrowly convergent further subsequence. Its second moments converge by $\int_{[0,T] \times \mA^I} t^2 + d_I(a,\ma_0^I)^2 \, \eta(\di t, \di a) = T^3/3 + \sum_{i \in I} \int_{[0,T] \times \mA} d(a,\ma_0)^2 \, \fm_i(\eta)(\di t, \di a)$. Hence convergence holds in $\cY^T_2(\mA^I)$, with the limit remaining in~$K$. Thus $K$ is compact, as claimed.

    Let $V^n,V$ denote the left- and right-hand sides of~(\ref{eq:common_kernel_convergence}), respectively. The joint laws of $(\eta^{n,I},V^n)$ are tight, so there is a subsequence, which we do not relabel, converging $\sG_T$-stably in law. The control-state marginal of the limit has the same $\sG_T$-law as $V$. Applying~\cite[Theorem~8.17]{kallenberg2021foundations} jointly with the projection onto $(\Omega,\sG_T,P)$ realizes the joint limit on a further probability extension $\mathfrak{r} : (\hOmega, \hsF_T,\hP) \rightarrow (\Omega,\sF_T,P)$,
    \begin{equation}
        \label{eq:common_kernel_extended_convergence}
        (\eta^{n,I}, V^n) \Longrightarrow_{\sG_T} (\heta^I, V \circ \fr).
    \end{equation} 
    Hence $\fm_i(\heta^I) = \eta^{\rho^i} \circ \fr$ $\hP$-a.s. By Lemma~\ref{lemma:predictable_disintegration}, disintegrate $\heta^I(\di t, \di a) = \di t \, \hlambda_t(\di a)$ with $\hlambda$ predictable for the filtration generated by $\hbG$, $\heta^{I}_{\leqslant}$, where $\hsG_t \coloneq \fr^{-1}(\sG_t)$. With $\ofr(\homega,t) \coloneq (\fr(\homega),t)$, uniqueness of disintegration gives $(\pr^I_i)_\# \hlambda = \rho^i \circ \ofr$ $\hP \otimes \di t$-a.e. 
    
    \textit{Step 2. Martingale properties on original space.} Choose finite $J_\ell \subset J$ with $J_\ell \uparrow J$. For each $\ell \geqslant 1$, set $\mS^{i,J_\ell} \coloneq \oplus_{j \in J_\ell} S^{i,j}, i \in I$, and let $\mS^{\ell,\stk}$ be their stacked system on $\mA^I$, with diffusion coefficient denoted $\sigma^\ell$. By (\ref{eq:system_stacked_state})--(\ref{eq:system_direct_sum_state}), $X^{U^n;\mS^{\ell,\stk}} = \col_{i \in I, j \in J_\ell} X^{u^{n,i};\mS^{i,j}}$; thus with $Y^\ell \coloneq \col_{i \in I, j \in J_\ell} Y^{i,j}$, we have $(\eta^{U^n}, X^{U^n;\mS^{\ell,\stk}}) \Rightarrow_{\sG_T} (\heta^I, Y^\ell \circ \fr)$ by~(\ref{eq:common_kernel_extended_convergence}). Proposition~\ref{prop:martingale_characteristics} gives $\hbG \hookrightarrow \hbH^\ell$, where $\hsH^\ell_t \coloneq \hsG_t \vee \sigma(\heta_{\leqslant t}^I, (Y^{i,j} \circ \fr)_{\cdot \wedge t} : i \in I, j \in J_\ell)$, and moreover, with $L^\ell \coloneq \col_{i \in I, j \in J_\ell} L^{i,j}$, $L^\ell \circ \fr$ is a continuous local $\hbH^\ell$-martingale satisfying
    \begin{equation}
        \label{eq:common_kernel_equality_extension}
        \di [L^\ell \circ \fr]_t = \bigg( \int_{\mA^I} F^\ell_t(a) \circ \fr \, \hlambda_t(\di a) \bigg) \dint t,
    \end{equation}
    where $F^\ell_t(a) \coloneq \sigma^\ell(t,Y_t^\ell,a) \sigma^\ell(t,Y^\ell_t,a)^\top$. Now fix $(i,j)\in I \times J$ and consider, for $R > 0$, $\tau_R \coloneq \inf\{ t \geqslant 0 : |L^{i,j}_t| \geqslant R\} \wedge T$. Define $\bH^\ell$ by $\sH^\ell_t \coloneq \sG_t \vee \sigma(\eta^{\rho^i}_{\leqslant t}, Y^{i,j}_{\cdot \wedge t}: i \in I, j \in J_\ell)$. Note $\fr^{-1}(\sH^\ell_t) \subset \hsH_t^\ell$ because $\fm_i(\heta^I) = \eta^{\rho^i} \circ \fr$. For $\ell$ large enough with $j \in J_\ell$, $(L^{i,j})^{\tau_R} \circ \fr$ is a bounded $\hbH^\ell$-martingale, and for any $H \in L^\infty(\Omega,\sH^\ell_s, P)$, $H \circ \fr$ is $\hsH_s^\ell$-measurable. So for $t \geqslant s$, $E[H(L^{i,j}_{t \wedge \tau_R} - L^{i,j}_{s \wedge \tau_R})] = \hE[(H\circ \fr)(L^{i,j}_{t \wedge \tau_R} - L^{i,j}_{s \wedge \tau_R})\circ \fr] = 0$. The argument in Step~3 of the proof of Proposition~\ref{prop:martingale_characteristics} shows $L^{i,j} \in \cM_{\loc,0}^c(\bH)$. Moreover, $[L^\ell] \circ \fr = [L^\ell \circ \fr]$. The same argument without stopping, applied to a bounded \cadlag $\bG$-martingale, proves $\bG \hookrightarrow \bH$.    
    
    \textit{Step 3. Predictable kernel on original space.} Define $\mu \coloneq (\ofr, \id_{\mA^I})_\#((\hP \otimes \di t)(\di \hz)\, \hlambda_{\hz}(\di a))$ on $(\oOmega \times \mA^I, \sP^{\bH} \otimes \sB(\mA^I))$. Its $\oOmega$-marginal is $\oP$, and since $\mA^I$ is Polish, disintegration in~\cite[Theorem~3.4]{kallenberg2021foundations} gives $\mu(\di z, \di a) = \oP(\di z) \lambda_z(\di a)$ for an $\bH$-predictable $\lambda$, and 
    \begin{equation}
        \label{eq:common_kernel_disintegration_identity}
        \int_{\oOmega} \int_{\mA^I} g(z,a) \, \lambda_z(\di a) \, \oP(\di z) = \int_{\hOmega \times [0,T]} \int_{\mA^I} g\big(\ofr(\hz), a\big) \, \hlambda_{\hz}(\di a) \, (\hP \otimes \di t)(\di \hz)
    \end{equation}
    for all integrable $\sP^\bH \otimes \sB(\mA^I)$-measurable $g$. As $(\pr^I_i)_\# \hlambda = \rho^i \circ \ofr$, the $(\oOmega, a_i)$-marginal of $\mu$ is $\oP(\di z) \, \rho^i_z(\di a_i)$. Thus $(\pr^I_i)_\# \lambda = \rho^i$ for $i \in I$. Hence $\lambda_z \in \Gamma((\rho^i_z)_{i \in I})$ for $\oP$-a.e.~$z$ and $\int_{\mA^I} d_I(a,\ma_0^I)^2 \, \lambda_z(\di a) = \sum_{i\in I} \int_\mA d(a,\ma_0)^2 \, \rho^i_z(\di a) < \infty$. For 
    $$
        C^\ell_t \coloneq \int_0^t \int_{\mA^I} F^\ell_s(a) \, \lambda_s(\di a) \dint s,
    $$
    $\tr C^\ell_t = \sum_{i \in I} \sum_{j \in J_\ell} \int_0^t \int_\mA |\sigma_{\mS^{i,j}}(s,Y^{i,j}_s,a)|^2 \, \rho^i_s(\di a) \dint s = \tr [L^\ell]_t$ by (\ref{eq:common_kernel_equality_extension}). Thus $C^\ell$ is a well-defined continuous $\bH$-adapted process of finite variation, and moreover, its increments are positive semidefinite. Thus with $\tau^\ell_k \coloneq \inf\{t \geqslant 0 : \tr [L^\ell]_t \geqslant k\} \wedge T$, both $[L^\ell]^{\tau^\ell_k}$ and $(C^\ell)^{\tau^\ell_k}$ have total variation at most $k$. For bounded $\bH$-predictable $H$, (\ref{eq:common_kernel_equality_extension})--(\ref{eq:common_kernel_disintegration_identity}) give
    \begin{align*}
        E\bigg[\int_0^{\tau^\ell_k} H_s \dint [L^\ell]_s \bigg] =
        \hE\bigg[ \int_0^{\tau^\ell_k \circ \fr} (H_s \circ \fr) \int_{\mA^I} F^\ell_s(a) \circ \fr \, \hlambda_s(\di a) \dint s\bigg] = E\bigg[ \int_0^{\tau^\ell_k} H_s \dint C^\ell_s \bigg].
    \end{align*}
    Taking $H = 1_G 1_{(s,t]}$ for $G \in \sH_s, s \leqslant t$, shows that $([L^\ell] - C^\ell)^{\tau^\ell_k}$ is an $\bH$-martingale. As it is also a continuous process of finite variation, null at zero, it vanishes, and thus $[L^\ell] = C^\ell$. As any $(i,j),(i',j') \in I \times J$ belong to $I \times J_\ell$ for large enough $\ell$, this completes the proof.
\end{proof}

By the pairwise marginal approximation property, we may apply Lemma~\ref{lemma:common_kernel} with $|I| = 2$ and $|J|=1$ to the fixed state assignment; (\ref{eq:common_kernel_equality}) then permits us to use the arguments of Lemma~\ref{lemma:ambient_filtration_characteristics} to promote the martingale characteristics of Proposition~\ref{prop:martingale_characteristics} to all of $\bF$.

\begin{lemma}
    \label{lemma:ambient_filtration_characteristics_ordinary}
    We have $\bG \hookrightarrow \bF$, and for every $(\rho,\mS) \in \cQ \times \cS$, $L^{\rho;\mS} \in \cM_{\loc,0}^c(\bF;\bR^{m_\mS})$,
    $$
        \di [L^{\rho;\mS}]_t = \bigg( \int_\mA \sigma_\mS(t,Y^{\rho;\mS}_t, a) \sigma_\mS(t,Y^{\rho;\mS}_t, a)^\top \, \rho_t(\di a) \bigg) \dint t,
    $$
    and $\di [L^{\rho;\mS},Z]_t = \osigma_\mS(t,Y^{\rho;\mS}_t,\rho_t) \dint [W,Z]_t$ for all $Z \in \cM_{\loc,0}(\bG)$.
\end{lemma}

\begin{proof}
    Let $\pi^{F,s}$ be as in~(\ref{eq:pi^F,s}), with $\rho^{F,s}$ its induced relaxed control. Definition~\ref{def:pair_approximation} and Lemma~\ref{lemma:common_kernel} applied to $(\rho,\mS), (\rho^{F,s}, \mP^0)$ generate $\bH^{F,s} \subset \bF$ and yield $\bG \hookrightarrow \bH^{F,s}$ and $L^{\rho;\mS} \in \cM_{\loc,0}^c(\bH^{F,s};\bR^{m_\mS})$. As $\eta^{\rho^{F,s}}((s,s'] \times \{\ma_1\}) = (s'-s) 1_F$ for $s' > s$, right-continuity of $\bH^{F,s}$ gives $F \in \sH^{F,s}_s$. Continuing as in Lemma~\ref{lemma:ambient_filtration_characteristics} completes the proof. 
\end{proof}

Given a pair $(\rho,\rho') \subset \cQ$, we find a kernel $\gamma^{\rho,\rho'}$ on $\mA^2$ for which~(\ref{eq:common_kernel_equality}) holds for all $\mS,\mS' \in \cS$. Denote by $\BL(\mA)$ and $\BL(\mA;\bRd)$ the bounded Lipschitz functions on $\mA$ into $\bR$ and $\bRd$. Fix a countable unital $\bQ$-algebra $\cA_0 \subset \BL(\mA)$ with $\sigma(\cA_0) = \sB(\mA)$, e.g.~the $\bQ$-algebra generated by 1 and $a \mapsto (1-d(a,a_k)/s)^+$ for $k \geqslant 1, s \in \bQ_{>0}$, where $(a_k)_{k \geqslant 1} \subset \mA$ is dense. Set
$$
    \cD \coloneq \spn_\bQ\{g \epsilon_i : g \in \cA_0, \, 1 \leqslant i \leqslant d\} \subset \BL(\mA;\bRd).
$$
The functional monotone-class theorem shows $\cA_0$ determines probability measures on $\mA$, $\cA_0 \otimes \cA_0 \coloneq \{(a,a') \mapsto g(a)g'(a') : g,g' \in \cA_0\}$ determines probability measures on $\mA^2$, and $\cD$ is dense in $L^2(\mA,\mu;\bRd)$ for all $\mu \in \cP(\mA)$. For $f \in \BL(\mA;\bRd)$ define the \textit{probe systems}
$$
    \mP^f \coloneq (1,0,0,\sigma_{\mP^f}) \in \cS, \quad \sigma_{\mP^f}(t,x,a) \coloneq f(a)^\top.
$$
Note $\mP^0 = (1,0,0,0)$. Define $\Delta : a \mapsto (a,a)$ on $\mA$ and $\sw : (a,a') \mapsto (a',a)$ on $\mA^2$.

\begin{proposition}
    \label{prop:gamma^rho,rho'}
    For every $\rho,\rho' \in \cQ$, there is a unique, up to $\oP$-a.e.~equality, $\bF$-predictable kernel $\gamma^{\rho,\rho'}$ with $\gamma^{\rho,\rho'}_z \in \Gamma(\rho_z,\rho'_z)$ $\oP$-a.e., such that, for every $\mS,\mS' \in \cS$, 
    \begin{equation}
        \label{eq:pairwise_action_coupling_bracket}
        \di [L^{\rho;\mS}, L^{\rho';\mS'}]_t = \bigg(\int_{\mA^2} \sigma_\mS(t,Y^{\rho;\mS}_t, a)\sigma_{\mS'}(t,Y^{\rho';\mS'}_t, a')^\top \, \gamma^{\rho,\rho'}_t(\di a, \di a')\bigg) \dint t.
    \end{equation}
    With $\sH^{\rho,\rho';\cD}_t = \sG_t \vee \sigma(\eta^{\rho}_{\leqslant t}, \eta^{\rho'}_{\leqslant t}, Y^{\rho;\mP^f}_{\cdot \wedge t}, Y^{\rho';\mP^f}_{\cdot \wedge t} : f \in \cD)$, $\gamma^{\rho,\rho'}$ has an $\bH^{\rho,\rho';\cD}$-predictable version. Moreover, $\gamma^{\rho,\rho} = \Delta_\# \rho$ and $\gamma^{\rho',\rho} = \sw_\# \gamma^{\rho,\rho'}$. 
\end{proposition}

\begin{proof}
    Enumerate $\cD = (f_k)_{k \geqslant 1}$, and apply Definition~\ref{def:pair_approximation} to $(\rho,\oplus_{\ell\leqslant k} \mP^{f_\ell}), (\rho', \oplus_{\ell \leqslant k} \mP^{f_\ell})$ for each $k \geqslant 1$. Padding the corresponding finite tuples by zero paths in a bounded countable-product metric and using direct-sum compatibility and Lemma~\ref{lemma:stable_convergence_facts}(iii), and then Lemma~\ref{lemma:stable_convergence_facts}(ii) to remove the padding, yields $u^{n,1}, u^{n,2} \in \cU(\bK^n)$ with
    $$
        \label{eq:gamma^rho,rho'_convergence}
        \big(\eta^{u^{n,1}}, \eta^{u^{n,2}}, (X^{u^{n,1};\mP^{f_\ell}})_{\ell \geqslant 1}, (X^{u^{n,2};\mP^{f_\ell}})_{\ell \geqslant 1} \big) \Longrightarrow_{\sG_T} \big( \eta^{\rho}, \eta^{\rho'}, (Y^{\rho;\mP^{f_\ell}})_{\ell \geqslant 1}, (Y^{\rho';\mP^{f_\ell}})_{\ell \geqslant 1} \big).
    $$
    Lemma~\ref{lemma:common_kernel} gives an $\bH^{\rho,\rho';\cD}$-predictable $\gamma^{\rho,\rho'}$ satisfying $\gamma^{\rho,\rho'}_z \in \Gamma(\rho_z,\rho_z')$ $\oP$-a.e.~and~(\ref{eq:pairwise_action_coupling_bracket}) for $(\rho, \mP^f),(\rho',\mP^{f'})$, for all $f,f' \in \cD$. With $f = g\epsilon_1, f' = g'\epsilon_1$, the determining property of $\cA_0\otimes \cA_0$ shows $\gamma^{\rho,\rho'}$ is the unique $\bF$-predictable kernel with this property. Fix $\mS,\mS' \in \cS$, and apply Definition~\ref{def:pair_approximation} to $(\rho,\oplus_{\ell=1}^k \mP^{f_\ell} \oplus \mS), (\rho', \oplus_{\ell=1}^k \mP^{f_\ell} \oplus \mS')$. The same argument gives an $\bF$-predictable $\tgamma$ satisfying~(\ref{eq:pairwise_action_coupling_bracket}) and the same probe identities as $\gamma^{\rho,\rho'}$; by uniqueness, $\tgamma = \gamma^{\rho,\rho'}$, and hence~(\ref{eq:pairwise_action_coupling_bracket}) holds for $\gamma^{\rho,\rho'}$. For $f,f' \in \cD$, by direct-sum compatibility and Lemma~\ref{lemma:ambient_filtration_characteristics_ordinary} applied to $\mP^f \oplus \mP^{f'}$, $\di [Y^{\rho;\mP^f}, Y^{\rho;\mP^{f'}}]_t = (\int_\mA f(a)^\top f'(a) \, \rho_t(\di a)) \dint t$. Thus (\ref{eq:pairwise_action_coupling_bracket}) and the determining property of $\cA_0 \otimes \cA_0$ yield $\gamma^{\rho,\rho} = \Delta_\# \rho$. The identity $\gamma^{\rho',\rho} = \sw_\# \gamma^{\rho,\rho'}$ follows similarly, using symmetry of quadratic covariation.
\end{proof}

\begin{remark}
    As $[Y^{\rho;\mS}, Y^{\rho';\mS'}] = [L^{\rho;\mS}, L^{\rho';\mS'}]$, (\ref{eq:pairwise_action_coupling_bracket}) is the analogue of (\ref{eq:action_coupling_state_characteristics}) for $|I| = 2$.
\end{remark}

\subsection{The minimal common noise}

We establish a cylindrical martingale representation of the relaxed state assignment under Definition~\ref{def:pair_approximation}, analogously to Theorem~\ref{theorem:cylindrical_converse}.

In Section~\ref{sec:4.3}, for $\pi,\pi' \in \cR$ and $f, f' \in \BL(\mA;\bRd$), with $N^{\pi,f}_t \coloneq Y^{\pi;\mP^f}_t - \int_0^t \rho^\pi_s(f)^\top \dint W_s$, 
\begin{equation}
    \label{eq:N^pi_characteristics}
    \di [N^{\pi,f}, N^{\pi',f'}]_t = \bigg( \int_0^1 f\big( \pi_t(r) \big)^\top f'\big(\pi'_t(r) \big) \dint r - \rho^\pi_t(f)^\top \rho^{\pi'}_t(f') \bigg) \dint t
\end{equation}
by Lemma~\ref{lemma:ambient_filtration_characteristics}. With $\alpha_t^{\pi,f} \coloneq f \circ \pi_t - \rho^\pi_t(f) \in \mH_0$, (\ref{eq:N^pi_characteristics}) can be equivalently written as
\begin{equation}
    \label{eq:N^pi_characteristics_inner_product}
    \di [N^{\pi,f}, N^{\pi',f'}]_t = \langle \alpha^{\pi,f}_t, \alpha^{\pi',f'}_t \rangle_{\mH_0} \dint t.
\end{equation}
This is used implicitly in Theorem~\ref{theorem:cylindrical_converse}. Indeed, a particular choice of $\pi$ and $f$ recovered the residual driving directions, with the reconstructed driver being cylindrical Brownian by~(\ref{eq:N^pi_characteristics_inner_product}). 

Presently, for $\rho,\rho' \in \cQ$ and $f,f' \in \BL(\mA;\bRd)$, with $N^{\rho,f}_t \coloneq Y^{\rho;\mP^f}_t - \int_0^t \rho_s(f)^\top \dint W_s$,
$$
    \di [N^{\rho,f}, N^{\rho',f'}]_t = \bigg( \int_{\mA^2} f(a)^\top f'(a') \, \gamma_t^{\rho,\rho'}(\di a, \di a') - \rho_t(f)^\top \rho'_t(f') \bigg) \dint t \eqqcolon \kappa^{\rho,\rho'}_t(f,f') \dint t
$$
by Proposition~\ref{prop:gamma^rho,rho'}, and $N^{\rho,f} \in \cM_{\loc,0}^c(\bF)$ and $[N^{\rho,f}, Z] \equiv 0$ for all $Z \in \cM_{\loc,0}(\bG)$ by Lemma~\ref{lemma:ambient_filtration_characteristics_ordinary}. In particular, $[N^{\rho,f}]_t = \int_0^t \rho_s(|f|^2) - |\rho_s(f)|^2 \dint s \leqslant t\|f\|_\infty^2$, so $N^{\rho,f} \in \cM^{2,c}_0(\bF)$. Furthermore, for any $(\rho^i,f_i)_{i=1}^n \subset \cQ \times \BL(\mA;\bRd)$, the matrix $(\kappa^{\rho^i,\rho^j}_z(f_i,f_j))_{i,j=1}^n$ is positive semidefinite for $\oP$-a.e.~$z \in \oOmega$, being the covariance density of $(N^{\rho^i,f_i})_{i=1}^n$.

In view of~(\ref{eq:N^pi_characteristics_inner_product}), we introduce formal elements $[\rho,f]_z$, and aim to specify their \say{inner products} $z$-fiberwise by $\kappa_z$. The family of Hilbert spaces to which the $[\rho,f]_z$ belong should form a Hilbert field. Thus, pass to a nonempty subclass $\cC \subset \cQ$ and work in $\bH^\cC$, where
$$
    \sH^\cC_t \coloneq \sG_t \vee \sigma(\eta^{\rho}_{\leqslant t}, \, Y^{\rho;\mS}_{\cdot \wedge t} : \rho \in \cC, \, \mS \in \cS).
$$
As $\bG \subset \bH^\cC \subset \bF$, Lemma~\ref{lemma:ambient_filtration_characteristics_ordinary} gives $\bG \hookrightarrow \bH^\cC$. Each $\rho \in \cC$ has an $\bH^\cC$-predictable version by Lemma~\ref{lemma:predictable_disintegration}, and for $\rho,\rho' \in \cC$, the $\bH^{\rho,\rho';\cD}$-predictable version of $\gamma^{\rho,\rho'}$ given by Proposition~\ref{prop:gamma^rho,rho'} is also $\bH^\cC$-predictable. These versions are used throughout, and predictability refers to $\bH^\cC$. 

\begin{definition}
    \label{def:countable_residual_multiplicity}
    The assignment has \textit{countable residual multiplicity} on $\cC$ if 
    \begin{longlist}
        \item there are predictable representatives of $\kappa^{\rho,\rho'}(f,f')$ for $\rho,\rho' \in \cC, f,f' \in \BL(\mA;\bRd)$, such that, outside one predictable $\oP$-null set, the prescription 
        \begin{equation}
            \label{eq:form}
            \langle [\rho,f]_z, [\rho',f']_z\rangle_z \coloneq \kappa^{\rho,\rho'}_z(f,f'),
        \end{equation}
        is a symmetric positive-semidefinite bilinear form on $\spn\{[\rho,f]_z : \rho \in \cC, f \in \BL(\mA;\bRd)\}$. Write $K_z$ for the completion of the corresponding  quotient by the zero-norm subspace;
        \item there is a sequence $(\rho^n, f_n)_{n \geqslant 1} \subset \cC \times \BL(\mA;\bRd)$ with
        \begin{equation}
            \label{eq:totality_of_fundamental_sequence}
            K_z = \ospn\{[\rho^n,f_n]_z : n \geqslant 1\} \quad \oP\text{-a.e.}
        \end{equation}
    \end{longlist}
    Set $K_z \coloneq \{0\}$ on a common predictable $\oP$-null set outside of which (i)--(ii) hold everywhere. 
    We call 
    $\cK \coloneq (K_z)_{z \in \oOmega}$ the \textit{minimal residual covariance field} of the assignment on $\cC$.
\end{definition}

\begin{remark}
    \label{remark:countable_residual_multiplicity_automatic}
    For $\rho \in \cC, f,f' \in \BL(\mA;\bRd)$, by (\ref{eq:form}) and construction of $\kappa^{\rho,\rho}_z(f,f')$, 
    \begin{equation}
        \label{eq:bracket_continuity}
        \|[\rho,f]_z - [\rho,f']_z\|_z^2 = \rho_z(|f-f'|^2) - |\rho_z(f-f')|^2 \leqslant \rho_z(|f-f'|^2).
    \end{equation}
    Thus, by density of $\cD$ in $L^2(\mA,\rho_z;\bRd)$, we may assume $f_n \in \cD$ for all $n$ in (ii). In particular, Definition~\ref{def:countable_residual_multiplicity} is automatic when $\cC$ is countable. Indeed, we may take $(\rho^n,f_n)_{n \geqslant 1} = \cC \times \cD$ in (ii); intersecting over $n \in \bN$ the full-measure sets on which $(\kappa^{\rho^i,\rho^j}_z(f_i,f_j))_{i,j=1}^n$ is positive semidefinite yields (i), by~(\ref{eq:bracket_continuity}) and density of $\cD$ in $L^2(\mA,\rho_z;\bRd)$.
\end{remark}

Henceforth, assume Definition~\ref{def:countable_residual_multiplicity} holds on $\cC$, with $(\rho^n,f_n)_{n \geqslant 1}$ fixed in (ii) satisfying $f_n \in \cD$. A section $v : \oOmega \ni z \mapsto v_z \in K_z$ is \textit{predictable} if $\oOmega \ni z \mapsto \langle v_z, [\rho^n,f_n]_z \rangle_z \in \bR$ is predictable for every $n$. This makes $\cK$ a measurable field of Hilbert spaces over $(\oOmega,\sP^{\bH^{\cC}}, \oP)$, with $([\rho^n,f_n])_{n \geqslant 1}$ a fundamental sequence, in the usual sense of~\cite[Definition~IV.8.9]{Takesaki2002Theory}. In particular, norms and inner products of predictable sections are predictable. By~(\ref{eq:form}), the section $[\rho,f]$ is predictable for all $\rho\in \cC, f \in \BL(\mA;\bRd)$.

We construct a common driver of the assignment using a linear isometry $\cI$ on predictable sections of $\cK$ such that $\cI([\rho,f]) = N^{\rho,f}$. Write
$$
    L^2_{\pred}(\cK) \coloneq \int_{\oOmega}^\oplus K_z \, \oP(\di z)
$$
for the direct integral of $\cK$ over $(\oOmega,\sP^{\bH^\cC}, \oP)$, consisting of equivalence classes, up to equality $\oP$-a.e., of predictable sections with $\|v\|^2_2 \coloneq E[\int_0^T \|v_t\|^2 \dint t] < \infty$. For $\rho \in \cC, f \in \BL(\mA;\bRd)$, $[\rho,f] \in L^2_{\pred}(\cK)$ as $\|[\rho,f]_z\|^2 = \kappa^{\rho,\rho}_z(f,f)\leqslant \|f\|_\infty^2$. Let $\cE(\cK) \subset L^2_{\pred}(\cK)$ consist of the sections of the form  $v = \sum_{(\rho,f) \in F} \xi^{\rho,f}[\rho,f]$, where $F \subset \cC \times \BL(\mA;\bRd)$ is finite and each $\xi^{\rho,f} : \oOmega \rightarrow \bR$ is bounded and predictable. Define 
\begin{equation}
    \label{eq:integral_elementary}
    \cI_t(v) \coloneq \sum_{(\rho,f) \in F} \int_0^t \xi^{\rho,f}_s \dint N^{\rho,f}_s.
\end{equation}
As $N^{\rho,f} \in \cM^{2,c}_0(\bF)$ and $\xi^{\rho,f}$ is bounded, $\cI(v) \in \cM^{2,c}_0(\bF)$. By construction of $\kappa$, we have $\di [\cI(v), \cI(v')]_t = \langle v_t, v_t' \rangle \dint t$ for $v,v' \in \cE(\cK)$. In particular, $\cI(v)$ is independent of the choice of representative of~$v$, hence it is well defined. For $v \in \cE(\cK)^\perp$, $\int_{\oOmega} \xi_z \langle v_z, [\rho^n,f_n]_z\rangle \, \oP(\di z) = 0$ for all $n$ and bounded predictable~$\xi$, and taking $\xi = \sgn \langle v,[\rho^n,f_n]\rangle$ shows $\langle v_z, [\rho^n,f_n]_z\rangle = 0$ for all $n$, so that $v_z = 0$ for $\oP$-a.e.~$z$ by~(\ref{eq:totality_of_fundamental_sequence}). Thus $\cE(\cK) \subset L^2_{\pred}(\cK)$ is dense, and so $\cI$ extends uniquely to a linear isometry
\begin{equation}
    \label{eq:I}
    \cI : L^2_{\pred}(\cK) \longrightarrow \cM_0^{2,c}(\bH^\cC) \cap \cM_0^{2,c}(\bF), \qquad \di [\cI(v), \cI(v')]_t = \langle v_t, v_t' \rangle \dint t,
\end{equation}
where the codomain carries the terminal $L^2$-norm; the quadratic covariation identity in~(\ref{eq:I}) follows from Doob's inequality and the Kunita--Watanabe inequality, which moreover give $[\cI(v),Z] \equiv 0$ for every $Z \in \cM_{\loc,0}(\bG)$ since $[N^{\rho,f}, Z] \equiv 0$. By~(\ref{eq:integral_elementary}), $\cI([\rho,f]) = N^{\rho,f}$.

A \textit{strongly predictable isometric embedding} $J : \cK \rightarrow \mH_0$ is a family of linear isometries $J_z: K_z \rightarrow \mH_0, z \in \oOmega$, such that $\oOmega \ni z \mapsto J_z v_z \in \mH_0$ is predictable for each predictable section $v$ of $\cK$. Such $J$ exists, e.g.\ $J_z v_z \coloneq \sum_{n=1}^\infty \langle v_z, v^n_z \rangle \me_n$, with $(\me_n)_{n \geqslant 1}$ the orthonormal basis of $\mH_0$ fixed in Section~\ref{sec:4}, and $(v^n)_{n \geqslant 1}$ obtained from fiberwise zero-preserving Gram--Schmidt applied to $([\rho^n,f_n])_{n \geqslant 1}$. Their nonzero values form an orthonormal basis of $K_z$ for $\oP$-a.e.~$z$, and the $v^n$ are predictable as the Gram--Schmidt procedure uses finitely many predictable inner products and normalizations at each step; see also~\cite[Proposition~14]{Yetter2005MeasurableCategories}.  We thus fix an arbitrary strongly predictable isometric embedding $J : \cK \rightarrow \mH_0$, and correspondingly define
\begin{equation}
    \label{eq:Q}
    Q  : \oOmega \longrightarrow \cL(\mH), \quad Q_z \coloneq I_d \oplus J_z J_z^*.
\end{equation}
Then $Q$ is predictable, and $Q_z$ is the orthogonal projection of $\mH$ onto $\bRd \oplus \ran J_z$, keeping all ground directions and only the directions in $\mH_0$ used by $\cK$.

For $S \in \cS, t \in [0,T], x \in \bR^{m_\mS}, \mu \in \cP_2(\mA)$, define the operator
$$
    \Sigma_\mS(t,x,\mu) : L^2(\mA,\mu;\bRd) \longrightarrow \bR^{m_\mS}, \quad \Sigma_\mS(t,x,\mu) h \coloneq \int_\mA \sigma_\mS(t,x,a) h(a) \, \mu(\di a).
$$
This satisfies $\Sigma_\mS(t,x,\mu)^*y = \sigma_\mS(t,x,\cdot)^\top y$ and $\|\Sigma_\mS(t,x,\mu)\|_{\HS}^2 = \int_\mA |\sigma_\mS(t,x,a)|^2 \, \mu(\di a)$. The measure argument distinguishes it from the analogous profile operator of Section~\ref{sec:4}.

To express the controlled diffusions in $\mH$, consider for $\rho \in \cC$ its action-integrand field 
\begin{equation}
    \label{eq:H^rho}
    \cH^\rho = (H^\rho_z)_{z \in \oOmega}, \quad H^\rho_z \coloneq L^2(\mA,\rho_z;\bRd).
\end{equation}
Let $[f]_{\rho_z}$ be the equivalence class of $f \in \BL(\mA;\bRd)$ in $H^\rho_z$. Then $\ospn\{[f]_{\rho_z} : f \in \cD\} = H^\rho_z$ for each $z \in \oOmega$, and $z \mapsto \langle [f]_{\rho_z}, [f']_{\rho_z}\rangle = \rho_z(f^\top f')$ is predictable; thus $\cH^\rho$ is a predictable Hilbert field with fundamental sequence $([f]_{\rho})_{f \in \cD}$. Define
\begin{equation}
    \label{eq:T^rho}
    T^\rho_z[f]_{\rho_z} \coloneq \rho_z(f) + J_z[\rho,f]_z \in \mH, \quad f \in \cD.
\end{equation}
As $J_z$ is an isometry, (\ref{eq:form}) shows $\|T^\rho_z[f]_{\rho_z}\|^2 = |\rho_z(f)|^2 + \kappa^{\rho,\rho}_z(f,f) = \rho_z(|f|^2) = \|[f]_{\rho_z}\|^2$, thus~(\ref{eq:T^rho}) extends uniquely to a strongly predictable isometric embedding 
$$
    T^\rho : \cH^\rho \longrightarrow \mH.
$$
Note that $T^\rho_z c = c$ for $c \in \bRd$, $Q_zT^\rho_z = T^\rho_z$, and for any $\rho,\rho' \in \cC$, $\oP$-a.e.,
\begin{equation}
    \label{eq:T^rho_covariation}
    \langle T^\rho_z h, T^{\rho'}_zh'\rangle = \int_{\mA^2} h(a)^\top h'(a') \, \gamma^{\rho,\rho'}_z(\di a, \di a'), \quad h \in H^\rho_z, \, h' \in H^{\rho'}_z,
\end{equation}
by density of $\cD$. The composite $\Sigma_\mS(t,x,\rho_z) (T^\rho_z)^*$ acts on $\mH$ for every $\rho \in \cC$.

\begin{theorem}
    \label{theorem:ordinary_relaxed_assignment_common_noise}
    Suppose $(\rho,\mS) \mapsto Y^{\rho;\mS}$ satisfies Definition~\ref{def:pair_approximation} and Definition~\ref{def:countable_residual_multiplicity} on~$\cC$. Then $\bG \hookrightarrow \bF$, and for every strongly predictable isometric embedding $J : \cK \rightarrow \mH_0$, there is a cylindrical continuous square-integrable $\bF$-martingale $\mM$ in $\mH$, adapted to $\bH^\cC$, with covariance density $Q$ in~(\ref{eq:Q}), lying above $W$, such that $[\mM(h),Z]\equiv 0$ for all $h \in \mH_0$, $Z \in \cM_{\loc,0}(\bG)$, and for every $(\rho,\mS) \in \cC \times \cS$, $Y^{\rho;\mS}$ is the unique $\bF$-strong solution of
    \begin{equation}
        \label{eq:Y^rho,S_SDE}
        Y^{\rho;\mS}_t = x_0^\mS + \int_0^t \ob_\mS(s,Y^{\rho;\mS}_s, \rho_s) \dint s + \int_0^t \Sigma_\mS(s, Y^{\rho;\mS}_s, \rho_s)(T^\rho_s)^* \dint \mM_s.
    \end{equation}
\end{theorem}

\begin{proof}
    By Lemma~\ref{lemma:ambient_filtration_characteristics_ordinary}, $\bG \hookrightarrow \bF$. Set $\mM(c+h) \coloneq c^\top W + \cI(J^* h) \in \cM^{2,c}_0(\bH^\cC) \cap \cM_0^{2,c}(\bF)$ for $c + h \in \bRd \oplus \mH_0 = \mH$. Predictability of $J^* h$, $\|J^*_z h\|\leqslant \|h\|$, (\ref{eq:I}), and Doob's inequality show that $\mM$ is well defined and continuous into $L^2(\Omega;C[0,T])$, and is therefore a cylindrical continuous square-integrable $\bF$-martingale in $\mH$, adapted to $\bH^\cC$. Moreover, $\mM$ lies above $W$, and the orthogonality assertion follows from that of $\cI$. Finally, $\mM$ has covariance density $Q$ because $\di[\mM(c+h), \mM(c'+h')]_t = (c^\top c' + \langle J^*_t h, J^*_t h' \rangle) \dint t = \langle Q_t(c+h), c' + h'\rangle \dint t$.
    
    Abbreviate $\Phi_t \coloneq \Sigma_\mS(t,Y^{\rho;\mS}_t,\rho_t)(T^\rho_t)^*$. Taking adjoints shows $\Phi_t Q_t = \Phi_t$, as $Q_t T^\rho_t = T^\rho_t$. For $\rho' \in \cC, f \in \cD$, the construction of $N^{\rho',f}$, Proposition~\ref{prop:gamma^rho,rho'}, and Lemma~\ref{lemma:ambient_filtration_characteristics_ordinary} give
    $$
        \di [L^{\rho;\mS}, N^{\rho',f}]_t = \bigg( \int_{\mA^2} \sigma_\mS(t,Y^{\rho;\mS}_t,a) f(a') \, \gamma^{\rho,\rho'}_t(\di a, \di a') - \osigma_\mS(t,Y^{\rho;\mS}_t, \rho_t) \rho_t'(f) \bigg) \dint t.
    $$
    By (\ref{eq:T^rho})--(\ref{eq:T^rho_covariation}), the quantity in parentheses is
    $\Phi_t J_t[\rho',f]_t$; taking bounded predictable linear combinations of the fundamental sections and passing to limits by Kunita--Watanabe gives $\di [L^{\rho;\mS}, \cI(v)]_t = \Phi_t J_t v_t \dint t$ for $v \in L^2_{\pred}(\cK)$. Taking $v = J^* h$ for $h \in \mH_0$  gives $\di [L^{\rho;\mS}, \mM(h)]_t = \Phi_t Q_t h \dint t = \Phi_t h \dint t$. Also, $\di [L^{\rho;\mS},\mM(c)]_t = \osigma_\mS(t,Y^{\rho;\mS}_t,\rho_t) c \dint t= \Phi_t c \dint t$ for $c \in \bRd$, by Lemma~\ref{lemma:ambient_filtration_characteristics_ordinary} and $T^\rho c = c$. Hence,
    \begin{equation}
        \label{eq:ordinary_relaxed_assignment_common_noise_1}
        \di [L^{\rho;\mS}, \mM(h)]_t = \Phi_t h \dint t, \quad h \in \mH.
    \end{equation}
    Moreover, $\di [L^{\rho;\mS}]_t = \Sigma_\mS(t,Y^{\rho;\mS}_t,\rho_t)\Sigma_\mS(t,Y^{\rho;\mS}_t,\rho_t)^* \dint t = \Phi_t Q_t \Phi_t^* \dint t$ by Lemma~\ref{lemma:ambient_filtration_characteristics_ordinary} and as $T^\rho_t$ is an isometry with $Q_t T^\rho_t = T^\rho_t$. So by Lemma~\ref{lemma:cylindrical_integral_identification}, $L^{\rho;\mS} = \Phi \bigcdot \mM$, which is precisely~(\ref{eq:Y^rho,S_SDE}). Composition with $(T^\rho_t)^*$ preserves Hilbert--Schmidt norms, hence strong well-posedness of~(\ref{eq:Y^rho,S_SDE}) follows as in Proposition~\ref{prop:cylindrical_well_posed}, with $Y^{\rho;\mS}$ its unique $\bF$-strong solution.
\end{proof}

Consider the \textit{residual multiplicity} of the assignment on $\cC$,
$$
    k_{\cC} \coloneq \esssup_{z \in \oOmega} \dim K_z \in \bN_0 \cup\{\infty\}.
$$
Even when $\cC$ is a singleton, $k_\cC$ need not be finite.

\begin{example}
    \label{example:infinite_residual_multiplicity}
    Take $\mA = [0,1]$, $d = 1$, and $\cC = \{\rho\}$ with $\rho_t(\di a) = \di a$. The bounded Lipschitz functions $f_n(a) \coloneq \sqrt{2}\cos(2\pi n a), n \geqslant 1$, are centered and orthonormal in $L^2[0,1]$. Thus, the identity $\gamma^{\rho,\rho} = \Delta_\# \rho$ gives $\langle [\rho,f_n]_z, [\rho,f_m]_z \rangle = \delta_{nm}$ $\oP$-a.e., so that $k_\cC = \infty$.
\end{example}

When $k_\cC = \infty$, no finite-dimensional continuous local martingale is a common residual driver for the assignment on $\cC \times \cS$. 

\begin{lemma}
    Suppose, for some $k \in \bN$, there are $M \in \cM_{\loc,0}^c(\bF;\bR^k)$ and $\bF$-predictable $D^n : \oOmega \rightarrow \bR^{1 \times k}$ such that $N^{\rho^n,f_n} = D^n \bigcdot M$ for all $n \geqslant 1$. Then $k_\cC \leqslant k$.
\end{lemma}

\begin{proof}
    For each $\ell \in \bN$, let $C^\ell \coloneq (\kappa^{\rho^i,\rho^j}(f_i,f_j))_{i,j=1}^\ell$ be the Gram matrix of the first $\ell$ elements of the fixed fundamental sequence. With $N^\ell \coloneq \col_{n=1}^\ell N^{\rho^n,f_n}$ and $V^\ell \coloneq \col_{n=1}^\ell D^n$, it holds $N^\ell = V^\ell \bigcdot M$ and $\di [N^\ell]_t = C^\ell_t \dint t$. Let $P^\ell_t$ be the orthogonal projection onto $\ran V^\ell$. Then $P^\ell \bigcdot N^\ell = (P^\ell V^\ell) \bigcdot M = N^\ell$, thus $C^\ell_t \dint t = \di [P^\ell \bigcdot N^\ell, N^\ell]_t = P^\ell_t C^\ell_t \dint t$. So $C^\ell = P^\ell C^\ell$ $\oP$-a.e., and therefore $\ran C^\ell_t \subset \ran P^\ell_t = \ran V^\ell_t$. Since $V^\ell_t$ is $\ell \times k$-dimensional, $\rank C^\ell_t \leqslant k$. Hence~(\ref{eq:totality_of_fundamental_sequence}) gives $\dim K_z = \sup_{\ell \geqslant 1} \rank C^\ell_z \leqslant k$ $\oP$-a.e., which completes the proof.
\end{proof}

\subsection{Finite-family coherence and common labels}

When the approximation property in Definition~\ref{def:pair_approximation} is upgraded to all finite families in $\cC$, the $\rho \in \cC$ come from coherent parametrizations inducing the action-pair couplings of Proposition~\ref{prop:gamma^rho,rho'}, and the underlying dynamics  can be written in the form~(\ref{eq:SDE_cylindrical}). Retain Definition~\ref{def:pair_approximation}, and Definition~\ref{def:countable_residual_multiplicity} on~$\cC$.

\begin{definition}
    \label{def:finite_family_approximation}
    The assignment has the \textit{finite-family marginal approximation property} on $\cC$ if, for any finite family $(\rho^i,\mS^i)_{i \in I} \subset \cC \times \cS$, there are $(u^{n,i})_{i \in I} \subset \cU(\bK^n)$ such that
    $$
        \big( (\eta^{u^{n,i}})_{i\in I}, (X^{u^{n,i}; \mS^i})_{i\in I}\big) \Longrightarrow_{\sG_T} \big( (\eta^{\rho^i})_{i\in I}, (Y^{\rho^i; \mS^i})_{i\in I}\big)
    $$
    in $\cY^T_2(\mA)^I \times \prod_{i\in I} C([0,T];\bR^{m_{\mS^i}})$.
\end{definition}

\begin{definition}
    A family $\Pi = (\pi^\rho)_{\rho \in \cC} \subset \cR(\bH^\cC)$ is a \textit{common-label lift} of the assignment on~$\cC$ if, for every $\rho,\rho' \in \cC$,
    \begin{equation}
        \label{eq:common_label_lift_pushforward_identity}
        (\pi^\rho_z,\pi^{\rho'}_z)_\# \dint r = \gamma^{\rho,\rho'}_z \quad \oP\text{-a.e.},
    \end{equation}
    where the $\gamma^{\rho,\rho'}$ are as in Proposition~\ref{prop:gamma^rho,rho'}.
    Since $\gamma^{\rho,\rho}_z = \Delta_\# \rho_z$, $(\pi^\rho_z)_\# \dint r = \rho_z$ $\oP$-a.e.
\end{definition}

Definition~\ref{def:finite_family_approximation} yields compatible action kernels on countable families, with the remaining relaxed controls realized as functions of the countable generating family by Definition~\ref{def:countable_residual_multiplicity}.

\begin{proposition}
    \label{prop:common_label_lift}
     Suppose $(\rho,\mS) \mapsto Y^{\rho;\mS}$ satisfies Definition~\ref{def:pair_approximation}, and Definitions~\ref{def:countable_residual_multiplicity} and~\ref{def:finite_family_approximation} on $\cC$. Then the assignment admits a common-label lift on $\cC$.
\end{proposition}

\begin{proof}
    \textit{Step 1. Compatible kernels on countable families.} Let $I \subset J\subset \cC$ be nonempty countable. Equip the space of predictable probability kernels on $\mA^J$ (analogously, on $\mA^I$), modulo equality $\oP$-a.e., with the Hausdorff topology generated by
    \begin{equation}
        \label{eq:coupling_family_transition_topology_generator}
        \lambda \longmapsto \int_{\oOmega}\int_{\mA^J} \zeta_z(a) \, \lambda_z(\di a) \, \oP(\di z)
    \end{equation}
    for all $\sP^{\bH^\cC} \otimes \sB(\mA^J)$-measurable $\zeta : \oOmega \times \mA^J \rightarrow \bR$, continuous in $a \in \mA^J$ and dominated in absolute value by a function in $L^1(\oP)$. Define $\pr^J_I : \mA^J \rightarrow \mA^I, (a_\rho)_{\rho \in J} \mapsto (a_\rho)_{\rho \in I}$; the induced map $\lambda \mapsto (\pr^J_I)_\# \lambda$ is continuous for the above topology. Consider the subspace $\fP_J$ of those $\lambda_z \in \cP(\mA^J)$ with marginals $(\pr^J_\rho)_\# \lambda = \rho, \rho \in J$, and the further subspace $\fC_J$ of those $\lambda$ additionally satisfying $(\pr^J_\rho,\pr^J_{\rho'})_\# \lambda = \gamma^{\rho,\rho'}, \rho,\rho' \in J$. Now $\fP_J$ is nonempty since it contains the product kernel, and moreover it is compact. Indeed, for each $\rho \in J$, the finite Borel measure $\nu_\rho(A) \coloneq \int_{\oOmega} \rho_z(A) \, \oP(\di z)$ is tight since $\mA$ is Polish. For fixed $\varepsilon > 0$, choose positive $(\varepsilon_\rho)_{\rho \in J}$ with $\sum_{\rho \in J} \varepsilon_\rho \leqslant \varepsilon$ and nonempty compacts $F_\rho \subset \mA$ satisfying $\nu_\rho(F^c_\rho) < \varepsilon_\rho$. Then $F^{\varepsilon} \coloneq \prod_{\rho \in J} F_\rho$ is compact and $\sup_{\lambda \in \fP_J} \int_{\oOmega} \lambda_z(\mA^J \backslash F^\varepsilon) \, \oP(\di z) \leqslant \varepsilon$. By~\cite[Remark~2.4 and Theorem~2.2]{Balder1989Prohorov}, $\fP_J$ is relatively compact; continuity of $(\pr^J_\rho)_\#$ yields compactness of $\fP_J$. 
    
    Each constraint $(\pr^J_{\rho,\rho'})_\# \lambda = \gamma^{\rho,\rho'}$ is closed. Choose finitely many of them, and let $F \subset J$ be the set of relaxed controls they use. Enumerate $\cD = (f_j)_{j \geqslant 1}$ and apply Definition~\ref{def:finite_family_approximation} to $(\rho, \oplus_{j=1}^k \mP^{f_j})_{\rho \in F} \subset \cC \times \mS$ for each $k \geqslant 1$; then diagonalize as in Proposition~\ref{prop:gamma^rho,rho'}. Lemma~\ref{lemma:common_kernel} supplies a predictable coupling $\lambda^F$ on $\mA^F$ whose pairwise marginals agree with the prescribed $\gamma^{\rho,\rho'}$ $\oP$-a.e., by the determining property of $\cA_0 \otimes \cA_0$. Then $\lambda^F_z \otimes (\otimes_{\rho \in J \backslash F} \rho_z)$ is an element of $\fP_J$ and satisfies the chosen constraints. Thus, the constraints $(\pr^J_{\rho,\rho'})_\# \lambda = \gamma^{\rho,\rho'}$ have the finite intersection property, and therefore their intersection $\fC_J$ is nonempty and compact.  

    \textit{Step 2. Common marginal of generating controls.} Let $([\rho^n,f_n])_{n \geqslant 1}$ be the fixed fundamental sequence of $\cK$ and set $\cC_0 \coloneq \{\rho^n : n \geqslant 1\}$ and $E \coloneq \mA^{\cC_0}$. For every $\rho \in \cC$, continuity shows $(\pr^{\cC_0 \cup \{\rho\}}_{\cC_0})_\# \fC_{\cC_0 \cup\{\rho\}}$ is a nonempty compact, hence closed, subset of the compact Hausdorff space $\fC_{\cC_0}$. They have the finite intersection property, since for any finite $F \subset \cC$, $\fC_{\cC_0 \cup F}$ is nonempty and projects into $\fC_{\cC_0 \cup \{\rho\}}$ for all $\rho \in F$. So take
    $$
        \mu \in \bigcap_{\rho \in \cC} (\pr^{\cC_0 \cup \{\rho\}}_{\cC_0})_\# \fC_{\cC_0\cup\{\rho\}}.
    $$
    \textit{Step 3. Representing remaining relaxed controls.} Put $v^n \coloneq [\rho^n,f_n]$ and let $\xi^n : \oOmega \times E \rightarrow \bR^d$ be given by $\xi^n_z(e) \coloneq f_n(e_{\rho^n}) - \rho^n_z(f_n)$. As $\mu \in \fC_{\cC_0}$, $(v^n)_{n \geqslant 1}$ and $(\xi^n_z)_{n \geqslant 1}$ have the same Gram matrix for $\oP$-a.e.~$z$. Thus by (\ref{eq:totality_of_fundamental_sequence}), $v^n_z \mapsto \xi^n_z$ extends uniquely to a linear isometry 
    $$
        U_z : K_z \longrightarrow L^2(E,\mu_z;\bRd).
    $$
    Fix $\rho \in \cC$. If $\rho \notin \cC_0$, Step~2 gives $\lambda^\rho \in \fC_{\cC_0 \cup\{\rho\}}$ with $(\pr^{\cC_0 \cup \{\rho\}}_{\cC_0})_\# \lambda^\rho = \mu$. If $\rho \in \cC_0$, take $\lambda^\rho \coloneq (\id_E, \pr^{\cC_0}_\rho)_\# \mu$. In either case, $\lambda^\rho_z$ has $a$-marginal $\rho_z$ and $(a,e_{\rho^n})$-marginal $\gamma^{\rho,\rho^n}_z$. So $\|f(a) -\rho_z(f) - (U_z v)(e)\|_{L^2(E \times \mA, \lambda_z^\rho; \bRd)} = \|[\rho,f]_z - v\|$ for $f \in \cD, v \in \spn\{v^n_z : n \geqslant 1\}$. As $\cD$ and the set of finite rational linear combinations of $([\rho^n,f_n])_{n \geqslant 1}$ are countable, the preceding identity holds, for each fixed $\rho$, outside a common $\oP$-null set. Approximating $[\rho,f]_z$ by such $v$ by~(\ref{eq:totality_of_fundamental_sequence}), and noting that $\lambda^\rho_z$ has $E$-marginal $\mu_z$, yield 
    \begin{equation}
        \label{eq:common_label_lift_U_identity}
        f(a)-\rho_z(f) = (U_z[\rho,f]_z)(e) \quad \text{in } L^2(E \times \mA, \lambda_z^\rho;\bRd), \quad f \in \cD.
    \end{equation}
    Disintegrate $\lambda^\rho_z(\di e, \di a) = \mu_z(\di e) \, q^\rho_z(e,\di a)$ $\oP$-a.e., where $q^\rho$ is $\sP^{\bH^\cC} \otimes \sB(E)$-measurable \cite[Corollary~3.6]{kallenberg2021foundations}. Choosing $f = g\epsilon_1$ for $g \in \cA_0$ in~(\ref{eq:common_label_lift_U_identity}) shows that $g(a)$ is $\lambda^\rho_z$-a.e.~a function of $e$. Thus $g(a) = \int_{\mA} g(a') \, q^\rho_z(e,\di a')$ $\lambda^\rho_z$-a.e. Since $\mA_0$ is countable and determining, $q_z^\rho(e,\cdot) = \delta_a$ $\lambda_z^\rho$-a.e. Let $\chi : \cP(\mA) \rightarrow \mA$ be the Borel inverse of the Dirac embedding on its image, extended by $\ma_0$ elsewhere. Then $a^\rho_z(e) = \chi(q^\rho_z(e,\cdot))$ is jointly measurable, and also $\lambda^\rho_z = (\id_E, a^\rho_z)_\# \mu_z$ for $\oP$-a.e.~$z$. Hence (\ref{eq:common_label_lift_U_identity}) is equivalently $f \circ a^\rho_z - \rho_z(f) = U_z[\rho,f]_z$ in $L^2(E, \mu_z;\bRd)$ for all $f \in \cD$. Thus for fixed $\rho,\rho' \in \cC$, for $\oP$-a.e.~$z$,
    $$
        \langle f \circ a^\rho_z - \rho_z(f), f' \circ a^{\rho'}_z - \rho'_z(f')\rangle_{L^2(E,\mu_z;\bRd)} = \kappa^{\rho,\rho'}_z(f,f'), \quad f,f' \in \cD.
    $$
    Comparing with the construction of $\kappa$ and using the determining property of $\cA_0 \otimes \cA_0$ shows $(a^\rho_z, a^{\rho'}_z)_\# \mu_z = \gamma^{\rho,\rho'}_z$ for $\oP$-a.e.~$z$. There is a $\sP^{\bH^\cC}\otimes \sB$-measurable map $\Psi : \oOmega \times [0,1] \rightarrow E$ with $(\Psi_z)_\# \dint r = \mu_z$~\cite[Lemma~4.22]{kallenberg2021foundations}. Setting $\pi^\rho_z \coloneq a^\rho_z\circ \Psi_z$ completes the proof. 
\end{proof}

Assume Definition~\ref{def:finite_family_approximation} on $\cC$, and fix a common-label lift $\Pi = (\pi^\rho)_{\rho \in \cC}$ on $\cC$. With $([\rho^n,f_n])_{n \geqslant 1}$ as in Definition~\ref{def:countable_residual_multiplicity}(ii), the predictable $\mH_0$-valued processes $f_n \circ \pi^{\rho^n}_z - \rho^n_z(f_n)$, $n \geqslant 1$, have the same Gram matrix as $[\rho^n,f_n]_z, n \geqslant 1$, by~(\ref{eq:common_label_lift_pushforward_identity}). Hence the prescription
$$
    J^\Pi_z [\rho^n,f_n]_z \coloneq f_n \circ \pi^{\rho^n}_z - \rho^n_z(f_n)
$$
extends uniquely to a strongly predictable isometric embedding $J^\Pi : \cK \rightarrow \mH_0$. For each fixed $\rho \in \cC$, $J^\Pi_z[\rho,f]_z = f \circ \pi^\rho_z - \rho_z(f)$ for all $f \in\cD$ $\oP$-a.e. Define
\begin{equation}
    \label{eq:Q^Pi}
    Q^\Pi : \oOmega \longrightarrow \cL(\mH), \qquad Q^\Pi_z \coloneq I_{\bRd} \oplus J^\Pi_z (J^\Pi_z)^*.
\end{equation}
Taking $J = J^\Pi$ in the construction of $T^\rho$ in~(\ref{eq:T^rho}) gives $T^\rho_z[f]_{\rho_z} = f \circ \pi^\rho_z$ for all $f \in \cD$ $\oP$-a.e. Density and continuity of $T^\rho_z$ yield $T^\rho_z h = h\circ \pi^\rho_z$ for all $h \in H^\rho_z$ $\oP$-a.e. 

\begin{theorem}
    \label{theorem:ordinary_relaxed_assignment_common_label_dynamics}
    Suppose $(\rho,\mS) \mapsto Y^{\rho;\mS}$ satisfies Definition~\ref{def:pair_approximation}, and Definitions~\ref{def:countable_residual_multiplicity} and~\ref{def:finite_family_approximation} on $\cC$. Then $\bG \hookrightarrow \bF$. For any common-label lift $\Pi = (\pi^\rho)_{\rho \in \cC}$ on $\cC$ (which exists), there is a cylindrical continuous square-integrable $\bF$-martingale $\mM^\Pi$ in $\mH$, adapted to $\bH^\cC$, with covariance density $Q^\Pi$ given by~(\ref{eq:Q^Pi}), lying above~$W$, such that $[\mM^\Pi(h), Z] \equiv 0$ for all $h\in \mH_0, Z \in \cM_{\loc,0}(\bG)$, and for every $(\rho,\mS) \in \cC \times \cS$, $Y^{\rho;\mS}$ is the unique $\bF$-strong solution of
    \begin{equation}
        \label{eq:Y^rho,S_SDE_parametrized}
        Y^{\rho;\mS}_t = x_0^\mS + \int_0^t \ob_\mS(s,Y^{\rho;\mS}_s, \rho_s) \dint s + \int_0^t \Sigma_\mS(s, Y^{\rho;\mS}_s, \pi^\rho_s) \dint \mM^\Pi_s.
    \end{equation}
\end{theorem}

\begin{proof}
    Immersion is Lemma~\ref{lemma:ambient_filtration_characteristics_ordinary}, and a common-label lift $\Pi$ exists by Proposition~\ref{prop:common_label_lift}. For $y \in \bR^{m_\mS}$, $(\Sigma_\mS(t,x,\rho_t)(T^\rho_t)^*)^* y = T^\rho_t(\sigma_\mS(t,x,\cdot)^\top y) = \sigma_\mS(t,x,\pi^\rho_t)^\top y = \Sigma_\mS(t,x,\pi^\rho_t)^* y$. Thus applying Theorem~\ref{theorem:ordinary_relaxed_assignment_common_noise} with $J = J^\Pi$ supplies $\mM^\Pi$ and identifies~(\ref{eq:Y^rho,S_SDE}) with (\ref{eq:Y^rho,S_SDE_parametrized}).
\end{proof}

The driver $\mM^\Pi$ is cylindrical Brownian motion precisely when $J^\Pi$ is surjective $\oP$-a.e. In Example~\ref{example:infinite_residual_multiplicity}, taking $\pi^\rho(r) \coloneq r$ gives $J^\Pi [\rho,f] = f - \rho(f)$ for $f \in \cD$, so $J^\Pi$ is surjective as $\cD$ is dense in $L^2[0,1]$. In any case, $\mM^\Pi$ can always be completed to a cylindrical Brownian motion on a further Wiener extension.

\begin{corollary}
    \label{cor:complete_M_to_B}
    On a further Wiener extension, $\mM^\Pi$ can be completed to a cylindrical Brownian motion $\mB^\Pi$ in $\mH$, lying above $W$, with $\mB^\Pi|_{\mH_0}$ independent of $\sG_T$, such that~(\ref{eq:Y^rho,S_SDE_parametrized}) remains valid with $\mB^\Pi$ in place of $\mM^\Pi$.
\end{corollary}

\begin{proof}
    Adjoin a cylindrical Brownian motion $\mB^0$ in $\mH$, independent of $\sF_T$, and define
    $$
        \mB^\Pi(h) \coloneq \mM^\Pi(h) + \int_0^\cdot \langle (I_\mH- Q^\Pi_s)h, \di \mB^0_s \rangle, \quad h \in \mH.
    $$
    Since the two summands are strongly orthogonal and $Q^\Pi_z$ is an orthogonal projection, $\mB^\Pi$ is well defined and $\di [\mB^\Pi(h), \mB^\Pi(h')]_t = \langle h,h' \rangle \dint t$, hence $\mB^\Pi$ is cylindrical Brownian motion. As $Q^\Pi_z c = c$ for $c \in \bRd$, $\mB^\Pi$ lies above $W$; by independence of $\mB^0$ from $\sF_T$, $[\mB^\Pi(h),Z] \equiv 0$ for all $h \in \mH_0, Z \in \cM_{\loc,0}(\bG)$, and the further Wiener extension preserves immersion of~$\bG$. Thus Lemma~\ref{lemma:conditional_independence_BM} applied to $(\mB^\Pi(\me_j))_{j=1}^n$ for each $n$ shows that $\mB^\Pi|_{\mH_0}$ is independent of $\sG_T$. The identities $\Phi_t Q_t = \Phi_t$ and $\Sigma_\mS(t,x,\pi^\rho_t) = \Sigma_\mS(t,x,\rho_t)(T^\rho_t)^*$  established in the proofs of Theorems~\ref{theorem:ordinary_relaxed_assignment_common_noise} and~\ref{theorem:ordinary_relaxed_assignment_common_label_dynamics} give $\Sigma_\mS(t,x,\pi^\rho_t) (I_\mH - Q^\Pi_t) = 0$, and therefore
    $$
        \int_0^t \Sigma_\mS(s,Y^{\rho;\mS}_s,\pi^\rho_s) \dint \mB^\Pi_s = \int_0^t \Sigma_\mS(s,Y^{\rho;\mS}_s,\pi^\rho_s) \dint \mM^\Pi_s.
    $$
    This completes the proof.
\end{proof}

Thus, for a state assignment indexed by ordinary relaxed controls, the relevant control class admits coherent parametrizations under the standing assumptions. Moreover, the assignment is expressed as the strong solution map of dynamics of the form~(\ref{eq:SDE_cylindrical}).

\section{Martingale measures and strong formulations}
\label{sec:6}
We explain the relaxed dynamics studied in Sections~\ref{sec:4}--\ref{sec:5} through action-space martingale measures.

Fix a filtered extension $(\Omega,\sF_T,\bF,P)$ of the ground basis with $\bG \hookrightarrow \bF$. For $\rho \in \cQ$, retain the field $\cH^\rho = (H^\rho_z)_{z \in \oOmega}$, $H^\rho_z = L^2(\mA,\rho_z;\bRd)$ of~(\ref{eq:H^rho}). A section $f$ of $\cH^\rho$ is \textit{predictable} if $z \mapsto \langle f_z, [g]_{\rho_z} \rangle_{H^\rho_z}$ is $\bF$-predictable for all $g \in \cD$. Let $L^2_{\pred,\loc}(\cH^\rho)$ be the predictable sections~$f$, up to equality $\oP$-a.e.\ (equivalently, the $\sP^\bF \otimes \sB(\mA)$-measurable $f : \oOmega \times \mA \rightarrow \bRd$, up to equality $\oP(\di z) \, \rho_z(\di a)$-a.e.), with $\int_0^T \rho_s(|f_s|^2) \dint s < \infty$ a.s. Let $L^2_{\pred}(\cH^\rho)$ be the subspace for which $E[\int_0^T \rho_s(|f_s|^2) \dint s] < \infty$. 

We identify an action-space martingale measure with its stochastic-integral map. 

\begin{definition}
    \label{def:martingale_measure}
    A continuous orthogonal $d$-component $\bF$-\textit{martingale measure} with intensity $\eta^\rho$ is a linear map $\fM : L^2_{\pred,\loc}(\cH^\rho) \rightarrow \cM^{c}_{\loc,0}(\bF)$ satisfying
    \begin{equation}
        \label{eq:martingale_measure}
        \di[\fM(f), \fM(f')]_t = \rho_t(f_t^\top f_t') \dint t, \quad f,f' \in L^2_{\pred,\loc}(\cH^\rho).
    \end{equation}
\end{definition}

For $f \in L^2_{\pred}(\cH^\rho)$, $E[|\fM_T(f)|^2] = E[\int_0^T \rho_s(|f_s|^2) \dint s] < \infty$. So $\fM$ restricted to $L^2_{\pred}(\cH^\rho)$ is an isometry into $\cM^{2,c}_0(\bF)$ for the terminal $L^2$-norm, and is continuous
into $L^2(\Omega;C[0,T])$ by Doob's inequality. For $f \in L^2_{\pred,\loc}(\cH^\rho)$ and bounded predictable $\xi : \oOmega \rightarrow \bR$,
\begin{equation}
    \label{eq:martingale_measure_predictable_multiplication}
    \fM(\xi f) = \xi \bigcdot \fM(f)
\end{equation}
since $[\fM(\xi f) - \xi \bigcdot \fM(f)] \equiv 0$ by~(\ref{eq:martingale_measure}). In particular, stopping commutes with $\fM$, and thus $\fM$ is determined by its restriction to $L^2_{\pred}(\cH^\rho)$. For $f \in L^2_{\pred,\loc}(\cH^\rho)$, write 
$$
    \fM_t(f) \eqqcolon \int_0^t \int_\mA f_s(a)^\top \, \fM(\di s, \di a),
$$
and apply this notation row-wise for matrix-valued integrands. This coincides with the usual martingale-measure integral of Walsh~\cite{Walsh1986}. Indeed, for $C \in \sB(\mA), 1 \leqslant i \leqslant d$, the set-indexed components $\fM^i_t (C) \coloneq \fM_t(1_C \epsilon_i)$ satisfy $[\fM^i(C), \fM^j(C')]_t = \delta_{ij} \int_0^t \rho_s(C \cap C') \dint s$, and $\sigma$-additivity in $L^2$ follows by the isometry following (\ref{eq:martingale_measure}). Conversely, stochastic integration against such set-indexed $\fM$ satisfies Definition~\ref{def:martingale_measure}; see~\cite[Proposition~I-6]{ElKarouiMeleard1990}.

For a prescribed pair $(\rho,\fM)$ and $\cS \in \mS$, where $\fM$ has intensity $\eta^\rho$, consider
\begin{equation}
    \label{eq:SDE_martingale_measure}
    Y_t = x_0^\mS + \int_0^t \ob_\mS(s,Y_s,\rho_s) \dint s + \int_0^t \int_\mA \sigma_\mS(s,Y_s,a) \,\fM(\di s, \di a).
\end{equation}
The linear growth and state-Lipschitz assumptions on $b_\mS,\sigma_\mS$, together with~(\ref{eq:martingale_measure}), give strong well-posedness of~(\ref{eq:SDE_martingale_measure}) by the argument in Proposition~\ref{prop:strict_well_posed}. The unique $\bF$-strong solution is determined by $(\rho,\fM)$, generally not by the intensity $\eta^\rho$ alone; thus we denote it $Y^{\rho,\fM;\mS}$. If $\rho = \delta_u$ for $u \in \cU$, $\fM^u_t(f) \coloneq \int_0^t f_s(u_s)^\top \dint W_s$ is canonical, and~(\ref{eq:SDE_martingale_measure}) reduces to~(\ref{eq:SDE_strict}). Also, $\di [\fM^u(f), Z]_t = f_t^\top(u_t) \dint [W,Z]_t$ for $Z \in \cM_{\loc,0}(\bG)$. This motivates the following condition.

\begin{definition}
    A martingale measure $\fM$ with intensity $\eta^\rho$ is \textit{ground-compatible} if  
    \begin{equation}
        \label{eq:martingale_measure_ground_compatible}
        \di [\fM(f), Z]_t = \rho_t(f_t)^\top \dint [W,Z]_t
    \end{equation}
    for every $f \in L^2_{\pred,\loc}(\cH^\rho)$ and $Z \in \cM_{\loc,0}(\bG)$.
\end{definition}

For a ground-compatible martingale measure $\fM$ and any $\bF$-predictable $\xi : \oOmega \rightarrow \bRd$ with $\int_0^T |\xi_s|^2 \dint s < \infty$ a.s., regarded as constant in the action variable, (\ref{eq:martingale_measure}) and (\ref{eq:martingale_measure_ground_compatible}) with $Z = W^i$ give $[\fM(\xi) - \xi^\top \bigcdot W] \equiv 0$. Hence $\fM(\xi) = \xi^\top \bigcdot W$. Taking $\xi = \rho(f)$ and linearity yield
\begin{equation}
    \label{eq:ground_compatible_decomposition}
    \fM(f) = \rho(f)^\top \bigcdot W + \fM\big( f - \rho(f) \big)
\end{equation}
for all $f \in L^2_{\pred,\loc}(\cH^\rho)$. By~(\ref{eq:martingale_measure_ground_compatible}), $[\fM(f-\rho(f)),Z] \equiv 0$ for $Z \in \cM_{\loc,0}(\bG)$. Applied row-wise to the diffusion, (\ref{eq:ground_compatible_decomposition}) gives the ground and residual martingales of Proposition~\ref{prop:martingale_characteristics}, with residual covariance density $R_\mS(t,Y_t,\rho_t)$ by~(\ref{eq:martingale_measure}). Also, if $\rho = \delta_u$, (\ref{eq:ground_compatible_decomposition}) forces $\fM = \fM^u$.

We now relate action-space martingale measures to cylindrical martingales.

\begin{lemma}
    \label{lemma:cylindrical_and_MM_integral_correspondence}
    Fix $\rho \in \cQ$ and a strongly predictable isometric embedding $T^\rho : \cH^\rho \rightarrow \mH$ with $T^\rho c = c$ for $c \in \bRd$. If $\mM$ is a cylindrical continuous local $\bF$-martingale in $\mH$ with covariance density $Q$ satisfying $QT^\rho = T^\rho$, then 
    \begin{equation}
        \label{eq:cylindrical_to_MM}
        \fM_t(f) \coloneq \int_0^t \langle T^\rho_s f_s, \di \mM_s\rangle, \quad f \in L^2_{\pred,\loc}(\cH^\rho),
    \end{equation}
    is a martingale measure with intensity $\eta^\rho$. If additionally $\mM$ lies above $W$ and $[\mM(h), Z] \equiv 0$ for all $h \in \mH_0, Z \in \cM_{\loc,0}(\bG)$, then $\fM$ is ground compatible. Conversely, for a martingale measure $\fM$ with intensity $\eta^\rho$,
    \begin{equation}
        \label{eq:MM_to_cylindrical}
        \mM_t(h) \coloneq \fM_t\big((T^\rho)^* h\big), \quad h \in \mH,
    \end{equation}
    is a cylindrical continuous square-integrable martingale with covariance density $Q \coloneq T^\rho(T^\rho)^*$, and~(\ref{eq:cylindrical_to_MM}) holds. If additionally $\fM$ is ground compatible, then $\mM$ lies above $W$ and $[\mM(h),Z] \equiv 0$ for all $h \in \mH_0, Z \in \cM_{\loc,0}(\bG)$.
\end{lemma}

\begin{proof}
    Define $\fM$ by~(\ref{eq:cylindrical_to_MM}). For $f \in L^2_{\pred,\loc}(\cH^\rho)$, $T^\rho f$ is predictable. Since $QT^\rho = T^\rho$ and $T^\rho_s$ is an isometry, $\int_0^T \|Q_s^{1/2} T^\rho_s f_s\|^2 \dint s = \int_0^T \|T^\rho_s f_s\|^2 \dint s = \int_0^T \rho_s(|f_s|^2) \dint s < \infty$ a.s., so $\fM : L^2_{\pred,\loc}(\cH^\rho) \rightarrow \cM^c_{\loc,0}(\bF)$ is a well-defined linear map. Also, $(T^\rho)^*_t Q_t T^\rho_t = I_{H^\rho_t}$ and~(\ref{eq:cylindrical_covariation}) give $\di [\fM(f),\fM(f')]_t= \langle T^\rho_t f_t, Q_tT^\rho_t f'_t \rangle \dint t = \langle f_t, f'_t\rangle \dint t = \rho_t(f_t^\top f'_t) \dint t$, which is~(\ref{eq:martingale_measure}). By linearity, $T^\rho_t f_t = \rho_t(f_t) + T^\rho_t(f_t - \rho_t(f_t))$, and $T^\rho_t(f_t-\rho_t(f_t)) \in \mH_0 = (\bRd)^\perp$ because $T_t^\rho$ is an isometry and $T^\rho_t c = c$ for $c \in \bRd$. So under the additional hypotheses on $\mM$, $[\fM(f-\rho(f)), Z] \equiv 0$ for every $f \in L^2_{\pred,\loc}(\bH^\rho)$ and $Z \in \cM_{\loc,0}(\bG)$ by predictable approximation, as in Proposition~\ref{prop:residual_cylindrical_correspondence}. Ground-compatibility of $\fM$ follows.

    Conversely, define $\mM$ by~(\ref{eq:MM_to_cylindrical}). For $h \in \mH$, $(T^\rho)^*h$ is predictable. As $T^\rho_s$ is an isometry, $(T^\rho_s)^*$ is a contraction, so $E[\int_0^T \|(T^\rho_s)^* h\|^2 \dint s] \leqslant T \|h\|^2 < \infty$. Hence $\mM : \mH \rightarrow \cM^{2,c}_0(\bF)$ is a well-defined linear map. By~(\ref{eq:martingale_measure}), $\di [\mM(h), \mM(h')]_t = \langle (T^\rho_t)^*h, (T^\rho_t)^* h' \rangle \dint t = \langle Q_t h, h' \rangle \dint t$. Since $Q_t$ is an orthogonal projection, Doob's inequality shows $\mM$ is continuous into $L^2(\Omega;C[0,T])$, thus is a cylindrical continuous square-integrable martingale with covariance density $Q$. Lemma~\ref{lemma:cylindrical_integral_identification} with $L = \fM(f)$ and $\Phi_t h = \langle T^\rho_t f_t, h \rangle$ recovers~(\ref{eq:cylindrical_to_MM}). If $\fM$ is ground-compatible, $\mM(c) = \fM((T^\rho)^*c) = \fM(c) = c^\top W$ by $(T^\rho)^* c = (T^\rho)^* T^\rho c$ and~(\ref{eq:ground_compatible_decomposition}). Also, for $h \in \mH_0 = (\bRd)^\perp$, $\rho_t((T^\rho_t)^* h) = 0$, so $\di [\mM(h), Z]_t = \rho_t((T^\rho_t)^* h)^\top \dint [W,Z]_t = 0$ for every $Z \in \cM_{\loc,0}(\bG)$, as required.
\end{proof}

We now identify the martingale measures associated with a relaxed state assignment.

\begin{proposition}
    \label{prop:relaxed_state_assignment_MM}
    Let $(\rho,\mS) \mapsto Y^{\rho;\mS}$ be a relaxed state assignment on $\cQ \times \cS$ satisfying Definition~\ref{def:pair_approximation}. For each $\rho \in \cQ$, there is a unique $\bF$-martingale measure $\fM^\rho$ with intensity $\eta^\rho$ satisfying $\fM^\rho(f) = Y^{\rho;\mP^f}$ for all $f \in \cD$. It is ground-compatible, and $Y^{\rho;\mS}$ is the unique solution of~(\ref{eq:SDE_martingale_measure}) for $\fM^\rho$, for all $\mS \in \cS$. For $\rho,\rho' \in \cQ, f \in L^2_{\pred,\loc}(\cH^\rho), f' \in L^2_{\pred,\loc}(\cH^{\rho'})$,
     \begin{equation}
        \label{eq:recovered_martingale_measure_covariations}
        \di [\fM^\rho(f), \fM^{\rho'}(f')]_t = \bigg(\int_{\mA^2} f_t(a)^\top f'_t(a') \, \gamma^{\rho,\rho'}_t(\di a, \di a') \bigg) \dint t,
    \end{equation}
    where $\gamma^{\rho,\rho'}$ is from Proposition~\ref{prop:gamma^rho,rho'}.
    If Definition~\ref{def:countable_residual_multiplicity} is also satisfied on a nonempty $\cC \subset \cQ$, the common driver $\mM$ and isometry $T^\rho$ of Theorem~\ref{theorem:ordinary_relaxed_assignment_common_noise} realize $\fM^\rho$ through (\ref{eq:cylindrical_to_MM}) for all $\rho \in \cC$. 
\end{proposition}

\begin{proof}
    By~(\ref{eq:martingale_measure_predictable_multiplication}), any two martingale measures that agree on $\cD$ agree on every finite sum $\sum_{j=1}^n \xi_j f_j$ where $\xi_j : \oOmega \rightarrow \bR$ are bounded predictable and $f_j \in \cD$. These sums are dense in $L^2_{\pred}(\cH^\rho)$, so the isometry following~(\ref{eq:martingale_measure}) shows uniqueness on $L^2_{\pred}(\cH^\rho)$, hence on $L^2_{\pred,\loc}(\cH^\rho)$. For nonempty $\cC \subset \cQ$ on which Definition~\ref{def:countable_residual_multiplicity} holds, Theorem~\ref{theorem:ordinary_relaxed_assignment_common_noise} provides $\mM,Q,T^\rho$ such that $Q T^\rho = T^\rho$ and~(\ref{eq:Y^rho,S_SDE}) is satisfied. Each $T^\rho$ is strongly $\bF$-predictable. Lemma~\ref{lemma:cylindrical_and_MM_integral_correspondence} gives a ground-compatible martingale measure $\fM^{\rho,\cC}$ for each $\rho \in \cC$ such that~(\ref{eq:cylindrical_to_MM}) holds. As $(\Sigma_\mS(t,x,\rho_t)(T^\rho_t)^*)^* y = T^\rho_t[\sigma_\mS(t,x,\cdot)^\top y]_{\rho_t}$ for $y \in \bR^{m_\mS}$, (\ref{eq:cylindrical_to_MM}) identifies~(\ref{eq:SDE_martingale_measure}) with~(\ref{eq:Y^rho,S_SDE}) for $\fM^{\rho,\cC}$, with $Y^{\rho;\mS}$ their unique $\bF$-strong solution by Theorem~\ref{theorem:ordinary_relaxed_assignment_common_noise}. In particular, $\fM^{\rho,\cC}(f) = Y^{\rho;\mP^f}$ for all $f \in \cD$. By Remark~\ref{remark:countable_residual_multiplicity_automatic}, Definition~\ref{def:countable_residual_multiplicity} holds on countable $\cC \subset \cQ$. So taking $\cC = \{\rho\}$ for each $\rho \in \cQ$ gives $\fM^\rho \coloneq \fM^{\rho,\{\rho\}}$ satisfying all the preceding properties. It only remains to establish~(\ref{eq:recovered_martingale_measure_covariations}). Taking $\cC = \{\rho,\rho'\}$ yields the same martingale measures $\fM^\rho,\fM^{\rho'}$ by uniqueness, and the common driver of Theorem~\ref{theorem:ordinary_relaxed_assignment_common_noise}, (\ref{eq:cylindrical_to_MM}), (\ref{eq:cylindrical_covariation}), and $Q T^{\rho'} = T^{\rho'}$ give $\di [\fM^\rho(f), \fM^{\rho'}(f')]_t = \langle T^\rho_t f_t, Q_t T^{\rho'}_t f'_t \rangle \dint t = \langle T^\rho_t f_t, T^{\rho'}_t f_t' \rangle \dint t,$
    which shows~(\ref{eq:recovered_martingale_measure_covariations}) by~(\ref{eq:T^rho_covariation}). This completes the proof.
\end{proof}

\begin{remark}
    In set-indexed notation, (\ref{eq:recovered_martingale_measure_covariations}) reads, for $C,C' \in \sB(\mA), 1 \leqslant i,j \leqslant d$, 
    $$
        \di [\fM^{\rho,i}(C), \fM^{\rho',j}(C')]_t = \delta_{ij} \gamma^{\rho,\rho'}_t(C \times C') \dint t.
    $$
\end{remark}

\begin{remark}
    \label{remark:parametrized_to_MM}
    Consider the parametrized cylindrical dynamics~(\ref{eq:SDE_cylindrical}). For each $\pi \in \cR$, define the strongly predictable isometric embedding $T^\pi : \cH^{\rho^\pi} \rightarrow \mH$ by $T^\pi_z f \coloneq f \circ \pi_z$. Then $T^\pi c = c$ for $c \in \bRd$. Lemma~\ref{lemma:cylindrical_and_MM_integral_correspondence} therefore gives the ground-compatible martingale measure
    \begin{equation}
        \label{eq:canonical_cylindrical_MM}
        \fM^\pi_t(f) \coloneq \int_0^t \langle f_s \circ \pi_s, \di \mB_s \rangle, \quad f \in L^2_{\pred,\loc}(\cH^{\rho^\pi}),
    \end{equation}
    with intensity $\eta^{\rho^\pi}$. Applying (\ref{eq:canonical_cylindrical_MM}) row-wise to the diffusion identifies~(\ref{eq:SDE_cylindrical}) with~(\ref{eq:SDE_martingale_measure}), and
    $$
        X^{\pi;\mS} = Y^{\rho^\pi,\fM^\pi;\mS}, \quad \mS \in \cS.
    $$
    In particular, $\fM^\pi(f) = X^{\pi;\mP^f}$ for $f \in \cD$. By the uniqueness argument in Proposition~\ref{prop:relaxed_state_assignment_MM}, $\fM^\pi$ is the unique continuous orthogonal martingale measure with intensity $\eta^{\rho^\pi}$ for which $X^{\pi;\mS}$ solves (\ref{eq:SDE_martingale_measure}) for all $\mS \in \cS$. Moreover, for any finite family $\Pi = (\pi^i)_{i \in I} \subset \cR$, (\ref{eq:canonical_cylindrical_MM}) shows 
    $$
        \di [\fM^{\pi^i}(f^i), \fM^{\pi^j}(f^j)]_t = \bigg(\int_{\mA^I} f^i_t(a_i)^\top f^j_t(a_j) \, \lambda^\Pi_t(\di a) \bigg) \dint t,
    $$
    for $f^i \in L^2_{\pred,\loc}(\cH^{\rho^{\pi^i}}), f^j \in L^2_{\pred,\loc}(\cH^{\rho^{\pi^j}})$, and $\lambda^\Pi_t = (\Pi_t)_\# \dint r$; compare with~(\ref{eq:action_coupling_state_characteristics}). By~(\ref{eq:canonical_cylindrical_MM}), even if $\rho^\pi = \rho^{\pi'}$, $\fM^\pi$ and $\fM^{\pi'}$ coincide if and only if $\pi = \pi'$ $\oP \otimes \di r$-a.e.
\end{remark}

\section{Conclusion}
\label{sec:7}
We regard~(\ref{eq:SDE_cylindrical}) as the unique, up to choice of cylindrical Brownian driver, universally coherent relaxed state assignment (Theorems~\ref{theorem:cylindrical_approximation} and~\ref{theorem:cylindrical_converse}). Parametrized relaxed controls (paired with systems) are natural inputs (Theorem~\ref{theorem:ordinary_relaxed_assignment_common_label_dynamics} and Corollary~\ref{cor:complete_M_to_B}).

Given a relaxed control, any ground-stable limit of a sequence of strict control-state pairs has a finite-dimensional representation on a further Wiener extension (Theorem~\ref{theorem:residual_representation}).

The states~(\ref{eq:SDE_cylindrical}) can also be described by martingale measures, but this carries the particular choice~(\ref{eq:canonical_cylindrical_MM}) of those measures, determined by the parametrizations and common driver.

%%%%%%%%%%%%%%%%%%%%%%%%%%%%%%%%%%%%%%%%%%%%%%
%% Single Appendix:                         %%
%%%%%%%%%%%%%%%%%%%%%%%%%%%%%%%%%%%%%%%%%%%%%%
%\begin{appendix}
%\section*{???}%% if no title is needed, leave empty \section*{}.
%\end{appendix}
%%%%%%%%%%%%%%%%%%%%%%%%%%%%%%%%%%%%%%%%%%%%%%
%% Multiple Appendixes:                     %%
%%%%%%%%%%%%%%%%%%%%%%%%%%%%%%%%%%%%%%%%%%%%%%
%\begin{appendix}
%\section{???}
%
%\section{???}
%
%\end{appendix}

%%%%%%%%%%%%%%%%%%%%%%%%%%%%%%%%%%%%%%%%%%%%%%
%% Support information, if any,             %%
%% should be provided in the                %%
%% Acknowledgements section.                %%
%%%%%%%%%%%%%%%%%%%%%%%%%%%%%%%%%%%%%%%%%%%%%%
\begin{acks}[Acknowledgments]
I am grateful to Prof.~Dylan Possamaï and Prof.~Josef Teichmann for fruitful discussions.
\end{acks}
\bibliographystyle{spr-ims-number} % Style BST file (spr-ims-number.bst or spr-ims-nameyear.bst)
\bibliography{references}       % Bibliography file (usually '*.bib')

@article{JMLR:v21:19-144,
  author  = {Wang, Haoran and Zariphopoulou, Thaleia and Zhou, Xun Yu},
  title   = {Reinforcement learning in continuous time and space: A stochastic control approach},
  journal = {Journal of Machine Learning Research},
  year    = {2020},
  volume  = {21},
  number  = {198},
  pages   = {1--34},
  url     = {https://jmlr.org/papers/v21/19-144.html}
}

@article{ElKarouiMeleard1990,
  author  = {El Karoui, Nicole and M{\'e}l{\'e}ard, Sylvie},
  title   = {Martingale measures and stochastic calculus},
  journal = {Probability Theory and Related Fields},
  year    = {1990},
  volume  = {84},
  number  = {1},
  pages   = {83--101},
  doi     = {10.1007/BF01288560}
}

@article{AldousEagleson1978,
  author  = {Aldous, David J. and Eagleson, G. K.},
  title   = {On mixing and stability of limit theorems},
  journal = {The Annals of Probability},
  year    = {1978},
  volume  = {6},
  number  = {2},
  pages   = {325--331},
  doi     = {10.1214/aop/1176995577}
}

@book{DaPratoZabczyk2014,
  author    = {Da Prato, Giuseppe and Zabczyk, Jerzy},
  title     = {Stochastic equations in infinite dimensions},
  series    = {Encyclopedia of Mathematics and its Applications},
  volume    = {152},
  edition   = {2},
  publisher = {Cambridge University Press},
  year      = {2014},
  isbn      = {978-1-107-05584-1},
  doi       = {10.1017/CBO9781107295513}
}

@book{kallenberg2021foundations,
  author    = {Kallenberg, Olav},
  title     = {Foundations of modern probability},
  series    = {Probability Theory and Stochastic Modelling},
  volume    = {99},
  edition   = {3},
  publisher = {Springer},
  address   = {Cham},
  year      = {2021},
  isbn      = {978-3-030-61871-1},
  doi       = {10.1007/978-3-030-61871-1}
}

@incollection{JacodMemin1981,
  author    = {Jacod, Jean and M{\'e}min, Jean},
  title     = {Sur un type de convergence interm{\'e}diaire entre la convergence en loi et la
               convergence en probabilit{\'e}},
  booktitle = {S{\'e}minaire de Probabilit{\'e}s XV, 1979/80},
  editor    = {Az{\'e}ma, Jacques and Yor, Marc},
  series    = {Lecture Notes in Mathematics},
  volume    = {850},
  pages     = {529--546},
  publisher = {Springer},
  address   = {Berlin, Heidelberg},
  year      = {1981},
  doi       = {10.1007/BFb0088392},
  language  = {french},
  note      = {With a correction in S{\'e}minaire de probabilit{\'e}s XVII (1983), pp.~509--511}
}

@book{Takesaki2002Theory,
  author    = {Takesaki, Masamichi},
  title     = {Theory of operator algebras {I}},
  series    = {Encyclopaedia of Mathematical Sciences},
  volume    = {124},
  publisher = {Springer-Verlag},
  address   = {Berlin, Heidelberg},
  year      = {2002},
  isbn      = {978-3-540-42248-8},
  url       = {https://link.springer.com/book/9783540422488},
  note      = {Operator Algebras and Non-Commutative Geometry, Vol. 5; reprint of the first 1979
               edition}
}

@article{Yetter2005MeasurableCategories,
  author  = {Yetter, David N.},
  title   = {Measurable categories},
  journal = {Applied Categorical Structures},
  year    = {2005},
  volume  = {13},
  number  = {5--6},
  pages   = {469--500},
  doi     = {10.1007/s10485-005-9003-6}
}

@book{JacodShiryaev2003,
  author    = {Jacod, Jean and Shiryaev, Albert N.},
  title     = {Limit theorems for stochastic processes},
  series    = {Grundlehren der mathematischen Wissenschaften},
  volume    = {288},
  edition   = {2},
  publisher = {Springer-Verlag},
  address   = {Berlin, Heidelberg},
  year      = {2003},
  isbn      = {978-3-540-43932-5},
  doi       = {10.1007/978-3-662-05265-5}
}

@book{bogachev2007measure,
  author    = {Bogachev, Vladimir I.},
  title     = {Measure theory},
  publisher = {Springer-Verlag},
  address   = {Berlin, Heidelberg},
  year      = {2007},
  isbn      = {978-3-540-34513-8},
  doi       = {10.1007/978-3-540-34514-5},
  note      = {Volumes I--II}
}

@article{araki1981inequality,
  author  = {Araki, Huzihiro and Yamagami, Shigeru},
  title   = {An inequality for {Hilbert--Schmidt} norm},
  journal = {Communications in Mathematical Physics},
  year    = {1981},
  volume  = {81},
  number  = {1},
  pages   = {89--96},
  doi     = {10.1007/BF01941801}
}

@article{VeraarYaroslavtsev2016,
  author  = {Veraar, Mark and Yaroslavtsev, Ivan},
  title   = {Cylindrical continuous martingales and stochastic integration in infinite dimensions},
  journal = {Electronic Journal of Probability},
  year    = {2016},
  volume  = {21},
  number  = {59},
  pages   = {1--53},
  doi     = {10.1214/16-EJP7}
}

@article{Yaroslavtsev2018,
  author  = {Yaroslavtsev, Ivan},
  title   = {{Brownian} representations of cylindrical continuous local martingales},
  journal = {Infinite Dimensional Analysis, Quantum Probability and Related Topics},
  year    = {2018},
  volume  = {21},
  number  = {2},
  pages   = {1850013},
  doi     = {10.1142/S0219025718500133}
}

@article{Balder1989Prohorov,
  author  = {Balder, Erik J.},
  title   = {On {Prohorov}'s theorem for transition probabilities},
  journal = {Travaux du S{\'e}minaire d'analyse convexe},
  year    = {1989},
  volume  = {19},
  issn    = {1158-3029}
}

@article{ElKarouiHuuNguyenJeanblancPicque1987,
  author  = {El Karoui, N. and Huu Nguyen, D. and Jeanblanc-Picqu{\'e}, M.},
  title   = {Compactification methods in the control of degenerate diffusions: Existence of an
             optimal control},
  journal = {Stochastics},
  year    = {1987},
  volume  = {20},
  number  = {3},
  pages   = {169--219},
  doi     = {10.1080/17442508708833443}
}

@article{HaussmannLepeltier1990,
  author  = {Haussmann, U. G. and Lepeltier, J. P.},
  title   = {On the existence of optimal controls},
  journal = {SIAM Journal on Control and Optimization},
  year    = {1990},
  volume  = {28},
  number  = {4},
  pages   = {851--902},
  doi     = {10.1137/0328049}
}

@incollection{Walsh1986,
  author    = {Walsh, John B.},
  title     = {An introduction to stochastic partial differential equations},
  booktitle = {{\'E}cole d'{\'E}t{\'e} de Probabilit{\'e}s de {Saint Flour XIV} -- 1984},
  editor    = {Hennequin, P. L.},
  series    = {Lecture Notes in Mathematics},
  volume    = {1180},
  pages     = {265--439},
  publisher = {Springer-Verlag},
  address   = {Berlin, Heidelberg},
  year      = {1986},
  doi       = {10.1007/BFb0074920}
}

@misc{ElKarouiTan2024PartII,
  author        = {El Karoui, Nicole and Tan, Xiaolu},
  title         = {Capacities, measurable selection and dynamic programming part {II}: Application
                   in stochastic control problems},
  year          = {2024},
  eprint        = {1310.3364v3},
  archivePrefix = {arXiv},
  primaryClass  = {math.OC},
  doi           = {10.48550/arXiv.1310.3364},
  url           = {https://arxiv.org/abs/1310.3364v3},
  note          = {Originally posted in 2013; revised October 2, 2024}
}

@article{FontbonaGuerinMeleard2010,
  author  = {Fontbona, Joaquin and Gu{\'e}rin, H{\'e}l{\`e}ne and M{\'e}l{\'e}ard, Sylvie},
  title   = {Measurability of optimal transportation and strong coupling of martingale measures},
  journal = {Electronic Communications in Probability},
  year    = {2010},
  volume  = {15},
  pages   = {124--133},
  doi     = {10.1214/ECP.v15-1534}
}

@article{BenderThuan2026,
  author  = {Bender, Christian and Thuan, Nguyen Tran},
  title   = {Continuous time reinforcement learning: A random measure approach},
  journal = {Stochastic Processes and their Applications},
  year    = {2026},
  volume  = {194},
  pages   = {104848},
  doi     = {10.1016/j.spa.2025.104848}
}

@article{Balder2001ws,
  author  = {Balder, Erik J.},
  title   = {On {ws}-convergence of product measures},
  journal = {Mathematics of Operations Research},
  year    = {2001},
  volume  = {26},
  number  = {3},
  pages   = {494--518},
  doi     = {10.1287/moor.26.3.494.10581}
}

@article{VanNeervenVeraarWeis2007,
  author  = {van Neerven, Jan M. A. M. and Veraar, Mark C. and Weis, Lutz},
  title   = {Stochastic Integration in {UMD} Banach Spaces},
  journal = {The Annals of Probability},
  year    = {2007},
  volume  = {35},
  number  = {4},
  pages   = {1438--1478},
  doi     = {10.1214/009117906000001006}
}

@article{Warga1967,
  author  = {Warga, J.},
  title   = {Functions of Relaxed Controls},
  journal = {SIAM Journal on Control},
  year    = {1967},
  volume  = {5},
  number  = {4},
  pages   = {628--641},
  doi     = {10.1137/0305042}
}
%% or include bibliography directly:
% \begin{thebibliography}{}
% \bibitem{b1}
% \end{thebibliography}

\end{document}